\documentclass[12pt]{article}

\usepackage{amsmath,amsthm}
\usepackage{amssymb}
\usepackage{amsfonts}
\usepackage{amscd}

\newtheorem{thm}{Theorem}[section]
\newtheorem{theorem}[thm]{Theorem}

\newtheorem{corollary}[thm]{Corollary}
\newtheorem{lemma}[thm]{Lemma}

\newtheorem{proposition}[thm]{Proposition}

\theoremstyle{definition}
\newtheorem{definition}[thm]{Definition}

\newtheorem{remark}[thm]{Remark}
\newtheorem{observation}[thm]{Observation}

\newtheorem{Quantum line bundles}[thm]{Quantum line bundles}

\newcommand{\bk}{\Bbbk}

\newcommand{\C}{\mathbb{C}}
\newcommand{\Z}{\mathbb{Z}}
\newcommand{\N}{\mathbb{N}}

\newcommand{\bP}{\mathbb{P}}
\newcommand{\cO}{{\mathcal O}}
\newcommand{\cL}{{\mathcal L}}
\newcommand{\cV}{{\cal V}}
\newcommand{\bd}{{\overline{d}}}
\newcommand{\Ogh}{\cO\big(G\big/H\big)}
\newcommand{\Ogp}{\cO\big(G\big/P\big)}
\newcommand{\Oqgh}{\cO_q\big(G\big/H\big)}
\newcommand{\Oqgp}{\cO_q\big(G\big/P\big)}
\newcommand{\Oqghloc}{\cO_q^{\text{\it loc}}\big(G\big/H\big)}

\newcommand{\ep}{\epsilon}
\newcommand{\al}{\alpha}

\newcommand{\la}{\lambda}

\newcommand{\D}{\Delta}
\newcommand{\Dij}{\Delta_{ij}}
\newcommand{\tDelta}{\tilde{\Delta}}
\newcommand{\fh}{\mathfrak h}
\newcommand{\fg}{\mathfrak g}
\newcommand{\fp}{\mathfrak p}
\newcommand{\fgl}{\mathfrak{gl}}
\newcommand{\fsl}{\mathfrak{sl}}

\newcommand{\m}{\mathrm m}
\newcommand{\e}{\mathrm e}
\newcommand{\E}{\mathrm E}
\newcommand{\f}{\mathrm f}
\newcommand{\g}{\mathrm g}

\newcommand{\lra}{\longrightarrow}

\newcommand{\Ind}{{\hbox{Ind}}}

\begin{document}
{\ }
\vskip-2cm
\centerline{\LARGE \bf Quantization of Projective Homogeneous}
\bigskip
\centerline {\LARGE \bf Spaces and Duality Principle}

\vskip 1cm

\centerline{N. Ciccoli$^\natural$, R. Fioresi$^\flat$, F. Gavarini$^\#$
\footnote{N.~Ciccoli was partially supported by
TOK \emph{``Non commutative geometry and quantum groups''}, no.~MKTD-CT-2004-509794, University of Warsaw  ---  R.~Fioresi and F.~Gavarini were partially supported by the European RTN network ``LIEGRITS - Flags, Quivers and Invariant Theory in Lie Representation Theory'', no.~MRTN-CT 2003-505078, and by the Italian research program PRIN 2005 ``Moduli and Lie Theories''.  \\  
   \phantom{.}  \quad   {\it MSC 2000:} Primary 20G42, 14M17; Secondary 17B37, 81R50.  \\
   \phantom{.}  \quad   {\it Keywords:} Quantum Homogeneous Spaces.} }

\bigskip

 \centerline{\it $^\natural$ Dipartimento di Matematica,
              Universit\`a di Perugia }
\centerline{\it Via Vanvitelli 1, I-06123 Perugia, Italy}
\centerline{{\footnotesize e-mail: ciccoli@dipmat.unipg.it}}

\smallskip

\centerline{\it $^\flat$ Dipartimento di Matematica,
              Universit\`a di Bologna }
\centerline{\it Piazza di Porta San Donato, 5  ---
I-40127 Bologna, Italy}
\centerline{{\footnotesize e-mail: fioresi@dm.unibo.it}}

\smallskip

\centerline{\it ${}^\#$ Dipartimento di Matematica,
              Universit\`a di Roma ``Tor Vergata'' }
\centerline{\it Via della ricerca scientifica 1  --- I-00133 Roma,
Italy} \centerline{{\footnotesize e-mail:
gavarini@mat.uniroma2.it}}

\vskip1cm

\begin{abstract}
%
%
 We introduce a general recipe to construct quantum projective
homogeneous spaces, with a particular interest for the examples
of the quantum Grassmannians and the quantum generalized flag
varieties.  Using this construction, we extend the quantum
duality principle to quantum projective homogeneous spaces.
\end{abstract}

\tableofcontents

\section{Introduction}

A projective variety can be described via its homogeneous graded
coordinate ring. This ring is not an invariant associated to the
variety, but depends on a chosen embedding of the variety into
some projective space. Different embeddings will, in general,
produce non isomorphic graded rings.

\medskip

   When a projective variety is homogeneous, i.e. endowed with a
transitive action of an (affine) algebraic group on it, it can be
realized as quotient  of affine algebraic groups $ G\big/H \, $.
In this case  a projective embedding can be obtained via sections
of a line bundle on $ G\big/H $, uniquely given once  a character
of  $ H \, $ is specified.

\medskip

   If one approaches a quantization of this picture in the context of
   quantum groups the problem immediately arising is that
   standard quantum groups have a very limited set of quantum
   subgroups. This explains why usually the preferred approach
   goes through representation theoretic techniques.

\medskip

    An explanation of the lack of quantum subgroups, together with a way
   to circumvent this problem, is suggested by considering the
   semiclassical picture, i.e. in the context of algebraic Poisson
   groups. In such setting algebraic Poisson subgroups are quite
   rare too; however there is no need of an algebraic Poisson
   subgroup to cook up a Poisson quotient. The existence of a surjective Poisson map $G\to G\big/ H\, $
   is guaranteed  simply by requiring $H\,$ to be a coisotropic subgroup of
   $G\,$. This condition can be expressed by saying that the defining ideal of
   $ H \, $, in the function algebra of  $ G \, $,  is required to be a Poisson subalgebra
   rather than a Poisson ideal, as required for Poisson subgroups.

\medskip

   Let  $ \cO_q(G) $   be  a quantization of the affine
algebraic Poisson groups  $ G $. At the quantum level, a
quantization  $ \cO_q(H) $ of its {\sl coisotropic subgroup}  $ H
\, $ can be defined through conditions on the projection $ \, \pi
\colon \cO_q(G) \lra \cO_q(H) \, $. We will see this in full
detail in the sections \ref{class-set}, \ref{setting}.

\medskip

   Our first aim is to build a  quantum deformation  $ \Oqgh $  of the projective
variety  $ G \big/ H \, $,  i.e.~of its graded ring  $ \Ogh \, $,
subject to the following requirements:

\medskip

   {\it (1)} \,  there exists a one dimensional corepresentation of the quantum coiso\-tropic subgroup  $ \cO_q(H) $  which is a deformation of the  corepresentation of  $ \cO(H) $  corresponding to the character of
$ H $  which defines the line bundle giving the projective
embedding of  $ G\big/H \, $;

\smallskip

   {\it (2)} \,  a quantum analogue  $ \Oqgh $  to  $ \Ogh $  is defined as the subset   --- inside  $ \cO_q(G) $  ---   of ``semi-invariant functions'' with respect to the given corepresentation of  $ \cO_q(H) \, $;

\smallskip

   {\it (3)} \,  the subset  $ \Oqgh $  is a graded subalgebra of
$ \cO_q(G) \, $;

\smallskip

   {\it (4)} \,  the graded subalgebra  $ \Oqgh $  is a graded left coideal of  $ \cO_q(G) \, $,  so the coproduct in  $ \cO_q(G) $  induces a (left)
$ \cO_q(G) $--coaction  on  $ \Oqgh \, $,  and the latter can be
thought of as a  {\sl quantum homogeneous space}.

\smallskip

   {\it (5)} \,  the semiclassical limit of  $ \Oqgh $  is  $ \Ogh $   --- embedded into  $ \cO(G) $  ---   as a graded subalgebra, left coideal and graded Poisson subalgebra.

\medskip

   In other words, a quantum deformation of a projective homogeneous
space, embedded into some projective space, consists of the
deformation of the graded algebra associated to the embedding, in
such a way that the action of the group on the homogeneous space
is also naturally quantized.
                                           \par

\medskip

   We will work out the details of the construction for the case
of the Grassmannian and its Pl{\"u}cker embedding, that is when $
G $  is the special linear group and  $ \, H = P \, $  is a
maximal parabolic subgroup and we will sketch it in the more
general case of quantum flag varieties of simple Lie groups.

\medskip

   Our main motivation to develop this point of view is to adapt to
projective homogeneous spaces, the correspondence introduced by
Ciccoli and Gavarini  \cite{cg1} for coinvariant subalgebras. This
recipe allows to associate functorially to a quantum quasi-affine
homogeneous space another quantum homogeneous space, through a
generalization of the \textit{quantum duality principle} (QDP),
defined by Drinfeld for quantum groups. A part of the arguments in
\cite{cg1} does not directly apply to  {\sl projective}
homogeneous spaces, since it is based on the realization of the
ring of the homogeneous space as the set of coinvariant functions
inside the  ring of the quantum group acting on it. But this is
possible   --- as in the classical case --- if and only if the
homogeneous space is quasi-affine, which is not the case of
projective varieties.  The coordinate ring of the homogeneous
space is replaced by a graded ring inside the quantum group ring,
consisting of semi-coinvariants with respect to a one-dimensional
representation, which can be seen as a deformation of the line
bundle that classically determines the projective embedding. The
definitions introduced in section 3 will allow us to define a
quantum duality functor and obtain the QDP construction in this
more general setting. In the last chapter we will discuss
applications to quantum flag manifolds.


\section{The classical setting} \label{class-set}

   In this section we recall some Poisson geometry
(see  Ref.~\cite{kos}  for details).

\subsection{The affine case}

   Let  $ \bk $  be a fixed field of characteristic 0. When doing algebraic
geometry construction, we tacitly assume that  $ \bk $  be algebraically closed.
However, this assumption is  {\sl not\/}  needed for our quantum constructions.
                                            \par
   Let $G$ be an affine algebraic group over  $ \bk \, $.  We denote by
$ \cO(G) $  the algebra of regular functions of  $ G $,  in short its
``function algebra'', which is naturally a Hopf algebra (over  $ \bk \, $).
We denote  $ \, \fg := \text{\it Lie}(G) \, $  the tangent Lie algebra
of  $ G \, $,  and  similarly  $ \, \fh := \text{\it Lie}(H) \, $  for
any closed algebraic subgroup  $ H $  of  $ G \, $.

\medskip

  Assume that  $ G $  is a Poisson group: this means that  $ \cO(G) $  is
a  {\sl Poisson\/}  Hopf algebra, i.e.{} we have a Poisson bracket  $ \,
\{\ ,\ \} \colon \cO(G) \otimes \cO(G) \longrightarrow \cO(G) \, $  which
is compatible with the Hopf algebra structure.  Moreover,  $ \fg $  is a
Lie bialgebra, for some Lie cobracket  $ \, \delta \colon \fg \longrightarrow
\fg \otimes \fg \, $,  and the same holds for its dual space  $ \fg^* \, $,
these two Lie bialgebras structures being dual to each other.  Indeed, the
notion of Lie bialgebra is the infinitesimal counterpart of the notion of
Poisson group.  Since the dual $\fg^*$ of $\fg$ is itself a Lie bialgebra,
it follows that any connected algebraic group $G^*$ with $ \fg^*=Lie(G^*) $,
is a Poisson group on its own, called  {\sl (Poisson) dual\/}  to  $ G $.
We are going to see the example of the Poisson group $GL_n$ treated
in detail in section \ref{qdp-grass}.
                                         \par
\begin{definition}  \label{coisotropic}
   A (closed) subgroup  $ H $  of  $ G $  is called
{\sl coisotropic\/}  if its defining ideal  $ I(H) $  is a Poisson subalgebra
of  $ \cO(G) \, $.  Also,  $ H $  is called
a {\sl Poisson subgroup\/}  if the embedding
$ \, H \lhook\joinrel\relbar\joinrel\rightarrow G \, $  is a
Poisson map; this is equivalent to require  $ I(H) $  to be
a Poisson  {\sl ideal}. Hence a Poisson subgroup is coisotropic.
\end{definition}
                                           \par
  The following equivalent conditions give an infinitesimal characterization
for a  {\sl connected\/}  subgroup  $ H $  to be coisotropic (see  \cite{lu}):

\begin{proposition}  \label{coiso-inf}
 Let  $ G $  be an algebraic group and  $ H $  a (closed) subgroup
of the Poisson group  $ G \, $.  Then the following are equivalent:
 \vskip2pt
\noindent
 {\it (C-i)} \,  $ H $  is a coisotropic subgroup of
$ \, G \, $;
                                                 \par
\noindent
 {\it (C-ii)} \,  $ \delta(\fh) \subseteq \fh \wedge \fg
\, $,  that is  $ \fh $  is (a Lie subalgebra and) a Lie
coideal of  $ \, \fg \, $;
                                                 \par
\noindent
 {\it (C-iii)}  the orthogonal space  $ \, \fh^\perp $  is
(a Lie coideal and) a Lie subalgebra of  $ \, \fg^* \, $.
\end{proposition}

\medskip

\begin{remark}
  Note that, thanks to these characterizations, the infinitesimal
counterpart of the notion of coisotropic subgroup is that of Lie subalgebra
Lie coideal.  The latter notion is self-dual.
In fact, let $G^*$ be any connected Poisson
group dual to $G$.
If  $ H $  is coisotropic in  $ G $,
then any connected subgroup of  $ G^* $,  with tangent Lie algebra
$ \fh^\perp $,  is in turn a coisotropic subgroup of  $ G^* $,  called
``complementary dual'' to  $ H $  and denoted by  $ H^\perp $.
\end{remark}

   We now want to describe the notion of {\sl Poisson quotient}.

\medskip

\begin{definition}
   Let $ M $  be a Poisson affine variety, i.e.~an affine
variety whose function algebra  $ \cO(M) $  is a Poisson algebra.
Then  $ M $  is a  {\sl Poisson homogeneous  $ G $--space\/}  if
there is a (regular) transitive action  $ \, \phi \, \colon \, G \times M
\rightarrow M \, $  which is a Poisson map with respect to the product
Poisson structure on  $ G \times M \, $.  We say that a Poisson homogeneous
$ G $--space  $ M $  is a  {\sl Poisson quotient\/}  if there is a coisotropic
closed Lie subgroup  $ H_M $  of  $ G $  such that  $ \, G \big/ H_M \simeq M
\, $  and the projection  $ \, p_M \colon \, G \relbar\joinrel\twoheadrightarrow
G \big/ H_M \simeq M \, $  is a Poisson map.
\end{definition}

\medskip

   The following is a characterization of Poisson quotients  (cf.~\cite{za}).

\medskip

\begin{proposition}
   Let  $ \; \phi : G \times M \lra M \; $  be a homogeneous action of  $ \, G $
on  $ M \, $.  Then the following are equivalent:
 \vskip1pt
  \; {\it (PQ-i)} \, there exists  $ \, \bar{m} \in M \, $  whose
stabilizer  $ \, G_{\bar{m}} \, $  is a coisotropic subgroup of
$ \, G \, $;
                                                 \par
\noindent
  \; {\it (PQ-ii)} \, there exists  $ \, \bar{m} \in M \, $  such that
$ \, \phi_{\bar{m}} \, \colon \, G \longrightarrow M \, $,  $ \, g \mapsto \phi(g,\bar{m}) \, $,  \; is a Poisson map, that is to say
$ M $  is a Poisson quotient;
                                                 \par
\noindent
  \; {\it (PQ-iii)} \, there is  $ \, \bar{m} \in M \, $  such that
$ \, \{\bar{m}\} \, $  is a symplectic leaf of  $ M \, $.
\end{proposition}

\medskip

   For any  $ \, \bar{m} \in M \, $  with stabilizer  $ G_{\bar{m}} $
one has  $ \, M \sim G \big/ G_{\bar{m}} \, $  as affine  $ G $--varieties.
As  $ M $  is affine, this is equivalent to  $ \, \cO(M) \cong \cO\big(G/G_{\bar{m}}\big) \, $.  Finally,  $ \, \cO\big(G\big/G_{\bar{m}}\big) \cong {\cO(G)}^{G_{\bar{m}}} \, $,  the subalgebra of  $ G_{\bar{m}} $--invariants  in  $ \cO(G) $.  The same holds with  $ G_{\bar{m}} $  replaced by any subgroup  $ H $  whose coset space $ G\big/H $
%
%
 is affine.  We have then an additional characterization of (affine)
Poisson quotients.

\medskip

\begin{proposition}
   If  $ M $  is as above, then  $ M $  is a Poisson quotient if and only if there
exists  $ \, \bar{m} \in M \, $  such that  $ {\cO(G)}^{G_{\bar{m}}} $  is a Poisson subalgebra of  $ \cO(G) \, $.
\end{proposition}

\medskip

   In particular, if  $ H $  is a subgroup of  $ G $,
and $G/H$ is affine then the following are equivalent:
\\
1. $H$ is coisotropic,  \\
2. $G\big/H$ is a Poisson quotient, \\
3. $\Ogh = {\cO(G)}^H$  is a Poisson subalgebra of $\cO(G)$.


\medskip

\subsection{The projective case}  \label{proj-case}
We are now interested in the case when
the homogeneous  $ G $--variety  $ G\big/H $  is  {\sl projective},
i.e.~$ H $  is  {\sl parabolic}.  To describe it in algebraic terms
(the setting we need for quantum deformations), we require a specific
realization, namely an embedding into a projective space.
                                            \par
  Given a representation  $\rho$ of $H$ on some vector space $V$,  we can
construct a vector bundle associated to it, namely
  $$  \cV := G \times_H V = G \times V/ \simeq \, , \quad
(gh,v) \simeq (g, h^{-1}v) \, , \qquad \forall h \in H,
g \in G, v \in V .  $$
The space of global sections of this bundle is identified with the induced
module (see, e.g., \cite{ha} for more details)
  $$  H^0 \big( G\big/H, \cV \big)  \, = \,  \Ind^G_H(V)  \, = \,
\big\{ f \colon G \rightarrow V \,\big|\, f \text{\ is regular},
\, f(gh) = h^{-1}.f(g) \,\big\} \; .  $$

\medskip

\begin{definition}
   Let  $ \, \chi \colon \, H \lra \Bbbk^* \, $  be a character of  $ H $,
i.e.~a one dimensional representation of  $ H $  on  $ \, L \cong \Bbbk \, $.
Then  $ L^{\otimes n} $  is again a one dimensional representation of  $ H $
with character  $ \chi^n \, $.  Let  $ \, \cL^n := G \otimes_H L^{\otimes n} \, $.
Define
  $$  \displaylines{
   \Ogh_n  \; := \; H^0 \big( G\big/H, \cL^{n} \big)  \cr \cr
   \Ogh  \; := \;  {\textstyle \bigoplus_{n \geq 0}}
\, \Ogh_n \; \subseteq \;  \cO(G)  \quad . }  $$
 Then  $ \Ogh $  is a subalgebra of  $ \cO(G) \, $,  whose elements are called
{\it semi-invariants\/}. Note that now the notation
$ \Ogh $  has not the same
meaning as when  $ G\big/H $  is affine.
\end{definition}

\smallskip

   Assume now the bundle  $ \cL $  to be very ample.  In the present context,
this is the same (cf.~\cite{ha}, \S II.7) as saying that  $ \cL $  is generated
by a set of global sections  $ f_0 $,  $ f_1 $,  $ \dots $,  $ f_N \in \Ogh_1 \, $
--- in particular, the algebra  $ \Ogh $  {\sl is graded, generated in degree 1\/}
(by the  $ f_i $'s).  Then  $ \Ogh $  is the homogeneous coordinate ring of the
projective variety  $ G\big/H \, $,  with respect to the embedding  given via
the global sections of  $ \cL $  (see \cite{gh}, p.~176).

\smallskip

   We want to reformulate this classical construction in purely Hopf algebraic
terms, more suited to the quantum setting we shall presently deal with.
%

\vskip-1pt

\begin{remark}  \label{desc-semiinv}
 In algebraic terms, saying that  $ \, \chi \, \colon H \rightarrow \Bbbk^* \, $
is a character is the same as saying that it is a group-like element in the coalgebra
$ \cO(H) \, $.  The same holds for all powers  $ \chi^n $  ($ n \! \in \! \N
\, $). In fact
if  $ \chi $  is group-like,
then the same is true for all its powers  $ \chi^n $
since  $ \cO(H) $  is a Hopf algebra.  As the  $ \chi^n $'s  are
group-like, if they are pairwise different they also are linearly independent,
which ensures that the sum  $ \, \sum\limits_{n \in \N} \Ogh_n \, $   --- inside
$ \cO(G) $  ---   is a direct one.  Moreover, once the embedding is given, each
summand  $ \Ogh_n $  can be described in purely Hopf algebraic terms as
  $$  \displaylines{
   \Ogh_n  \, := \,  \big\{ f \in \cO(G) \,\big|\, f(gh) =
\chi^n\big(h^{-1}\big) f(g) \,\big\} \, = \,  \cr
   \hfill   = \, \Big\{\, f \in \cO(G) \,\Big|\, \big((\text{\it id}
\otimes \pi) \circ \Delta \big)(f) = f \otimes S\big(\chi^n\big) \Big\}
\qquad }  $$
 with  $ \, \pi \, \colon \cO(G) \relbar\joinrel\twoheadrightarrow \cO(H)
\, $  the standard projection,  $ S $  the antipode of  $ \cO(H) \, $.

\vskip4pt

   To simplify notation, we set  $ \, \lambda := S(\chi) \, $   --- the
character of  $ H $  which maps  $ \, h \in H \, $  to  $ \, \lambda(h)
= \chi\big(h^{-1}\big) \, $  ---  and we set  $ \; \Delta_\pi :=
(\text{\it id} \otimes \pi) \circ \Delta \; $,  \, so that
  $$  \Ogh_n  \, = \;  \Big\{\, f \in \cO(G) \;\Big|\,
\Delta_\pi(f) = f \otimes \lambda^n \,\Big\}  \quad .   \eqno (2.1)  $$
\end{remark}

%
%
%

\begin{proposition} \label{t}
Let $G/H$ be embedded into some projective space via some
very ample line bundle. Then there exists a $t \in \cO(G)$ such
that
  $$  \Delta_\pi(t) \, := \, \big((\text{\it id} \otimes \pi)
\circ \Delta \big)(t) \, = \, t \otimes \pi(t)   \eqno (2.2)  $$
  $$  \pi\big(t^m\big) \not= \pi\big(t^n\big) \quad \forall \;\; m \not= n \in \N   \eqno (2.3)  $$
  $$  \Ogh_n  \; = \;  \Big\{\, f \in \cO(G) \,\;\Big|\;
\Delta_\pi(f) = f \otimes \pi\big(t^n\big) \Big\}   \eqno (2.4)  $$
  $$  \Ogh  \; = \;  {\textstyle \bigoplus_{n \in \N}} \; \Ogh_n   \eqno (2.5)  $$
and $ \, \Ogh \, $  is generated in degree 1, namely by
$ \Ogh_1 \, $.
\end{proposition}

\begin{proof}
   If  $ \, f \in \Ogh_n \, $,  \, then
  $$  \displaylines{
   \pi(f)  \, = \,  \pi \Big( {\textstyle \sum_{(f)}}
\, \epsilon\big(f_{(1)}\big) \, f_{(2)} \Big)  \, = \,
(\epsilon \otimes \pi) \big( \Delta(f) \big)  \, = \,
(\epsilon \otimes \text{\it id}\,) \big( \Delta_\pi(f) \big)
\, = \,  \epsilon(f) \, \lambda^n  }  $$
%
%
Now, by assumption there exists a non-zero global section of the line
bundle on  $ G\big/H \, $,  i.e.~a regular function  $ \, t \in \Ogh_1
\setminus \{0\} \, $  on  $ G \, $
and $\ep(t) \neq 0$.  By the above (for  $ \, n = 1 \, $),
up to dividing out by  $ \epsilon(t) $,  we can assume that  $ \, \pi(t)
= \lambda \, $. The result follows immediately.

\vskip2pt



%

\end{proof}

\smallskip

   Notice that while  $ \, \lambda = \pi(t) \, $  {\sl is\/}  group-like,
$ t $  instead is something less, yet still has an ``almost
group-like property'', given by (2.2). This element $t$ and
its quantization will turn to be crucial for the quantum setting.

\medskip

\begin{remark}  \label{rem-Ogh}
 We point out that  $ \Ogh $  is a unital subalgebra,  {\sl and\/}  also
a (left) coideal of  $ \cO(G) \, $;  the latter reflects the fact that
$ G\big/H $  is a (left)  $ G $--space.  Thus, the restriction of the
comultiplication of  $ \cO(G) \, $,  namely
  $$  \Delta\big|_{\cO(G/H)} : \Ogh \lra \cO(G) \otimes \Ogh \quad ,  $$
is a coaction of  $\cO(G)$ on  $ \Ogh $,  which makes  $\Ogh $
into an  $ \cO(G) $--comodule  algebra, in the sense of  \cite{mo},
\S 4.1.  Moreover $\Ogh$ is {\sl graded\/}  and the coaction
$ \Delta\big|_{\cO(G/H)} $  is also {\sl graded\/}  w.r.t.~the
trivial grading on  $ \cO(G) \, $,  so that each  $ \Ogh_n $  is
indeed a coideal of  $ \cO(G) $  as well.
\end{remark}

\medskip

   As to Poisson geometrical properties in this projective setup, the
following characterization   --- which might be used to define the notion
of  {\sl Poisson quotient\/}  structure for the projective  $ G $--space
$ G\big/H $  ---   holds:

\vskip5pt

\begin{proposition}  \label{Pq-proj}
   Let  $ G $  be a Poisson algebraic group,  $ H $  a closed parabolic subgroup, and  $ \, t \in \cO(G) \, $  as in  Remark \ref{desc-semiinv}.
The following are equivalent:
 \vskip4pt
   (a) \quad  $ \big\{ I(H) \, , I(H) \big\} \subseteq I(H) \; $  --- that is,  $ H $  is coisotropic ---   {\sl and}  in addition  $ \,\; \phantom{\Big|} \big\{ t \, , \Ogh \big\} \subseteq I(H)
\; $;
 \vskip4pt
   (b) \quad  $ \big\{\, \Ogh_r \, , \, \Ogh_s \,\big\} \subseteq \Ogh_{r+s} \;\; $  for
all  $ r, s \in \N \, $   --- that is,  $ \, \Ogh $  is a  {\sl graded Poisson}  subalgebra
of  $ \, \cO(G) \; $.
\end{proposition}

\noindent
 {\it Proof.}  To simplify notation, we set  $ \; I_{H,n} \! := I(H) \bigcap
\Ogh_{\,n} \, $  for  $ \, n \in \N \, $.
 \vskip5pt
   {\it (a)}  $ \Longrightarrow \! $  {\it (b)\,:} \,  First of all, note that
(2.1) can be reformulated as
  $$  \Ogh_{\,n}  \; = \;  \Big\{\, f \in \cO(G) \,\;\Big|\; \Delta(f)
\, \in \, f \otimes t^n + \cO(G) \otimes I(H) \Big\}   \eqno (2.6)  $$
   \indent   Second, by Remark  \ref{rem-Ogh}{\it (a)},  each  $ \Ogh_{\,n} $
is a coideal of  $ \cO(G) \, $,  that is  $ \; \Delta \Big( \Ogh_{\,n} \Big)
\subseteq \cO(G) \otimes \Ogh_{\,n} \, $.  This along with (2.6) gives
  $$  \Delta(f) \, \in \, f \otimes t^n + \, \cO(G) \otimes I_{H,n}
\quad \qquad \forall \;\, f \in \Ogh_{\,n} \; .   \eqno (2.7)  $$
   \indent   Then, for any  $ \, f \in \Ogh_{\,r} \, $  and  $ \, \ell \in
\Ogh_{\,s} \, $,  \, we have
  $$  \displaylines{
   \Delta \big( \{f,\ell\} \big) \, = \, \big\{ \Delta(f) \, , \, \Delta(\ell) \big\}
\, \in \, \Big\{ f \otimes t^r + \, \cO(G) \otimes I_{H,r} \, , \, \ell \otimes t^s +
\, \cO(G) \otimes I_{H,s} \Big\} \, =  \cr
   = \; \big\{ f \otimes t^r \, , \, \ell \otimes t^s \,\big\} \, + \,
\Big\{ f \otimes t^r \, , \, \cO(G) \otimes I_{H,s} \Big\} \;\, +   \hfill  \cr
   + \,\; \Big\{ \cO(G) \otimes I_{H,r} \, , \, \ell \otimes t^s \Big\} \, +
\, \Big\{ \cO(G) \otimes I_{H,r} \, , \, \cO(G) \otimes I_{H,s} \Big\}
\;\, \subseteq  \cr
   \subseteq \,\; \big\{ f, \ell \big\} \otimes t^{r+s} \, + \, f \, \ell \otimes \big\{ t^r, t^s \big\} \, + \, \Big\{ f \, , \, \cO(G) \Big\} \otimes t^r \, I_{H,s} \;\, +
\hfill  \cr
   + \,\; f \, \cO(G) \otimes \Big\{ t^r \, , I_{H,s} \Big\} \, + \, \Big\{ \cO(G) , \, \ell \,\Big\} \otimes I_{H,r} \, t^s \, + \, \cO(G) \ell \, \otimes \, \Big\{ I_{H,r} , \, t^s \Big\} \;\, +  \cr
   \hfill   + \,\; \Big\{ \cO(G) , \, \cO(G) \Big\} \otimes I_{H,r} \, I_{H,s} \, + \, \cO(G) \, \cO(G) \otimes \Big\{ I_{H,r} \, , \, I_{H,s} \Big\} \; \subseteq  \cr
   \hfill   \subseteq \;  \big\{ f, \ell \big\} \otimes t^{r+s} \, + \, \cO(G) \otimes I_{H,r+s} }  $$
thanks to (2.7) and to  {\it (a)}.  Thus  $ \; \{f,\ell\} \in \Ogh_{r+s} \; $,  \, by (2.6) again.

 \vskip5pt

   {\it (b)}  $ \Longrightarrow \! $  {\it (a)\,:} \,  By assumption we havewe have  $ \, t \in {\mathcal{O}\big(G/H\big)}_1 \, $,  hence  $ \, \{t,f\} \in {\mathcal{O}\big(G/H\big)}_{1+n} \, $  for all  $ \, f \in {\mathcal{O}\big(G/H\big)}_n \, $  by  {\it (b)}.  In particular, this gives
  $$  \{t,f\}(h)  \; = \;  t^{1+n}(h) \, \{t,f\}(1_{\scriptscriptstyle G})  \; = \;  0
  \eqno \forall \;\; h \in H  \quad  $$
because any Poisson group structure is zero at the identity.  Eventually, this yields
$ \; \big\{ t \, , \Ogh \big\} \subseteq I(H) \, $,  \; q.e.d.
 \vskip5pt
   To prove that  $ \, \big\{ I(H) \, , I(H) \big\} \subseteq I(H) \, $,  \, we need some additional tools.
                                                        \par
   First, let  $ \cO(G)_{\,1_G} $  be the localization of  $ \cO(G) $  at  $ \, J := \text{\it Ker}\,(\epsilon_{\cO(G)}) \, $   --- a maximal ideal in  $ \cO(G) \, $.  This
is the stalk at the point  $ 1_G $  of the structure sheaf of  $ G \, $,  and the Poisson bracket of  $ \cO(G) $  canonically (and uniquely) extends   --- for  $ \, f_1, f_2 \in \cO(G) \, $,  $ \, y_1, y_2 \in J \, $,  $ n_1, n_2 \in \N \, $  ---   via the identity
  $$  \displaylines{
   \;\;   \big\{ f_1 \, y_1^{-n_1} , f_2 \, y_2^{-n_2} \big\}  \; = \;  \{f_1,f_2\} \, y_1^{-n_1} y_2^{-n_2}  \, - \,  n_1 \, f_1 \, y_1^{-n_1-1} \, \big\{ y_1 ,  f_2 \big\}
\, y_2^{-{n_2}}  \;\, -
\hfill  \cr
   \hfill   - \,\;  n_2 \, f_2 \, \big\{ f_1, f_2 \big\} \, y_1^{-n_1} \, y_2^{-n_2-1}
\, + \,  n_1 \, n_2 \, f_1 \, f_2 \, \big\{ y_1, y_2 \big\} \, y_1^{n_1-1} \, y_2^{n_2-1}   \;\; }  $$
to the local algebra  $ \cO(G)_{\,1_G} \, $.  The co\-unit  $ \epsilon $  of  $ \cO(G) $  uniquely extends to an algebra morphism from  $ \cO(G)_{\,1_G} $  to  $ \Bbbk \, $,  again denoted by  $ \epsilon \, $,  whose kernel is  $ J_{1_G} \, $,  the localization of  $ J $  inside  $ \cO(G)_{\,1_G} \, $.  Finally, we denote by  $ I(H)_{1_G} $  the localization,
inside  $ \cO(G)_{\,1_G} \, $,  of the ideal  $ I(H) $  of  $ \cO(G) \, $.
                                                        \par
   Second, let  $ \, X_t := \big\{\, t = 0 \,\big\} \, $  be the zero locus in  $ G\big/H $  defined by the vanishing of the divisor  $ t \, $,  and let  $ \; \varGamma := \big( G\big/H \big) \setminus X_t \; $.  This is an affine open dense subset of  $ G\big/H \, $,  whose algebra of regular functions is the graded localisation of  $ \Ogh $  by the multiplicative subset  $ \big\{ t^n \big\}_{n \in \N} \, $,  that is  $ \; \Ogh_{[t]} := \bigoplus_{n \in \N} t^{-n} \, \Ogh_{\,n} \; $.  Note that  $ \Ogh_{[t]} $  naturally embeds into
$ \cO(G)_{\,1_G} \, $,  because  $ \, t \in J^{G/H} := \Ogh \setminus J \, $.  Again,
the Poisson bracket of  $ \Ogh $,  induced by that of  $ \cO(G) \, $,  uniquely extends
to  $ \Ogh_{[t]} \, $,  and so the latter is a  {\sl graded\/}  Poisson subalgebra of
$ \cO(G)_{\,1_G} \, $;  thus  $ X_t $  is an affine Poisson variety.
%
%
 Also,  $ \epsilon $  induces an algebra morphism from  $ \Ogh_{[t]} $  to  $ \Bbbk
\, $,  whose kernel we denote by  $ J^{G/H}_{[t]} \, $.
                                                        \par
   Third, let  $ \Ogh_{\,\overline{1_G}} $  be the localization of  $ \Ogh_{[t]} $  at
$ J^{G/H}_{[t]} \, $:  by construction, this is the stalk at  $ \, \overline{1_G} = 1_G
\, H \, $  of the structure sheaf of  $ G\big/H \, $,  and the Poisson bracket of
$ \Ogh_{[t]} $  uniquely extends  to $ \Ogh_{\,\overline{1_G}} \, $.
%
%
                                                        \par
   Now, the maximal ideal  $ J_{1_G} $  in the local algebra  $ \cO(G)_{\,1_G} $  can be generated by a local system of parameters on  $ G $  at the point  $ 1_G \, $,  say  $ \, \{y_1, \dots, y_n\} \, $,  with  $ \, n := \text{\it dim}\,(G) \, $.  As  $ H $  is a closed subgroup of  $ G \, $,  we can choose this system of parameters in such a way that:
                                                        \par
   {\it (1)} \,  if  $ \, h := \text{\it dim}\,(H) \, $,  the image inside  $ \, \cO(H)_{1_H} \cong \cO(G)_{1_G} \big/ I(H)_{1_G} \, $  of  $ \, \{y_1, \dots, y_h \} \, $  is a local system of parameters on  $ H $  at the point  $ \, 1_H = 1_G \; $;
                                                        \par
   {\it (2)} \,  $ y_h = t-1 \; $;
                                                        \par
   {\it (3)} \,  $ \{y_{h+1}, \dots, y_n \} \, $  is a local system of parameters on
$ X_t $  at  $ 1_G \, $;  in particular, it generates in  $ \, \Ogh_{\,\overline{1_G}}
\, $  the ideal  $ \; J^{G/H}_{\,\overline{1_G}} := \text{\it Ker}\,(\epsilon) \bigcap \Ogh_{\,\overline{1_G}} \; $.
 \vskip4pt
   As a direct consequence of the above assumptions, the elements  $ y_{h+1} $,  $ \dots $,
$ y_n $  generate the ideal  $ I(H)_{1_G} $  inside  $ \cO(G)_{\,1_G} \, $.  Moreover, we have  $ \; y_i = f_i \, \ell_i^{-1} \; $  for some  $ \, f_i \in \Ogh \, $,  $ \, \ell_i \in \Ogh \setminus J^{G/H} \, $  (for  $ \, i = h+1, \dots, n \, $).  Then
  $$  \displaylines{
   \big\{ y_i , y_j \big\}  \; = \;  \big\{ f_i \, \ell_i^{-1} , f_j \, \ell_j^{-1} \big\}   \; = \;  \{f_i,f_j\} \, \ell_i^{-1} \, \ell_j^{-1}  \, -  \big\{f_i,\ell_j\big\} \,
\ell_i^{-1} \, f_j \, \ell_j^{-2}  \, - \,   \hfill  \cr
   \hfill   - \,  f_i \, \big\{\ell_i,f_j\big\} \, \ell_i^{-2} \, \ell_j^{-1}  \, + \,
f_i \, f_j \, \big\{\ell_i^{-1},\ell_j^{-1}\big\} \, \ell_i^{-2} \, \ell_j^{-2} }  $$
which   --- by assumption  {\it (b)},  that yields  $ \, \big\{ \Ogh, \Ogh \big\} \subseteq \Ogh \, $  ---   gives  $ \, \big\{ y_i , y_j \big\} \in \big( y_{h+1}, \dots, y_n \big) = I(H)_{1_G} \, $,  \, for all  $ i $,  $ j \, $.  This together with Leibnitz' rule implies that  $ \, \big\{ k_i \, y_i \, , \, k_j \, y_j \big\} \in \big( y_{h+1}, \dots, y_n \big) = I(H)_{1_G} \, $  for any  $ \, k_i, k_j \in \cO(G)_{\,1_G} \, $  (with  $ i $,  $ j = h+1, \dots, n \, $);  in turn,  $ \, I(H)_{1_G} = \big( y_{h+1}, \dots, y_n \big) \, $  satisfies  $ \, \big\{ I(H)_{1_G} \, , \, I(H)_{1_G} \big\} \subseteq I(H)_{1_G} \, $.
                                                        \par
   Eventually, as  $ \, I(H) = \cO(G) \bigcap I(H)_{1_G} \, $,  the above results give also
  $$  \big\{ I(H) , I(H) \big\}  \, \subseteq \,  \cO(G) \, {\textstyle \bigcap} \, \big\{ I(H)_{1_G} , I(H)_{1_G} \big\}  \, \subseteq  \cO(G) \, {\textstyle \bigcap} \, I(H)_{1_G}  = \,  I(H)   \qed  $$

\bigskip

\section{Quantum bundles and quantum homogeneous spaces}
\label{setting}

\subsection{Quantum groups}
We want to translate all
the framework of section  \ref{class-set}  into the quantum setup.
The first step is to introduce quantum groups, in the form of
quantum (or ``quantized'') function algebras, as follows.

\medskip

Let $G$ be an algebraic Poisson group, $\cO(G)$ its function algebra.

\begin{definition}  \label{quant-G}
   By  {\sl quantization\/}  of  $ \cO(G) $,  we mean a Hopf algebra
$ \cO_q(G) $  over the ground ring  $ \, \Bbbk_q := \Bbbk\big[q,q^{-1}\big]
\, $, where $q$ is an indeterminate, such that:
                                      \par
\hskip-9pt
   {\it (a)} \,  the  {\sl specialization\/}  of  $ \cO_q(G) $  at
$ \, q = 1 \, $,  that is  $ \, \cO_q(G) \Big/ (q\!-\!1) \, \cO_q(G)
\, $,  is isomorphic to  $ \cO(G) $  as a Poisson Hopf algebra;
                                       \par
\hskip-9pt
   {\it (b)} \,  $ \cO_q(G) $  is torsion-free, as a  $ \Bbbk_q $--module;
                                       \par
\hskip-9pt
   {\it (c)} \,  if  $ \, I_G := (q-1) \, \cO_q(G) + \text{\it Ker}\,
\big(\epsilon_{\cO_q(G)}\big) \, $,  then  $ \, \bigcap\limits_{n \geq 0}
\! I_G^{\,n} = \! \bigcap\limits_{n \geq 0} {(q-1)}^n \, \cO_q(G) \, $.
 \medskip
   We call  $ \cO_q(G) $  {\sl quantum\/  {\rm (or}
quantized\/{\rm)}  function algebra\/}  over  $ G \, $,
or  {\sl quantum deformation} of  $ G \, $,  or even simply
{\sl quantum group}.  It is standard terminology to say that the  {\sl Poisson\/}  Hopf algebra  $ \cO(G) $  is the  {\sl semiclassical limit\/}  of  $ \cO_q(G) \, $.

\smallskip

   Similarly, we say that a  $ \Bbbk_q $--algebra  $ \cO_q(X) $  is a
{\sl quantization\/}  of the commutative  $ \Bbbk $--algebra  $ \cO(X) $
if it is torsion-free and  $ \, \cO_q(X) \Big/ \! (q-1) \cO_q(X)
\cong \cO(X) \, $.  Then  $ \cO(X) $  is also a  {\sl Poisson\/}  algebra,
called  {\sl semiclassical limit\/}  of  $ \cO_q(X) \, $.
\end{definition}

\vskip2pt

\begin{remark}  \label{pois-bra} {\ }
                                                      \par
  {\it (1)} \,  The technical requirement in {\it (c)\/}  corresponds, in the
context of formal deformations, to require that the algebra is separated; we
also point out that it is satisfied by all quantum function algebras usually
considered in literature.  In any case, it will not be necessary till \S 5.
Moreover, both  {\it (b)\/}  and  {\it (c)\/}  above are automatically
satisfied when  $ \cO_q(G) $  is free as a  $ \bk_q $--module.

\smallskip

  {\it (2)} \, The classical algebra  $ \cO(G) $  inherits from  $ \cO_q(G) $
a Poisson bracket, given as follows: if  $ \, x $,  $ y \in \cO_q(G) \Big/
(q\!-\!1) \, \cO_q(G) \, \cong \cO(G) \, $,  \, then
 \vskip-3pt
  $$  \{x,y\} \; := \; \frac{\, x' \, y' - y' \, x' \,}{q-1}
\mod (q-1) \, \cO_q(G)  $$
 \vskip-3pt
\noindent
 for any lifts  $ \, x' $,  $ y' \in \cO_q(G) \, $  of  $ x $  and  $ y $
respectively.  One checks that this bracket is well-defined, and makes
$ \cO(G) $  into a Poisson Hopf algebra, so that  $ G $  is a Poisson
group.  But  $ G $  already had, by assumption, a Poisson group structure;
then, the requirement in  {\it (a)\/}  above that  $ \; \cO_q(G) \Big/
(q \! - \! 1) \, \cO_q(G) \, \cong \, \cO(G) \; $  as  {\sl Poisson\/}
Hopf algebra amounts to say, in particular, that the two Poisson group
structures of  $ G $  are isomorphic.
                                       \par
   On the other hand, if we start without asking  $ G $  to have
a Poisson group structure, then the previous analysis tells that if
a quantization  $ \cO_q(G) $  exists, then it automatically endows
$ \cO(G) $  with a Poisson algebra structure.  And similarly for a
quantization  $ \cO_q(X) $  of a commutative algebra  $ \cO(X) \, $.
\end{remark}

\vskip2pt

\subsection{Quantum subgroups and quantum coisotropic subgroups}
Our second
step is to introduce the notions of  {\sl quantum coisotropic subgroup\/}
and of  {\sl quantum subgroup},  the former being weaker than the latter.

\vskip1pt

\begin{definition} \label{qcoisosg}
 By  {\sl quantum coisotropic subgroup\/}  of  $ \cO_q(G) $
we mean a  $\bk_q$-co\-algebra  $ \cO_q(H) \, $,  along
with a projection  $ \, \pi : \cO_q(G)
\relbar\joinrel\twoheadrightarrow \cO_q(H) \, $,
such that
                                       \par
   {\it (a)} \;  $ \cO_q(H) $  is torsion-free, as a  $ \bk_q $--module;
                                       \par
   {\it (b)} \;  $ \pi $  is a  $ \bk_q $--coalgebra  (epi)morphism;
                                       \par
   {\it (c)} \;  $ \pi $  is an  $\cO_q(G)$--module  (epi)morphism,
where  $ \cO_q(H) $  has the  $ \cO_q(G) $--module  structure induced
by $\pi$, that is  $ \, f \cdot \pi(g) = \pi(fg) \, $.

\medskip

  {\sl If, in addition},  $ \cO_q(H) $ is a Hopf algebra and  $ \pi $
is a Hopf algebra morphism, and for  $ \; I_H := (q-1) \, \cO_q(H) +
\text{\it Ker}\,\big(\epsilon_{\cO_q(H)}\big) \; $  we have
   $$  {\textstyle \bigcap\limits_{n \geq 0}} \, I_H^{\,n}
\; = \;  {\textstyle \bigcap\limits_{n \geq 0}} \,
(q\!-\!1)^n \, \cO_q(H)   \leqno \indent \text{\it (e)}  $$
then we say that  $ \cO_q(H) $  is a  {\sl quantum subgroup\/}  of
$ G \, $.
                                         \par
   For later use, we introduce also the notation  $ \, I_q(H) := \text{\it Ker}\,(\pi) \; $.
\end{definition}

\medskip

\begin{remark}
                                             \hfill\break
 \smallskip   \indent
   {\it (1)} \,  $ I_q(H) := \text{\it Ker}\,(\pi) \, $  satisfies
$ \, I_q(H) \bigcap \, (q - 1) \cO_q(G) = (q - 1) \, I_q(H)
\; $.  So the specialization of  $ I_q(H) $  at  $ \, q = 1 \, $,
i.e.~$ \, I_1(H) := I_q(H) \big/ (q-1) \, I_q(H) \, $,
coincides with the image of  $ I_q(H) $  under the specialization  $ \,
\cO_q(G) \big/ (q-1) \, \cO_q(G) \cong \cO(G) \, $  of  $ \cO_q(G)
\, $,  which is  $ \, I_q(H) \Big/ \big( I_q(H) \bigcap \, (q\!-\!1)
\cO_q(G) \big) \, $.
                                                \par
 \smallskip
   {\it (2)} \, Conditions  {\it (b)\/}  and  {\it (c)\/}  imply that
$ I_q(H) $  is a two-sided coideal and a left ideal of  $ \cO_q(G) \, $.
Then, by  {\it (1)},  the specialization  $ I_1(H) $  is a coideal and
a (two-sided) ideal in the commutative ring $ \cO(G) \, $.  Moreover,
$ I_1(H) $  equals the kernel of  $ \, \pi_1 \colon \cO(G) = \cO_1(G)
\relbar\joinrel\twoheadrightarrow \cO_1(H) \, $,  the specialization  of
$ \pi $  at  $ \, q = 1 \, $,  where  $ \, \cO_1(H) := \cO_q(H) \big/
(q\!-\!1) \, \cO_q(H) \, $  is the specialization of  $ \cO_q(H) \, $.
So  $ \cO_1(H) $  admits the unique quotient Hopf algebra structure such
that  $ \pi_1 $  is the canonical Hopf algebra epimorphism.  In particular,
$ \cO_1(H) $  is the func\-tion algebra  $ \cO(H) $  of some
   \hbox{closed algebraic subgroup  $ H $  of  $ G \, $,
and  $ \, I_1(H) = $}
 $ \text{\it Ker}\,\big( \pi_1 \, \colon \cO(G)
\relbar\joinrel\twoheadrightarrow \cO(H) \big) = I(H) \, $,
whence the terminology and notation.
                                      \par
   In the Hopf algebra language, conditions {\it (b)\/}  and  {\it (c)\/}
are expressed
%
%
 by saying that  $ \cO_q(H) $  is an  {\sl  $ \cO_q(G) $--module
coalgebra},  that is a coalgebra and  $ \cO_q(G) $--module  such
that  $ \Delta_{\cO_q(H)} $  and  $ \epsilon_{\cO_q(H)} $  are
$ \cO_q(G) $--module  morphisms.
                                          \par
 \smallskip
   {\it (3)} \, Assumptions at the quantum level imply properties for
specializations.  So, the semiclassical specialization of a quantum
coisotropic subgroup is (the function algebra of) a  {\sl coisotropic\/}
subgroup, because  $ \, I_1(H) = \text{\it Ker}\,(\pi_1) \, $  is a
Poisson subalgebra of  $ \cO(G) \, $.  On the other hand, the
specialization of a quantum subgroup instead is (the function
algebra of) a  {\sl Poisson\/} subgroup.
                                      \par
   At the semiclassical level there are many examples of coisotropic
subgroups, among which only a few are Poisson subgroups.  This is a
key motivation to focus on the more general setting of quantum
coisotropic subgroups.
                                          \par
 \smallskip
   {\it (4)} \, A quantum coisotropic subgroup
$ \cO_q(H) $  is by no means a ``quantum group'', in the sense of Definition
\ref{quant-G},  unless it is a quantum subgroup.
\end{remark}

\medskip

\subsection{Quantum line bundles} \label{q-lin-bun}
We now want to carry across to the quantum
setting the notion of embedding  $ \, G\big/H \lhook\joinrel\longrightarrow
\mathbb{P}^N \, $  associated to a line bundle $\cL$
that we assume to be very ample.  The idea is to transfer to
this framework the description (2.4) of  $ \Ogh_n $  given in terms of an
element  $ \, t \in \cO(G) \, $  as in  Remark \ref{desc-semiinv}
and Proposition \ref{t}.  Thus,
the starting point will be a quantization of such an element  $ t \, $
that we will call a pre-quantum section.
                                             \par
   Given  $ G $  and  $ H $  as in  \S \ref{proj-case},
we assume that quantizations of them be given, i.e.~we are given
$ \cO_q(G) $,  $ \cO_q(H) $  and  $ \, \pi \, \colon \cO_q(G)
\relbar\joinrel\twoheadrightarrow \cO_q(H) \, $  as in Definitions
\ref{quant-G},  \ref{qcoisosg}.  To simplify notation, hereafter we
shall also write  $ \, \overline{\ell} := \pi(\ell) \in \cO_q(H)
\, $  for every  $ \, \ell \in \cO_q(G) \, $.
                                             \par
   Moreover, we assume that an element  $ \, t \in \cO(G) \, $  as in  Remark \ref{desc-semiinv},  and the corresponding closed embedding  $ \, G\big/H
\lhook\joinrel\longrightarrow \mathbb{P}^N \, $,  be given as in  Proposition
\ref{t}  (so, in particular,  $ t $  is a section of the line
bundle  $ \cL $  on  $ G \big/ H \, $).
                                             \par
   We define a quantization of the latter setup as follows:

\smallskip

\begin{definition}  \label{pre-q}
  We define  {\sl pre-quantization of\/}  $ \, t \, $,  or  {\sl
pre-quantum section\/}  of the line bundle  $ \cL \, $  on  $ G \big/ H $
(given by  $ t \, $),  any  $ \, d \in \cO_q(G) \, $  such that
 \vskip4pt
   {\it (a)}  \quad  $ \Delta_\pi(d) = d \otimes \pi(d) \; $,  \quad  i.e.
\quad  $ \Delta(d) \in \big( d \otimes d + \cO_q(G) \otimes I_q(H) \big) $
 \vskip4pt
   {\it (b)}  \qquad  $ d \mod (q\!-\!1) \, \cO_q(G)
\,\; = \;  t  \qquad  \big( \in \cO(G) \,\big) \; $
                                                   \par
\noindent
 with respect to the identification  $ \; \cO_q(G) \Big/ (q\!-\!1)
\, \cO_q(G) \, \cong \, \cO(G) \; $.
\end{definition}

\smallskip

\begin{remark} \label{character}  \quad
 \vskip3pt
   {\it (a)} \, Given a pre-quantum section  $ d \, $,  property
{\it (a)\/}  in Definition  \ref{pre-q}  implies that  $ \, \pi(d)
= \overline{d} \, $  is a group-like element in  $ \cO_q(H) $.
Therefore, it defines a one-dimensional corepresentation of
$ \cO_q(H) $,  namely
  $$  \rho_d : \, \bk_q \longrightarrow \bk_q \otimes_{\bk_q} \cO_q(H)
\;\; ,  \qquad  1 \mapsto 1 \otimes \overline{d}  $$
which
%
%
 gives back, modulo  $ (q-1) $,  the one-dimensional representation
of  $ \cO(H) $
  $$  \rho_\lambda : \, \bk \longrightarrow \bk \otimes_{\bk} \cO(H)
\;\; ,  \qquad  1 \mapsto 1 \otimes \lambda  $$
corresponding to the character  $ \, \lambda = \pi(t) \, $  of the
group  $ H $  we started from.
 \vskip3pt
   {\it (b)} \, In the classical setup, having the character  $ \lambda $
is equivalent to having a Hopf algebra morphism  $ \, \bk\big[x,x^{-1}\big]
\longrightarrow \cO(H) \, $,  given by  $ \, x^z \mapsto \lambda^z \;
(z \in \Z \,) $.  Indeed, this occurs because the powers  $ \lambda^z $
do exist in  $ \cO(H) $,  and are group-like because  $ \lambda $  is.
In fact, each one of them corresponds to a one-dimensional
corepresentation, namely the  $ z $--th  tensor power of
$ \rho_\lambda $
  $$  {\rho_\lambda}^{\otimes z} = \rho_{\lambda^z} :
\, \bk \longrightarrow \bk \otimes_{\bk} \cO(H) \;\; ,
\qquad  1 \mapsto 1 \otimes \lambda^z  $$
   \indent   On the other hand, in the quantum setup there is no natural analogue, since the powers
   $ {\overline{d}}^{\,z} $  are not even defined in  $ \cO_q(H) $   --- which is not an algebra! ---
   nor we can assume (would we define them in some way) that they are group-like.  This
means that we miss somehow the ``tensor powers'' of  $ \rho_d \,
$. In \cite{bcgst} one can find an example of a countable family
of group-like elements in a quantum coisotropic subgroups, which
are not obtained by projecting powers of the same element, but
quantize a classical character.
                                                        \par
   However, when  $ \cO_q(H) $  is a quantum subgroup instead, it is
a Hopf algebra, hence the group-like  $ \overline{d} $  is invertible,
and all powers  $ {\overline{d}}^{\,z} $  exist, and are group-like in
$ \cO_q(H) $.  So we do have all ``tensor power corepresentations''
  $$  {\rho_d}^{\otimes z} : \; \bk_q \longrightarrow
\bk_q \otimes_{\bk_q} \cO_q(H) \; ,  \qquad
1 \mapsto 1 \otimes {\overline{d}}^{\,z}  $$
which in turn means that {\sl having  $ \overline{d} $  is equivalent
to having a Hopf\/  $ \bk_q $--algebra  morphism\/}  $ \, \bk_q \big[ x,
x^{-1} \big] \longrightarrow \cO_q(H) \, $,  given by  $ \, x^z \mapsto
{\overline{d}}^{\,z} \; (z \in \Z \,) $.
                                                \par
   Moreover, notice also that  $ \, {\overline{d}}^{\,z} = \overline{d^z}
\, $  for all  $ \, z \geq 0 \, $,  so that  $ {\rho_d}^{\otimes z} $
for  $ \, z \geq 0 \, $  can be directly recovered from the element
$ d^z $  in  $ \cO_q(G) \, $;  thus in the end we can handle everything
working with the elements  $ \, d^n \in \cO_q(G) \, $,  $ \, n \in \N \, $.
\end{remark}

\vskip2pt

\begin{definition}  \label{def_semi-inv}
   Let  $ d \in \cO_q(G)$  be a pre-quantum section on  $ G \big/ H $.
 \vskip3pt
   {\it (a)} \, We call  {\sl  $ d $--semi-invariants of degree}
$ \, n \, $  the elements of the set
  $$  \begin{array}{rl}
   \Oqgh_n \,  &  := \;
 \Big\{\, \ell \in \cO_q(G) \;\Big|\, \Delta_\pi(\ell) =
\ell \otimes \pi\big(d^n\big) \Big\}  \; =  \\
%
%
               &  \hskip-3pt \phantom{|}= \;
 \Big\{\, \ell \in \cO_q(G) \;\Big|\, \Delta(\ell) \in
\big( \ell \otimes d^n + \cO_q(G) \otimes I_q(H) \big) \Big\}
\end{array}  $$
 \vskip3pt
   {\it (b)} \, We call  {\sl  $ d $--semi-invariants\/}  the elements
of the set
  $$  \Oqgh  \; := \;  {\textstyle \sum\limits_{n \in \N}} \Oqgh_n  $$
\end{definition}

\vskip1pt

   It is clear that each  $ \Oqgh_n $  is a  $ \bk_q $--submodule  of
$ \cO_q(G) $,  hence the same holds for  $ \Oqgh $.  We shall now see
some further properties of these modules, which eventually will tell
us that   --- under suitable, additional assumptions ---   we can
take  $ \Oqgh $  as a quantization of  $ \Ogh $.
%
%

\vskip11pt

\begin{lemma} \label{prop-d_semiinv-d}
  Let  $ \, d \in \cO_q(G) \, $  be a pre-quantum section on
$ \, G \big/ H \, $.  Then
  $$  d \in \Oqgh_1 \; ,  \;\; \text{i.e.  $ \; d $  is semi-invariant of degree 1.}
\leqno \indent \text{(a)}  $$
   \indent   (b) \; for any  $ \, n \in \N \, $,  and any
$ \, \ell \in \Oqgh_n \, $,  we have
  $ \;\; \overline{\ell} \, = \, \ep(\ell\,) \, \overline{d^n} \; $.
%
%
 \vskip7pt
   (c) \; the map  $ \, \pi : \cO_q(G)
\relbar\joinrel\twoheadrightarrow \cO_q(H) \, $  restricts
to a  $ \bk_q $--module  epimorphism
  $$  \pi' \, \colon \Oqgh \relbar\joinrel\twoheadrightarrow \text{\it Span}_{\,\bk_q}
\Big( \big\{\, \overline{d^n} \,\big\}_{n \in \N} \Big) \; .  $$
\end{lemma}

\begin{proof}
%
%
   The only statement which needs a proof here is  {\it (b)},  which quickly
follows applying  $ (\epsilon \otimes \text{\it id}\,) $  to both sides of
$ \; \Delta_\pi(\ell) = \ell \otimes \overline{d^n} \; $.
\end{proof}

\smallskip

\begin{remark}  \label{spec_semi-inv-n}
 The semi-invariants have a good arithmetic property, which ensures
that the specialization of  $ \, \Oqgh_n \subset \cO_q(G) \, $  (and
of  $ \Oqgh $)  at  $ \, q = 1 \, $  will be consistent with that of
$ \cO_q(G) $  itself.  Namely, given  $ \, n \in \N \, $,  since
$ \Oqgh_n $  is defined by  $ \bk_q $--linear  conditions we find
at once
  $$  \Oqgh_n \, {\textstyle \bigcap} \; c \, \cO_q(G)  \; = \;
c \, \Oqgh_n  \qquad  \forall \;\; c \in \bk_q  $$
and in particular,  $ \;\; \Oqgh_n \, {\textstyle \bigcap} \,
(q\!-\!1) \, \cO_q(G) \, = \, (q\!-\!1) \, \Oqgh_n \;\; $.

\smallskip

%
%
\end{remark}

\medskip

   Next result shows that each  $ \Oqgh_n $  is a left coideal
of  $ \cO_q(G) $,  hence it bears a structure of left
$ \cO_q(G) $--comodule.

\medskip

\begin{proposition}  \label{Oqgh-coid}
  Every  $ \Oqgh_n $  is a left coideal of  $ \, \cO_q(G) \, $,  that is
   $$  \Delta \Big( \Oqgh_n \Big) \, \subseteq \, \cO_q(G) \otimes \Oqgh_n
\qquad  \forall \;\, n \in \N  $$
so that  $ \, \Delta\big|_{\cO_q(G/H)_n} $  makes  $ \, \Oqgh_n $  into
a left  $ \cO_q(G) $--comodule.  Therefore,  $ \Oqgh $  is a left coideal of  $ \, \cO_q(G) \, $,  hence a left  $ \cO_q(G) $--comodule.
\end{proposition}

\begin{proof}  Let  $ \, \cO := \cO_q(G) \, $,  $ \, \cO_n := \Oqgh_n \, $,
and set  $ \, \cO' := \bk(q) \otimes_{\bk_q} \cO \, $,  $ \, \cO'_n :=
 \bk(q) \otimes_{\bk_q} \cO_n \, $.  Then  $ \, \cO \otimes_{\bk_q}
\cO \, $  naturally embeds into  $ \, \cO' \otimes_{\bk(q)} \cO' \, $,
because  $ \cO $  and  $ \cO'_n $  are torsion free as
$ \bk_q $--modules.  Using this embedding, given any
$ \, \ell \in \Oqgh_n \, $  we can always write  $ \,
\Delta(\ell\,) = \sum_i g'_i \otimes h'_i \, $  for some
$ \, g'_i \, , h'_i \in \cO' \, $  such that the  $ g'_i $'s
are all linearly independent, and similarly for the  $ h'_i $'s.
Then (taking a common denominator) there exists  $ \, c \in \bk_q \, $
such that
 $ \; \Delta(\ell\,) = c^{-1} \, \sum\limits_i \, g_i
\otimes h_i \; $,  \;
%
%
with the  $ g_i $'s  in  $ \cO_q(G) $  being linearly independent, and
the  $ h_i $'s  in  $ \cO_q(G) $  which are linearly independent too.
                                                     \par
   We shall now prove that
  $$  {\textstyle \sum\limits_i} \, g_i \otimes {\textstyle
\sum\limits_{(h_i)}} \, \big( h_i \big)_{\!(1)} \otimes
\overline{\big( h_i \big)_{\!(2)}}  \; = \;
{\textstyle \sum\limits_i} \, g_i
\otimes h_i \otimes \overline{d^n}  \quad .   \eqno (3.1)  $$
Indeed, the left-hand side in (3.1) is just the image of
$ \, c \, \ell \, $  via the map
  $$  \Big( \textit{id} \otimes \big( (\textit{id} \otimes \pi ) \circ
\Delta \big) \Big) \circ \Delta  \, = \,  \big( \textit{id} \otimes
\textit{id} \otimes \pi \big) \circ \big( \textit{id} \otimes
\Delta \big) \circ \Delta  $$
By coassociativity of  $ \Delta \, $,  the latter map coincides with
  $$  \big( \textit{id} \otimes \textit{id} \otimes \pi \big) \circ
\big( \textit{id} \otimes \Delta \big) \circ \Delta  =
\big( \Delta \otimes \pi \big) \circ \Delta  =
\big( \Delta \otimes \textit{id} \,\big) \circ
\big( \textit{id} \otimes \pi \big) \circ \Delta  =
\big( \Delta \otimes \textit{id} \,\big) \circ \Delta_\pi  $$
and now the last map applied to  $ \, c \, \ell \, $  gives
  $$  \big( \Delta \otimes \textit{id} \,\big) \big( \Delta_\pi (c\,\ell\,) \big)
\;\; {\buildrel \circledast \over =} \;\;
\big( \Delta \otimes \textit{id} \,\big)
\big( c \, \ell \otimes \overline{d^n} \,\big)  \; = \;
{\textstyle \sum_i} \, g_i \otimes h_i \otimes \overline{d^n}  $$
%
%
%
where  $ \, {\buildrel \circledast \over =} \, $  follows from the assumption
$ \, \ell \in \Oqgh_n \, $,  which implies  $ \, c \, \ell \in \Oqgh_n \, $
as well.  This eventually gives the right-hand side of (3.1), q.e.d.
                                                     \par
   Now, because of the linear independence of the  $ h_i $'s,  the
identity (3.1) implies that all the  $ h_i $'s  satisfy  $ \;
\sum_{(h_i)} \, \big( h_i \big)_{\!(1)} \otimes
\overline{\big( h_i \big)_{\!(2)}}  \; = \;
h_i \otimes \overline{d^n} \; $,  \; which means exactly
$ \, h_i \in \Oqgh_n \, $,  for every index  $ i \, $.  Thus we have
  $$  \Delta\big(c\,\ell\,\big)  \, = \,  {\textstyle \sum_i} \,
g_i \otimes h_i \;\; ,   \eqno g_i \in \cO_q(G) \, , \; h_i \in
\Oqgh_n \; ,  $$
the  $ g_i $'s,  resp.~the  $ h_i $'s,  being linearly independent,
and also
  $$  \Delta\big(c\,\ell\,\big) \in c \, \cO_q(G) \otimes \cO_q(G)
\;\; .  $$
These two conditions imply  $  \; h_i \, \in \, \Oqgh_n \, \bigcap \;
c \, \cO_q(G) \, = \, c \, \Oqgh_n \; $  thanks to Remark  \ref{spec_semi-inv-n}.
Therefore  $ \, \Delta(\ell\,) \in \cO_q(G) \otimes \Oqgh_n \, $.  Finally,
the claim for  $ \Oqgh $  follows from that for the  $ \Oqgh_n $'s.
\end{proof}

 \medskip

   The above construction provides us with a reasonable candidate
for a quantum analogue of  $ \Ogh \, $,  namely the space of the
$ d $--semi-invariants  $ \Oqgh \, $,  which we proved has many
important properties.  Nevertheless, we still would like  $ \Oqgh $
to verify three more key properties, namely:
 \vskip3pt
   {\it (a)} \;  $ \Oqgh $  is a  {\sl subalgebra\/}  of  $ \cO_q(G) \, $;
 \vskip3pt
   {\it (b)} \;  $ \Oqgh $  is a  {\sl graded\/}  object, its
$ n $--th  (for all  $ \, n \! \in \! \N \, $)  graded summand
being  $ \Oqgh_n \, $;
 \vskip3pt
   {\it (c)} \; the grading is compatible with all other structures, so
$ \Oqgh $  is a  {\sl graded\/}  $ \cO_q(G) $--comodule  algebra  (when
$ \cO_q(G) $  is given the trivial grading).
 \vskip7pt
   Indeed, we are still missing these properties so far.  In order to
have them, an additional property must be required to the pre-quantum
section  $ d $  we started from.  This is provided by the following
result.

 \medskip

\begin{proposition}  \label{cond-Oqgh-subalg}
   Set  $ \, I_q(H) := \text{\it Ker}\,(\pi) \, $,  and let
$ d $  be a pre-quantum section on  $ G\big/H \, $.  Then the
following properties are equivalent:
 \vskip5pt
(a) \,  $ \Oqgh_r \cdot \Oqgh_s \, \subseteq \, \Oqgh_{r+s} \;\; $
for all  $ \, r, s \in \N \, $,  hence, in particular,  $ \Oqgh $
is a  $ \bk_q $--subalgebra  of  $ \cO_q(G) \, $;
 \vskip5pt
   (b) \,  $ \overline{[d,f]} = \bar{0} \; $  in  $ \cO_q(H) \, $,
\; for all  $ \, f \in \Oqgh \, $;
 \vskip5pt
   (c) \,  $ \big[\, d\, , \, \Oqgh \,\big]
\, \subseteq \, I_q(H) \; $.
\end{proposition}

\begin{proof}  {\it (b)  $ \! \Longrightarrow \! $  (a)\,:}   For
any  $ \, r, s \in \N \, $,  pick  $ \, f \in \Oqgh_r \, $,  $ g \,
\in \Oqgh_s \, $.  Then   --- by Proposition  \ref{Oqgh-coid}  ---
for  $ \, \Delta(f) = \sum_{(f)} f_{(1)} \otimes f_{(2)} \, $  and
$ \, \Delta(g) = \sum_{(g)} g_{(1)} \otimes g_{(2)} \; $  we have
$ \, f_{(2)} \in \Oqgh_r \, $,  $ \, g_{(2)} \in \Oqgh_s \, $.  This,
along with assumption  {\it (b)\/}  to get the equality  $ \,
{\buildrel \circledast \over =} \, $,  yields the chain of
identities
  $$  \begin{array}{rl}  \hskip-15pt
  &  \hskip-9pt
    \Delta_\pi(fg) \; =
{\textstyle \sum\limits_{(f),(g)}} f_{(1)} \, g_{(1)} \otimes
\overline{f_{(2)} \, g_{(2)}} \; = {\textstyle \sum\limits_{(f),(g)}}
f_{(1)} \, g_{(1)} \otimes f_{(2)}.\overline{g_{(2)}}
\; =  \\
 &  \phantom{\bigg|} = \!
 {\textstyle \sum\limits_{(f),(g)}} \!\! \big( f_{(1)} \otimes
f_{(2)} \big).\big(g_{(1)} \otimes \overline{g_{(2)}}\big)
=  \Big( {\textstyle \sum_{(f)}} \, f_{(1)} \otimes f_{(2)}
\!\Big).\Big({\textstyle \sum_{(g)}} \, g_{(1)} \otimes
\overline{g_{(2)}}\Big)  \! =  \\
 &  \phantom{\bigg|} =
 \Big( {\textstyle \sum_{(f)}} \, f_{(1)} \otimes f_{(2)}
\Big).\Big( g \otimes \overline{d^s} \,\Big) =
{\textstyle \sum_{(f)}} \; f_{(1)} \, g \otimes
f_{(2)}.\,\overline{d^s}  \; =  \\
 &  \phantom{\bigg|} =
 {\textstyle \sum_{(f)}} \, f_{(1)} \, g \otimes
\overline{f_{(2)} \, d^s}  \;\;  {\buildrel \circledast
\over =} \;\,  {\textstyle \sum_{(f)}} \,
f_{(1)} \, g \otimes \overline{d^s \, f_{(2)}}  \; =  \\
 &  \phantom{\bigg|} =
 {\textstyle \sum_{(f)}} \, f_{(1)} \, g \otimes
\big( d^s \, f_{(2)} \big).\overline{1}  \, = \,
\Big(\! \big( 1 \otimes d^s \big) \cdot {\textstyle \sum_{(f)}} \;
f_{(1)} \otimes f_{(2)} \,\Big).\big( g \otimes \overline{1} \,\big) \, =
\end{array}  $$
  $$  \begin{array}{rl}  \hskip-15pt
 &  \phantom{|} =
\Big(\! \big( 1 \otimes d^s \big) \, \Delta(f) \Big).\big( g \otimes \overline{1}
\,\big)  \, = \,  \Big(\! \big( 1 \otimes d^s \big) \, \big( f \otimes d^r +
\varphi \otimes \eta \big) \Big).\big( g \otimes \overline{1} \,\big)  \, =  \\
 &  \phantom{\bigg|} =
\,  \big( 1 \otimes d^s \big).\big( f \, g \otimes \overline{d^r} +
\varphi \, g \otimes \overline{\eta} \,\big)  \, = \;
f \, g \otimes \overline{d^{s+r}}  \; = \;
f \, g \otimes \pi\big(d^{s+r}\big)
\end{array}  $$
 \vskip-3pt
\noindent
for some suitable  $ \, \varphi \in \cO_q(G) \, $,  $ \, \eta \in I_q(H) \, $,
with notation  $ \, (x \otimes y).\big( a \otimes \overline{b} \,\big) :=
(x \, a) \otimes \big( y.\,\overline{b} \,\big) \, $  referring to the action of
$ \, \cO_q(G) \otimes \cO_q(G) \, $  onto  $ \, \cO_q(G) \otimes \cO_q(H) \, $
induced by the action of  $ \cO_q(G) $  onto  $ \cO_q(H) \, $,  via  $ \pi \, $,
and onto itself, via left regular representation.  So  $ f \, g $  is also
$ d $--semi-invariant,  of degree  $ \, r + s \, $,  q.e.d.
 \vskip3pt
   {\it (a)  $ \! \Longrightarrow \! $  (b)\,:}  \; Assume that  {\it
(a)\/}  holds.  Then for  $ \, f \in \Oqgh_n \, $  we have  $ \;
d \, f \, $,  $ \, f \, d \, \in \Oqgh_{n+1} \, $,  and so
$ \; [d,f] \in \Oqgh_{n+1} \, $.  Then the identity
 \vskip2pt
\noindent
  $$ \overline{[d,f]}  \; = \;  \ep\big([d,f]\big) \,
\overline{d^{\,\partial(f)+1}}  $$
 \vskip-2pt
\noindent
 holds, by Proposition  \ref{prop-d_semiinv-d}.  But clearly  $ \,
\ep\big([d,f]\big) = 0 \, $,  hence  $ \; \overline{[d,f]} \, = \,
\overline{0} \, $,  that is  $ \, [d,f] \in I_q(H) \, $.  The
outcome is  $ \, \big[ d, \Oqgh \big] \subseteq I_q(H) \, $,
\, q.e.d.
 \vskip3pt
   {\it (b)  $ \! \Longleftrightarrow \! $  (c)\,:}  \,  This is just
a matter of rephrasing.
\end{proof}

\vskip0pt

\begin{definition}  \label{q-sect}
  We call  {\sl quantization of\/}  $ \, t \, $,  or  {\sl quantum
section\/}  (of the line bundle  $ \cL \, $)  on  $ G \big/ H \, $,
any pre-quantum section  $ d $  of  $ G\big/H $  (cf.~Definition
\ref{pre-q})  which satisfies any one of the equivalent conditions
in  Proposition \ref{cond-Oqgh-subalg}.
\end{definition}

\vskip3pt

%
%
   The following result gives a criterion to
detect quantum sections, and shows that for quantum subgroups they
are
 just
pre-quantum sections.

\medskip

\begin{proposition}  \label{crit_q-sect}  \quad
 \vskip2pt
\noindent
 {\it (a)}  Let  $ d $  be a pre-quantum section on
$ \, G \big/H \, $.
                                             \par
   If  $ \, \big[\, d\, , \, I_q(H) \,\big] \, \subseteq \, I_q(H) \; $,
\, then  $ d $  is a quantum section.
 \vskip4pt
\noindent
 {\it (b)}  Let  $ \cO_q(H) $  be a quantum subgroup.
                                             \par
   Then any pre-quantum section on  $ \, G\big/H \, $  is a quantum
section.
\end{proposition}

\begin{proof}  {\it (a)}  Pick any  $ \, f \in \Oqgh_r \, $,  $ \, g
\in \Oqgh_s \; $.  Definition  \ref{def_semi-inv}  gives
  $$  \begin{matrix}
   \Delta(f) \, = \, f \otimes d^r + f_1 \otimes \phi_1 \;\; ,
 &  \qquad  f_1 \in \cO_q(G) \; , \;\; \phi_1 \in I_q(H)  \\
   \Delta(g) \, = \, g \otimes d^s + g_1 \otimes \gamma_1 \;\; ,
 &  \qquad  g_1 \in \cO_q(G) \; , \;\; \gamma_1 \in I_q(H)
   \end{matrix}  $$
Therefore, for the product  $ \, f \, g \, $  we have
  $$  \displaylines{
   \Delta\big(f\,g\big) \, = \, \Delta(f) \, \Delta(g) \, = \,
\big( f \otimes d^r + f_1 \otimes \phi_1 \big) \,
\big( g \otimes d^s + g_1 \otimes \gamma_1 \big) \, =
\hfill  \cr
   \hfill   = \, f \, g \otimes d^{r+s} \,
+ \, f \, g_1 \otimes d^r \, \gamma_1 \,
+ \,  f_1 \, g \otimes \phi_1 \, d^s \, +
\, f_1 \, g_1 \otimes \phi_1 \, \gamma_1 }  $$
Now,  $ \, d^r \, \gamma_1 \, , \, \phi_1 \, \gamma_1 \in I_q(H) \, $
because  $ I_q(H) $  is a (left)  $ \cO_q(G) $--submodule,  and
  $$  \phi_1 \, d^s  \, = \,  d^s \, \phi_1
+ \big[ \phi_1 \, , d^s \big]  \, \in \,
\Big( I_q(H) + I_q(H) \Big)  \, = \,  I_q(H)  $$
because, in addition,  $ \, \big[d,I_q(H)\big] \subseteq I_q(H) \, $,
by assumption  {\it (d)}.  Thus
  $$  \Delta\big(f\,g\big)  \, \in \,  \Big( f \, g \otimes d^{r+s}
\, + \, \cO_q(G) \otimes I_q(H) \Big)  $$
which means exactly  $ \, f \, g \in \Oqgh_{r+s} \; $,  \, by
Definition  \ref{def_semi-inv}  again.  Thus condition  {\it (a)\/}
of Proposition  \ref{cond-Oqgh-subalg}  holds, hence we conclude by
Definition  \ref{q-sect}.
 \vskip5pt
   {\it (2)} \, If  $ \cO_q(H) $  is a quantum subgroup, then
$ I_q(H) $  is a two-sided ideal.  Therefore,  $ \, \big[\, d\, ,
\, I_q(H) \,\big] \, \subseteq \, I_q(H) \; $,  \, hence by  {\it
(1)\/}  we get the claim.
\end{proof}

\medskip

   The following result records yet another feature of quantum sections:

\medskip

\begin{lemma}  \label{d-powers}
   Let  $ d $  be a quantum section on  $ \, G\big/H \, $.
                                                \par
   Then  $ \, d^n \in \Oqgh_n \, $,  and  $ \, \overline{d^n} \, $
is group-like in  $ \cO_q(H) \, $,  \, for all  $ \, n \in \N \, $.
Moreover,  $ \, \text{\it Span}_{\,\bk_q} \Big( \big\{\, \overline{d^n}
\,\big\}_{n \in \N} \Big) $  is a  $ \bk_q $--subcoalgebra  of
$ \, \cO_q(H) \, $,  and
  $$  \text{\it Span}_{\,\bk_q} \Big(
\big\{\, \overline{d^n} \,\big\}_{n \in \N} \Big)
\; = \;  {\textstyle \bigoplus_{n \in \N}} \, \bk_q \,
\overline{d^n} \quad .  $$
\end{lemma}

\begin{proof}  By Definition  \ref{q-sect}  and condition  {\it (a)\/}
in Proposition  \ref{cond-Oqgh-subalg}  we have  $ \, d^n \in \Oqgh_n \, $,
for all  $ \, n \in \N \, $.  This means  $ \, \Delta\big(d^n\big)
\in \big( d^n \otimes d^n + \cO_q(G) \otimes I_q(H) \big) \, $,
whence   --- as  $ \, \pi \, \colon \cO_q(G)
\relbar\joinrel\twoheadrightarrow \cO_q(H) \, $
is a  {\sl coalgebra\/}  morphism ---   we get
  $$  \Delta\big( \overline{d^n} \,\big) \, = \,
\Delta\big(\pi\big(d^n\big)\big) \, = \,
(\pi \otimes \pi) \big( \Delta\big(d^n\big) \big)
\, = \, \pi\big(d^n\big) \otimes \pi\big(d^n\big)
\, = \, \overline{d^n} \otimes \overline{d^n}  $$
thus the  $ \overline{d^n} $'s  are group-like, and different from each other
because the  $ \, t^n = d^n{\big|}_{q=1} \, $  are.  Finally, this implies
that  $ \text{\it Span}_{\,\bk_q} \Big( \big\{\, \overline{d^n} \,\big\}_{n
\in \N} \Big) $  is a  $ \bk_q $--subcoalgebra  of  $ \cO_q(H) \, $,  and
also
%
%
 that the  $ \overline{d^n} $'s  are linearly independent,
which eventually gives  $ \; \text{\it Span}_{\,\bk_q} \Big( \big\{\,
\overline{d^n} \,\big\}_{n \in \N} \Big) \, = \, {\textstyle
\bigoplus_{n \in \N}} \, \bk_q \, \overline{d^n} \; $.
\end{proof}

\medskip

   Gathering all together the previous results, we can now show that
semi-invariants built out of a quantum section satisfy all properties
we look for:

\smallskip

\begin{theorem}  \label{Oqgh-graded}
 \vskip3pt
   Let  $ d $  be a quantum section on  $ \, G\big/H \, $.  Then
 \vskip3pt
   {\it (a)} \;  $ \Oqgh \, $  is a  {\sl graded\/}  $ \bk_q $--module,
its  $ n $--th  graded summand  ($ \, n \! \in \! \N \, $)  being
$ \Oqgh_n \, $;
 \vskip3pt
   {\it (b)} \;  $ \Oqgh \, $  is a  {\sl subalgebra\/}  of
$ \, \cO_q(G) \, $;
 \vskip3pt
   {\it (c)}  the grading in  {\it (a)}  is compatible with all other
structures of  $ \, \Oqgh \, $,  so that  $ \, \Oqgh $  is a  {\sl graded}
$ \, \cO_q(G) $--comodule  algebra,  where we take on  $ \, \cO_q(G) $
the trivial grading;
 \vskip3pt
   {\it (d)} \, for every  $ \, c \in \bk_q \, $,  we have
$ \;\; \Oqgh \, {\textstyle \bigcap} \, c \, \cO_q(G) \, = \,
c \, \Oqgh \;\; $.
%
%
                                            \par
\noindent
 In particular,  $ \;\; \Oqgh \, {\textstyle \bigcap} \,
(q\!-\!1) \, \cO_q(G) \, = \, (q\!-\!1) \, \Oqgh \;\; $.

\end{theorem}

\begin{proof}  {\it (a)} \, We must simply prove that the sum  $ \,
{\textstyle \sum\limits_{n \in \N}} \Oqgh_n \, $  is  {\sl direct},
so that  $ \; \Oqgh := {\textstyle \sum\limits_{n \in \N}} \Oqgh_n =
{\textstyle \bigoplus\limits_{n \in \N}} \, \Oqgh_n \; $.  But this
is an easy consequence of Lemma  \ref{d-powers}.
                                            \par
   Indeed, let  $ \; \sum_{n \in \N} c_n \, f_n = 0 \; $  a linear
dependence relation, with  $ \, f_n \in \Oqgh_n \, $  and  $ \, c_n \in
\bk_q \, $  (almost all zero) for every  $ \, n \in \N \, $.  Applying
$ \, \Delta_\pi \, $  to this relation we get  $ \; \sum_{n \in \N} c_n
\, f_n \otimes \overline{d^n} = 0 \; $.  But the  $ \overline{d^n} $'s,
by Lemma  \ref{d-powers},  are linearly independent; thus  $ \, c_n = 0
\, $  for all  $ n \, $,  q.e.d.
 \vskip3pt
   {\it (b)} \,  This follows directly from Definition  \ref{q-sect}  and
Proposition  \ref{cond-Oqgh-subalg}.
 \vskip3pt
   {\it (c)} \,  This follows (again) from Definition  \ref{q-sect}  and
Proposition  \ref{cond-Oqgh-subalg},  from Proposition  \ref{Oqgh-coid},  and
from the Hopf algebra axioms.
 \vskip3pt
   {\it (d)} \,  This easily follows from the identity  $ \; \Oqgh =
{\textstyle \bigoplus_{n \in \N}} \, \Oqgh_n \; $,  \, given by
claim  {\it (a)},  and from Remark  \ref{spec_semi-inv-n}.
\end{proof}

\vskip2pt

\begin{corollary}  \;  Let  $ d $  be a quantum section on
$ \, G\big/H \, $  (in the sense of Definition  \ref{q-sect}).
Then the restriction of  $ \,\; \pi \, \colon \cO_q(G)
\relbar\joinrel\twoheadrightarrow \cO_q(H) \; $  yields
an epimorphism of  {\sl graded  $ \cO_q(G) $--module
coalgebras}
  $$  \pi' \; \colon \; \Oqgh
\relbar\joinrel\relbar\joinrel\twoheadrightarrow
\text{\it Span}_{\,\bk_q} \Big( \big\{\, \overline{d^n}
\,\big\}_{n \in \N} \Big) \, = \; {\textstyle \bigoplus_{n \in \N}}
\, \bk_q \, \overline{d^n} \;\;\; .  $$
\end{corollary}

\begin{proof}  By Lemma  \ref{prop-d_semiinv-d}  we
know that  $ \pi' $  is a well-defined epimorphism of  {\sl graded\/}
$ \bk_q $--modules.  The rest follows from  $ \Oqgh $  being a
subalgebra of  $ \cO_q(G) \, $,  and  $ \pi $  being a morphism
of  $ \cO_q(G) $--module  coalgebras.
\end{proof}

\medskip

   Last aspect to take into account is the behavior of  $ \Oqgh $
under specialization at  $ \, q=1 \, $.  The last part of claim
{\it (d)\/}  in  Theorem \ref{Oqgh-graded}  ensures that such
specialization is consistent with that of  $ \cO_q(G) $:
in other words, the embedding  $ \, \Oqgh
\lhook\joinrel\longrightarrow \cO_q(G) \, $  gets down
under specialization to an embedding  $ \; \cO_1\big(G\big/H\big)
:= \Oqgh \Big/ (q-1) \, \Oqgh \lhook\joinrel\longrightarrow \cO(G)
\; $.  The next result tells us something about the specialized
space  $ \cO_1\big(G\big/H\big) $  itself.

\medskip

\begin{proposition}  \label{prop-spec-oqgh}
  $ \, \cO_1\big(G\big/H\big) := \Oqgh \Big/ (q-1) \, \Oqgh \; $  is a graded Poisson subalgebra of  $ \; \Ogh = \bigoplus_{n \in \N} \Ogh_n \; $  and a graded left coideal of  $ \cO(G) \, $,  \, with  $ \; \cO_1\big(G\big/H\big)_n := \Oqgh_n \Big/ (q\!-\!1) \, \Oqgh_n \;\; $  as  $ \, n $--th  graded summand  ($ \, n \in \N \, $).  In particular, it is a left  $ \, \cO(G) $--comodule  algebra.
\end{proposition}

\begin{proof}  For all  $ \, n \in \N \, $,  the restriction to
$ \Oqgh_n $  of the specialization map
  $$  p_1 \, : \, \cO_q(G) \,
\relbar\joinrel\relbar\joinrel\twoheadrightarrow \,
\cO_q(G) \Big/ (q-1) \, \cO_q(G) \; \cong\; \cO(G)  $$
has kernel  $ \; \Oqgh_n \, {\textstyle \bigcap} \, (q-1) \, \cO_q(G)
\, = \, (q-1) \, \Oqgh_n \; $,  \, by  Remark \ref{spec_semi-inv-n}.
This in turn ensures also that the restriction of  $ p_1 $  to  $ \;
\Oqgh = {\textstyle \bigoplus_{n \in \N}} \Oqgh_n \; $  has kernel
$ \; {\textstyle \bigoplus_{n \in \N}} \, (q-1) \, \Oqgh_n \; $,
\, so its image is  $ \, \Oqgh \Big/ (q-1) \, \Oqgh \, $,  \,
i.e.~just the specialization of  $ \Oqgh \, $.  So
  $$  \cO_1\big(G\big/H\big)  \; := \;  \Oqgh \Big/ (q-1) \, \Oqgh
\; = \;  p_1\Big(\Oqgh\Big)  $$
where the right-hand side is a subalgebra of  $ \, p_1\Big(\cO_q(G)\Big)
= \cO(G) \, $.  Moreover, we have also that the specialization maps
preserves the grading, namely
  $$  p_1\Big(\Oqgh\Big) \; = \; p_1 \Big(\, {\textstyle
\bigoplus_{n \in \N}} \, \Oqgh_n \,\Big) \; = \; {\textstyle
\bigoplus_{n \in \N}} \, p_1\Big(\Oqgh_n\Big)  $$
so that  $ \, \cO_1\big(G\big/H\big) \, $  is  {\sl graded},  with
$ n $--th  graded summand
  $$  p_1\Big(\Oqgh_n\Big)  \; = \;  \Oqgh_n \Big/ (q-1) \, \Oqgh_n
\; =: \;  \cO_1\big(G\big/H\big)_n  $$
   \indent   Now  Theorem \ref{Oqgh-graded}  implies at once that
$ \cO_1\big(G\big/H\big) $  is a graded subalgebra left coideal inside
$ \; p_1\big(\cO_q(G)\big) = \cO(G) \, $,  \, hence a graded (left)
$ \cO(G) $--comodule  algebra.  In addition, the identity
  $$  \Oqgh \; {\textstyle \bigcap} \; (q-1) \, \cO_q(G)  \; = \;
(q-1) \, \Oqgh  $$
implies also that the Poisson bracket defined in  $ \Ogh $  starting
from its quantization  $ \Oqgh $   --- see  Remark \ref{pois-bra}  ---
coincides with the restriction to  $ \Ogh $  of the Poisson bracket
similarly induced on  $ \cO(G) $  from  $ \cO_q(G) \, $.  Therefore,
$ \cO_1\big(G\big/H\big) $  is also a Poisson subalgebra of  $ \cO(G) \, $.
 \vskip3pt
   We are only left to prove that the embedding of  $ \cO_1\big(G\big/H\big) $  into  $ \cO(G) $  maps  $ \cO_1\big(G\big/H\big) $  into  $ \Ogh $  respecting the grading
on either side, that is
  $$  \cO_1\big(G\big/H\big)_n \; \subseteq \, {\Ogh}_n   \quad \qquad \forall \;\, n \in \N   \eqno (3.2)  $$
Now, the left-hand side of (3.2) is  $ \; \Oqgh_n \Big/ (q-1) \, \Oqgh_n
\; $, \; with
  $$  \Oqgh_n  \; := \;
\Big\{\, \ell \in \cO_q(G) \;\Big|\, \Delta_\pi(\ell)
= \ell \otimes \pi\big(d^n\big) \Big\}  $$
(cf.~Definition \ref{def_semi-inv}),  while the right-hand side, by (2.4),
is
  $$  \Ogh_n  \; = \;  \Big\{\, f \in \cO(G) \,\;\Big|\;
\Delta_\pi(f) = f \otimes \pi\big(t^n\big) \Big\}  $$
But all specialization maps commute with the coproducts  $ \Delta $  and
with the (quantum and classical) maps  $ \pi \, $,  and the specialization
of each  $ d^n $  is nothing but  $ t^n \, $.  Therefore, we conclude that
(3.2) holds.
\end{proof}

\medskip

   Finally, we can define our ``quantum projective homogeneous spaces''.
\smallskip

\begin{definition}  \label{q-prhspace}
   Let  $ G $  be an algebraic Poisson group,  $ H $  a
coisotropic subgroup, as in  \S \ref{proj-case},  and let  $ \,
\cO_q(G) \, $,  $ \cO_q(H) $  and  $ \, \pi \, \colon \cO_q(G)
\relbar\joinrel\twoheadrightarrow \! \cO_q(H) \, $  be given
(cf.~Definitions  \ref{quant-G}, \ref{qcoisosg}).
Let  $ d $ be a quantum section on  $ G\big/H \, $  (see Definition  \ref{q-sect}, Proposition \ref{t}), in particular
 \vskip4pt
   {\it (a)}  $ \quad  \Delta_\pi(d) \, = \, d \otimes \pi(d) $
 \vskip4pt
   {\it (b)}  $ \quad  d \, \equiv \, t \mod (q-1) \; $,  \; where  $ t $  is a non-zero section of the very ample line bundle on  $ G \big/ H $  giving the embedding into some projectve space.
 \vskip4pt
   Then given
  $$  \begin{array}{rl}
   \Oqgh_n  &  := \;
 \Big\{\, \ell \in \cO_q(G) \;\Big|\, \Delta_\pi(\ell) =
\ell \otimes \pi\big(d^n\big) \Big\}  \; =  \\
               &  \hskip-4pt \phantom{\bigg|}= \;
 \Big\{\, \ell \in \cO_q(G) \;\Big|\, \Delta(\ell) \in
\big( \ell \otimes d^n + \cO_q(G) \otimes I_q(H) \big) \Big\}
\end{array}  $$
we say that
  $$  O_q(G/H)  \; := \;  {\textstyle \bigoplus}_{n \in \N} \, \Oqgh_n  $$
is a  {\sl quantization of}  $ \, \Ogh \, $  if
 \vskip-5pt
  $$  \Oqgh \Big/ (q\!-\!1) \, \Oqgh  \,\; \cong \;\,  \Ogh  $$
 \vskip-1pt
\noindent
 as graded  $ \cO(G) $--module  algebras and as Poisson algebras
over  $ \bk \, $.  We will then refer to $\Oqgh$ as {\sl quantum projective homogeneous space}.
 \vskip3pt
   In particular, we have seen that any such  $ \Oqgh $  has the following properties (Theorem \ref{Oqgh-graded}):
                                                  \par
   \hskip3pt   {\it (I)} \,  it is a graded subalgebra of  $ \, \cO_q(G) \, $;
 \vskip2pt
   {\it (II)} \,  It is a left coideal of  $ \cO_q(G) \, $,  hence a left
$ \cO_q(G) $--comodule  via
  $$  \Delta|_{\cO_q(G/H)}: \cO_q(G/H) \lra \cO_q(G) \otimes
\cO_q(G/H)  $$
\end{definition}

\vskip3pt

\begin{remark}  \label{rem-qprhsp}
   As a matter of fact, the only additional property required
in  Definition \ref{q-prhspace}  is that the embedding of  $ \;
\Oqgh \Big/ (q-1) \, \Oqgh \; $  into  $ \, \Ogh \, $  provided
by  Proposition \ref{prop-spec-oqgh}  be  {\sl onto}.  But actually, as both these are  {\sl graded\/}  algebras, and  $ \Ogh $  is generated in degree one,  {\it this is equivalent to require (only) that is  {\sl onto}  the embedding}
 \vskip-4pt
  $$  \Oqgh_1 \Big/ (q-1) \, \Oqgh_1 \;
\lhook\joinrel\relbar\joinrel\relbar\joinrel\relbar\joinrel\longrightarrow
\; \Ogh_1  $$
 \vskip-1pt
   \indent   This requirement might be seen as the quantum analogue of the
requirement   --- at the semiclassical level, see \S \ref{proj-case}  ---
of having enough global sections of the line bundle  $ \cL $  on  $ G\big/H $
so to have an embedding of  $ G\big/H $  into  $ \mathbb{P}^N \, $.
\end{remark}

\vskip1pt

   In section \ref{Q-Grass} we show that quantum Grassmannians and quantum generalized flag varieties are examples of quantum projective homogeneous space.

\bigskip

\section{The Quantum Duality Principle (QDP)}

\subsection{The QDP philosophy} \label{qdp-philosophy}
   The quantum duality principle (QDP) is a two-fold recipe which allows
to obtain a quantum group dual, in an appropriate sense, to a given one.
                                                   \par
   In \cite{cg1} Ciccoli and Gavarini extended this principle to quantum
formal homogeneous spaces.  Their result goes as follows.
                                                   \par
  Given a Poisson group  $ G \, $,  we consider pairs  $ \big(H,G/H\big) $
where  $ H $  is a coisotropic subgroup of  $ G $  and  $ G\big/H $  is
the corresponding homogeneous space.  At a local level, i.e.~in the setup
of formal geometry, any such pair can be described in algebraic terms by
any one of the following four objects:
  $$  U(\fh) \; ,  \quad  U(\fg) \fh \; ,  \quad  I(H) \; ,
\quad  \Ogh   \eqno (4.1)  $$
where hereafter  $ \, \fh := \text{\it Lie}\,(H) \, $  and  $ \, \fg :=
\text{\it Lie}\,(G) \, $,  and the notation is standard, but for  $ \Ogh $
which here denotes the algebra of regular functions on the formal algebraic variety  $ G\big/H \, $.  The main result in  \cite{cg1}  is a four-fold functorial recipe which, from a quantization of each object in (4.1), constructs a quantization of one of the four object of the similar quadruple which describes the ``dual'' pair  $ \Big( H^\perp, \, G^* \big/ H^\perp \Big) $,  where ``dual'' refers to Poisson duality.
                                                   \par
   If we try to do the same at a global level  (cf.~\cite{cg2}),  i.e.~not
restricting ourselves to the framework of formal geometry,
%
%
 then something changes when handling
the algebra  $ \Ogh \, $.  Namely, the latter is meaningful   --- that
is, it permits to get back the pair  $ \big( H, G/H \big) $
---   only if
$ G\big/H $  is a quasi-affine variety.  This is the case, in particular,
if  $ G\big/H $  is affine, and instead it is not if  $ \, G\big/H $  is projective.  Therefore, in the latter case one describes  $ G\big/H $  taking, instead of the algebra of regular functions, the algebra of (algebraic) sections of a line bundle on  $ G\big/H $  realizing an embedding in a projective space, i.e.~what is denoted  $ \Ogh $
in  \S \ref{class-set}.  Once this
%
%
 is settled, one can
consider a quantization  $ \Oqgh $  and try to cook up a
suitable analogue of the recipe of  \cite{cg1}, \cite{cg2}.

\smallskip

  With this program in mind, we want to associate to any
quantum homogeneous  $ G $--space  $ \Oqgh $   --- in the sense of
Definition \ref{q-prhspace}  ---   a (local) quantization  $ U_q \big(
\fh^\perp \big) $  of the dual  $ G^* $--space  (actually, of the dual
coisotropic subgroup), right in the spirit of the QDP.

\bigskip

   {\bf Warning:}  In order to make our statements simpler, from now on  {\sl we take as ground ring the local ring}
  $$  \bk_q'  \; := \;  \big( \bk_q \big)_{(q-1)}  \; = \;
\text{\sl localization of\/  $ \bk_q $  at the ideal generated
by  $ (q \! - \! 1) \, $.}  $$
Therefore,  {\sl hereafter we shall tacitly extend scalars
from\/  $ \bk_q $  to\/  $ \bk'_q $  for all\/  $ \bk_q $--modules  and
$ \bk_q $--algebras  we have considered so far},  with no further mention.

\bigskip

   Let  $ G $  be an affine algebraic group and  $ H $  a closed
coisotropic parabolic subgroup, i.e.~$ G\big/H $  is a projective
homogeneous space.  Let  $ \cO_q(G) \, $,  $ \cO_q(H) $  and
$ \Oqgh $  be quantum deformations of  $ \cO(G) \, $,  $ \cO(H) $
and  $ \Ogh $  as defined in Section 3.  In particular,  $ \Oqgh $
is built out of a specific quantum section  $ d $  on  $ G\big/H \, $.
Also,  $ \pi(d) = \overline{d} \in \cO_q(H) \, $  is (non-zero)
group-like, hence  $ \, \ep(d) = \ep\big(\,\overline{d}\,\big)
= 1 \, $,  and  $ d $  specializes to  $ \; d\big|_{q=1} =
t \in \Ogh \subseteq \cO(G) \, $.

\medskip

   In addition,  {\sl we make the following assumption:}

 \vskip7pt

   \centerline{\it  $ d \, $  is an Ore element in the algebra
$ \cO_q(G) \, $.}

 \vskip7pt

 This property will allow us to enlarge the algebras  $ \cO_q(G) $
and  $ \Oqgh $  by the formal inverse  $ d^{-1} \, $.  Geometrically,
it corresponds to ask to have   --- besides a quantization of  $ \Ogh $
---   a quantization of a Zariski neighbourhood of the identity element;
more precisely, it means that we have a quantization of the function
algebra  $ \cO \big( X_t \big) $  of the affine variety  $ \, X_t \, $,
the complement in  $ G\big/H $  of the divisor defined by the function
$ \, t = d \big|_{q=1} \; $.  This property is satisfied in the examples of the Grassmannian and the flag varieties  (cf.~\S \ref{Q-Grass}), with a suitable choice of  $ d \, $.

\medskip

   Let us define
  $$  \Oqghloc  \; := \;  \Oqgh\big[d^{-1}\big]_{proj}
\; \subseteq \;  \cO_q(G)\big[d^{-1}\big]  $$
where the localization is a  \textit{projective localization},
i.e.~we take the elements of degree zero in the ring  $ \,
\Oqgh\big[d^{-1}\big] \, $,  where  $ d^{-1} $  is given
degree  $ -1 \, $.

\subsection{The QDP functor} \label{qdp-functor}

   We now recall the definition of the functor  $ \, \cO_q(G) \mapsto
{\cO_q(G)}^\vee \, $,  which sends quantized function algebras (of
Poisson groups) to quantized universal enveloping algebras (of
Lie bialgebras).  More precisely,  $ {\cO_q(G)}^\vee $  is a
quantization of  $ U\big(\fg^*\big) $,   where  $ \fg^* $  is
the Lie bialgebra dual to  $ \fg \, $.  For more details and
proofs, we refer the reader to  \cite{ga2},  Theorem 2.2 and
Theorem 4.7.
 \vskip4pt
   {\sl Remark:}  the overall assumption in  \cite{ga2}  for  $ G $  is to be connected; nevertheless, this condition is not needed in the proof of Theorem 4.7 therein.

\medskip

\begin{definition}  \label{OqGvee}
  \, Let  $ \, J_G := \text{\it Ker}\, \big(\, \ep \, \colon \cO_q(G)
\lra \bk'_q \,\big) \, $  be the augmentation ideal of  $ \cO_q(G) \, $.
Also, let  $ \; I_G:= J_G + (q-1) \, \cO_q(G) \; $.  We define
  $$  {\cO_q(G)}^\vee  \, := \,  {\textstyle \sum\limits_{n \geq 0}}
\, {(q-1)}^{-n} I_G^n  \, = \,  {\textstyle \sum\limits_{n \geq 0}}
{\big( {(q-1)}^{-1} I_G \, \big)}^n  = \,
 {\textstyle \bigcup\limits_{n \geq 0}}
{\big( {(q-1)}^{-1} I_G \, \big)}^n  \;\; .  $$
   \indent   This is a well defined  $ \bk'_q $--subalgebra  of
$ \; \bk(q) \otimes_{\bk'_q} \cO_q(G) \; $.  Notice also that
  $$  {\cO_q(G)}^\vee  \; = \;  {\textstyle \sum_{n \geq 0}} \,
{(q-1)}^{-n} J_G^n  \; = \;  {\textstyle \sum_{n \geq 0}}
\, {\big( {(q-1)}^{-1} J_G \, \big)}^n \;\;\; .  $$
\end{definition}

\medskip

   The results of  \cite{ga2}   --- in particular,  Theorem 4.7 therein ---  tell us that  $ {\cO_q(G)}^\vee $  is a quantization of  $ U\big(\fg^*\big) \, $,  that is  $ \; {\cO_q(G)}^\vee \Big/ (q-1) \, {\cO_q(G)}^\vee \cong \, U\big(\fg^*\big) \; $  as co-Poisson Hopf algebras.  Our idea is to ``restrict to  $ \cO_q(G/H) \, $'',
somehow, the definition of  $ {\cO_q(G)}^\vee $,  so as to
define  $ {\cO_q(G/H)}^\vee $.  To begin with, let
  $$  J_{G/H}  \; := \;  \text{\it Ker}\,\big(\ep|_{\cO_q(G/H)}\big)
\quad .  $$
   \indent   Notice that  $ \ep $  extends uniquely to
$ \, \cO_q(G)[d^{-1}] \, $,  \, so we can define also
  $$  J_{G/H}^{\,loc}  \; := \;
\text{\it Ker}\,\big(\ep|_{\cO_q^{\,loc}(G/H)}\big)  \quad .  $$

\medskip

\begin{definition} \label{cOqghvee}  We define
  $$  \cO_q(G/H)^\vee  \; := \;  {\textstyle \sum_{n \geq 0}} \,
{(q-1)}^{-n} \big(J_{G/H}^{\,loc}\big)^n  \; = \;  {\textstyle
\sum_{n \geq 0}} \, {\big( {(q-1)}^{-1} J_{G/H}^{\,loc} \,
\big)}^n \;\;\; ,  $$
the unital  $ \bk'_q $--subalgebra  of  $ \, \bk(q) \otimes_{\bk'_q}
\Oqghloc \, $  generated by  $ \, {(q\!-\!1)}^{-1} J_{G/H}^{\text{\it loc}} \, $,  or   --- what amounts to be the same ---   by  $ \, {(q-1)}^{-1} I_{G/H}^{\text{\it loc}} \, $,  where by definition we set  $ \, I_{G/H}^{loc} := J_{G/H}^{\text{\it loc}} + (q-1) \, \Oqghloc \, $.
\end{definition}

\medskip

   Indeed, one can check that the construction  $ \, \Oqgh \mapsto
\Oqgh^\vee $  is functorial, in a natural sense (see  Remark \ref{functoriality}  later on).

\bigskip

   In order to study  $ \Oqgh^\vee $,  we need a rather explicit
description of it.  In turn, this requires a description of
$ \cO_q(G)^\vee $,  which we take from  \cite{ga2}.

\smallskip

   Let  $ J_1 $  be the augmentation ideal of  $ \cO(G) \, $,  namely
  $$  J_1  \; := \;  \text{\it Ker}\,\big(\,\epsilon \, \colon
\cO(G) \longrightarrow \bk \,\big)  \; = \;  J_G \! \mod (q \! - \! 1)
\, \cO_q(G)  $$
so that  $ \, J_1 \Big/ J_1^{\,2} = \fg^* \, $,  \, the cotangent Lie
bialgebra of  $ G \, $.  Let  $ \, \{y_1,\dots,y_n\} \, $  be a subset
of  $J_1$  whose image in the local ring of  $ G $  at  $ e \, $  (the
unit element of the group  $ G \, $)  is a local system of parameters;
in particular,  $ \, n = \text{\it dim}\,(G) \, $.  Define  $ \, \{j_1,
\dots, j_n\} \, $  as a pull-back of  $ \, \{y_1, \dots, y_n\} \, $
to  $ J_G \, $.

\medskip

\begin{theorem} \, {\sl (see  \cite{ga2},  Lemma 4.1)}
                                          \par
\noindent
 \; (a)  The set of ordered monomials  $ \, \Big\{ j^{\,\underline{e}} := \prod_{s=1}^n j_s^{\,e_s} \,\Big|\, \underline{e} := (e_1,\dots,e_n) \in \N^{\,n} \Big\} \, $  is a  $ \bk'_q $--pseudobasis  (or  {\sl topological\/}  basis) of  $ \, \widetilde{\cO_q}(G) \, $, the  $ I_G $--adic  completion
of  $ \, \cO_q(G) \, $.  In other words, each element of
$ \, \widetilde{\cO_q}(G) $  has a unique expansion as a
formal infinite  $ \bk'_q $--linear  combination of the
$ j^{\,\underline{e}} $'s.  In particular, $ \widetilde{\cO_q}(G) $
is generated   --- as a  {\sl topological}  $ \bk'_q $--algebra  ---
by  $ \, \{j_1, \dots, j_n\} \, $.
                                          \par
\noindent
 \; (b)  The  $ (q\!-\!1) $--adic  completion  $ \widehat{\cO_q}(G)^\vee $
of  $ \, {\cO_q(G)}^\vee $  admits the set of ordered monomials  $ \, \Big\{
{(q\!-\!1)}^{-(e_1 + \cdots + e_n)} \prod_{s=1}^n j_s^{\,e_s} \,\Big|\,
(e_1, \dots, e_n) \in \N^{\,n} \Big\} \, $  as a  $ \bk'_q $--pseudobasis.
In particular,  $ \widehat{\cO_q}(G)^\vee $  is generated   --- as a  {\sl
topological}  $ \bk'_q $--algebra  ---   by  $ \; \big\{\, j^{\,\vee}_s \!
:= {(q\!-\!1)}^{-1} j_s \,\big|\, s = 1, \dots, n \,\big\} \; $.
\end{theorem}

\medskip

   The description of  $ \Oqgh^\vee $  goes much along the
same lines.

\medskip

   Recall that  $ \, t \in \Ogh \, $  is the specialization of
$ \, d \in \Oqgh \, $.  We then consider  $ X_t \, $,  the open
subvariety of  $ G\big/H $  where  $ \, t \neq 0 \, $.  On this
variety, choose functions  $ \, l_1, \dots, l_{n-h} \, $   --- where
$ \, h = \text{\it dim}\,(H) \, $  ---   such that the set
  $$  \big\{\, {x}_s := {l}_s \!\! \mod (q\!-\!1) \, \cO_q\big(G/H\big)
\;\big|\; s = 1, \dots, n\!-\!h \,\big\}  $$
yields, in the localization of  $ \cO \big( X_t \big) $  at  $ \, eH \in
X_t \; \big(\! \subseteq G/H \big) \, $,  a local system of parameters
at  $ \, eH \, $.

\begin{theorem}  \quad
%
%
                                          \par
 \; (a)  The set of ordered monomials
  $$  \bigg\{\; {\textstyle \prod\limits_{s=1}^{n-h}} \, l_s^{\,e_s}
\;\bigg|\; (e_1,\dots,e_{n-h}) \in \N^{\,n-h} \,\bigg\}  $$
is a  $ \bk'_q $--pseudobasis  of  $ \widetilde{\cO_q^{\text{\it loc}}}
\big(G\big/H\big) \, $,  the latter being the  $ I_{G/H}^{\,loc} $--adic
completion of  $ \, \Oqghloc \, $,  where  $ \, I_{G/H}^{\,\text{\it loc}}
:= J_{G/H}^{\,\text{\it loc}} + (q\!-\!1) \, \Oqghloc \; $.
                                                    \par
   In particular,  $ \widetilde{\cO_q^{\text{\it loc}}}\big(G\big/H\big) $
is  ({\sl topologically})   generated by $ \, \{l_1,\dots,l_{n-h}\} \, $.
 \vskip5pt
 \; (b)  The  $ (q-1) $--adic  completion  $ \widehat{\cO_q} \big( G\big/H
\big)^\vee $  of  $ \, \Oqgh^\vee $  admits as a  $ \bk'_q $--pseudobasis
the set of ordered monomials
  $$  \bigg\{\; {(q-1)}^{-(e_1 + \cdots + e_{n-h})}
{\textstyle \prod\limits_{s=1}^{n-h}} \, l_s^{\,e_s} \;\bigg|\;
(e_1,\dots,e_{n-h}) \in \N^{\,n-h} \,\bigg\}  \quad .  $$
   \indent   In particular,  $ \widehat{\cO_q}\big(G/H\big)^\vee $
is  ({\sl topologically})  generated by the set
  $$  \Big\{\, l^{\;\vee}_s := {(q\!-\!1)}^{-1} l_s
\;\Big|\; s = 1, \dots, n-h \,\Big\}  \quad .  $$
\end{theorem}

\begin{proof}  The argument follows the proof of Theorem 4.9 in
\cite{ga3}.  In this theorem is treated the general setting of a
quantization  $ \cO_q(V) $  of any Poisson affine variety  $ V $
with a distinguished point on it, given by a character  $ \chi $
of  $ \cO_q(V) \, $,  the kernel of  $ \chi $  playing the role of
$ J_G $  above.  Here we apply all this to  $ \, V = X_t \; \big(\!
\subseteq {G/H} \big) \, $  with  $ \, \ep\big|_{\cO(X_t) }\, $  as
character on  $ \, \cO(X_t) = \Oqghloc \; $.
\end{proof}

\smallskip

   The next lemma plays a crucial role in the construction of the QDP.

\medskip

\begin{lemma}
 The quantum section  $ \, d \in \Oqgh \, $  enjoys the following
properties:
                                                 \par
   (a) \,  $ d $  is invertible in  $ \widetilde{\cO_q^{\text{\it
   loc}}}
\big(G\big/H\big) \, $,  with  $ \; d^{-1} = \sum\limits_{k=0}^{+\infty}
\, {\big( 1 - d \,\big)}^k \; $;
                                                 \par
   (b) \,  $ d $  is invertible in  $ \widehat{\cO_q}\big(G\big/H\big)^\vee
\, $,  with  $ \; d^{-1} = \sum\limits_{k=0}^{+\infty} \, {(q\!-\!1)}^k
\bigg(\! \displaystyle{\frac{\,(1-d)\,}{\,(q-1)}} \bigg)^k \, $.
\end{lemma}

\begin{proof}
Observe that  $ \, \ep(d) = 1 \, $  implies that
  $$  d  \, = \,  1 - \big( 1 - d \,\big) \, \in \,
\big( 1 + J_{G/H}^{\,\text{\it loc}} \big) \,
\subseteq \, \big( 1 + I_{G/H}^{\,loc} \big)  \quad ;  $$
this gives the invertibility in  $ \widetilde{\cO_q^{\text{\it
loc}}} \big(G\big/H\big) \, $.
                                                   \par
   Similarly  since  $ \, J_{G/H}^{\,\text{\it loc}} \subseteq
(q\!-\!1) \, \cO_q\big(G\big/H\big)^\vee \,$,  the identity
  $$  d  \, = \,  1 + \big( d - 1 \big) \, \in \, \big( 1 +
J_{G/H}^{\,\text{\it loc}} \big) \, \subseteq \, \Big( 1 +
(q\!-\!1) \, \cO_q\big(G\big/H\big)^\vee \Big)  $$
also ensures that  $ d $ is invertible in  $ \widehat{\cO_q}
\big(G\big/H\big)^\vee \, $.
                                                   \par
   In both cases, the explicit formula for  $ d^{-1} $  follows by
taking the limit of the geometric series, namely  $ \; {(1-x)}^{-1}
= \sum_{k=0}^{+\infty} \, x^k \, $,  \, applied to  $ \; x =
1 - d \; $.
\end{proof}

\vskip1pt

\begin{proposition}  \label{emb-compl}
   There are natural embeddings
$$
 \widetilde{\cO_q}\big(G\big/H\big)
 \lhook\joinrel\relbar\joinrel\relbar\joinrel\longrightarrow
 \, \widetilde{\cO_q}(G) \;\; ,   \quad \qquad
 \widehat{\cO_q}\big(G\big/H\big)
 \lhook\joinrel\relbar\joinrel\relbar\joinrel\longrightarrow
 \, \widehat{\cO_q}(G)
$$
which both are extensions of the embedding  $ \, \Oqgh
\lhook\joinrel\relbar\joinrel\longrightarrow \, \cO_q(G)
\, $.  Moreover, via these embeddings the pseudobases for the
(topological) algebras on  $ G\big/H $  identify with subsets
of the corresponding ones for the (topological) algebras on  $ G \, $.
\end{proposition}

\begin{proof}
   By construction, we have  $ \, l_k = j_s \big/ d^{c_s} \, $  for
some  $ \, j_s \in J_{G/H} \, $,  $ \, c_s \in \N \, $  ($ \, s = 1,
\dots, n-h \, $).  Since  $ d $  is invertible (in both cases), the
previous analysis tells us that we can replace the  $ l_k $'s  with
the  $ j_k $'s  ($ \, k = 1, \dots, n-h \, $)  in the descriptions of
$ \widetilde{\cO_q}\big( G\big/H \big)^\vee $  and  $ \widehat{\cO_q}
\big(G\big/H\big)^\vee $  given above (i.e.~in the
$ \bk'_q $--pseudobases  and as topological generators).
But then, as
  $$  \big\{\, j_s \,\big|\, s = 1, \dots, n\!-\!h \,\big\} \, \subseteq
\, J_{G/H} \, \subseteq \, \Oqgh \, \subseteq \, \cO_q(G)  $$
we can always complete  $ \, \big\{\, j_s \,\big|\, s = 1, \dots, n\!-\!h
\,\big\} \,$  to a set  $ \, \big\{\, j_r \,\big|\, r = 1, \dots, n \,\big\}
\, $  such that  $ \, \big\{\, y_r := j_r \mod (q\!-\!1) \, \cO_q(G)
\;\big|\; r = 1, \dots, n \,\big\} \, $  yield, in the localization
of  $ \cO(G) $  at  $ \, e \in G \, $,  a local system of parameters at
$ \, e \, $.  Thus using the latter we can describe  $ \widetilde{\cO_q}
(G) $  and  $ \widehat{\cO_q}(G) $  as explained above.
\end{proof}

\medskip

   From now on we shall use these embeddings
%
%
 to identify  $ \widetilde{\cO_q}\big(G\big/H\big) $
and  $ \widehat{\cO_q}\big(G\big/H\big) $  with a subalgebra of
$ \widetilde{\cO_q}(G) $  and of  $ \widehat{\cO_q}(G)  $  respectively.

\vskip15pt

\begin{lemma}  \label{coqgv-strict}
 $ \quad \displaystyle{ \Oqgh^\vee \, {\textstyle \bigcap}
\; (q-1) \, \cO_q{(G)}^\vee \; = \; (q-1) \, \Oqgh^\vee} $
%
%
\end{lemma}

\begin{proof} Let us choose a subset  $ \{j_1,\dots,j_n\} $  in
$ J_{G/H}^{\,loc} $  as explained above for the description of
$ \Oqgh^\vee $.  Then, mapping  $ \cO_q(G)^\vee $  and
$ \Oqgh^\vee $  into their  $ (q\!-\!1) $--adic
completions, and exploiting the descriptions of the latter
via pseudobases given above, we easily get the claim.
\end{proof}

\medskip

   Next result is that  $ \Oqgh^\vee $  is a quantization of
$ U\big(\fh^\perp\big) \, $:

\medskip

\begin{theorem} \label{quant1}
  $ \Oqgh^{\vee} $  is a quantization of  $ U\big(\fh^\perp\big) $
as a  $ \bk $--algebra   --- subalgebra of  $ \, U(\fg^*) $  ---
where  $ \, \fh = \text{\it Lie}\,(H) \, $,  $ \, \fg = \text{\it
Lie}\,(G) \, $.
\end{theorem}

\begin{proof}  By assumption,  $ H $  is coisotropic in  $ G \, $.
Therefore,  $ \, \fh = Lie(H) \, $  is a Lie coideal (and subalgebra)
of $ \, \fg = Lie(G) \, $,  and  $ \fh^\perp $  is a Lie subalgebra
(and coideal) of  $ \fg^* \, $  (see  Proposition \ref{coiso-inf}).
Thus the claim does make sense.
 \vskip3pt
   In order to prove the statement, we proceed much like in the proof
of the fact that  $ \; \cO_q(G)^\vee \Big/ (q-1) \, \cO_q(G)^\vee \cong
U(\fg^*) \; $   --- cf.~\cite{ga2},  Theorem 4.7.  The arguments being
the same, we briefly recall them.
 \vskip5pt
   Fix again a special subset  $ \, \{ j_1, \dots, j_n \} \, $  of
$ J_G $  as we did in the proof of  Proposition \ref{emb-compl},  in
particular with  $ \, j_1, \dots, j_{n-h} \in J_{G/H} \, $.  Also, set
notation:
  $$  \displaylines{
   \cO_1(G)^\vee \, := \; \cO_q(G)^\vee \Big/ (q-1) \, \cO_q(G)^\vee
\;\, ,  \qquad  J_G^\vee \, := \; {(q\!-\!1)}^{-1} J_G \; \subseteq
\; \cO_q(G)^\vee  \cr
   j^{\,\vee} \, := \; {(q-1)}^{-1} j \quad \forall \;\; j \in J_G \;\, ,
\qquad  \mathfrak{t} \; := \; J_G^\vee \mod (q-1) \, \cO_q(G)^\vee
\;\; . }  $$
Taking into account that the specializations at  $ \, q = 1 \, $  of any
$ \bk'_q $--module  and of its  $ (q-1) $--adic  completion are the
same, the above discussion gives that
  $$  \bigg\{\; {\textstyle \prod\limits_{s=1}^n}
\big(j_s^{\,\vee}\big)^{\,e_s} \!\! \mod (q-1) \, \cO_q(G)^\vee
\,\;\bigg|\;\, (e_1,\dots,e_n) \in \N^{\,n} \,\bigg\}  $$
is a  $ \bk $--basis  of  $ \, \cO_1(G)^\vee \, $.  Similarly,
$ \, \big\{ j_1^{\,\vee}, \dots, j_n^{\,\vee} \big\} \, $
is a  $ \bk $--basis  of  $ \, \mathfrak{t} \, $.
                                              \par
   Now,  $ \, j_\mu \, j_\nu - j_\nu \, j_\mu \in (q \! - \! 1) \,
J_G \, $  (for  $ \, \mu, \nu \in \{1,\dots,n\} \, $)  implies that
  $$  j_\mu \, j_\nu - j_\nu \, j_\mu \; = \; (q - 1) \, {\textstyle
\sum_{s=1}^n} \, c_s \, j_s \, + \, {(q - 1)}^2 \, \gamma_1 \, +
\, (q - 1) \, \gamma_2  $$
for some  $ \, c_s \in \bk'_q \, $,  $ \, \gamma_1 \in J_G \, $  and
$ \, \gamma_2 \in J_G^{\,2} \, $.  Therefore
  $$  \displaylines{
   \quad   \big[ j_\mu^\vee, j_\nu^\vee \,\big] := j_\mu^\vee \,
j_\nu^\vee - j_\nu^\vee \, j_\mu^\vee = {\textstyle \sum_{s=1}^n}
\, c_s \, j_s^\vee + \gamma_1 + {(q\!-\!1)} \, \gamma_2^\vee
\equiv   \hfill  \cr
  \hfill   \equiv {\textstyle \sum_{s=1}^n} \, c_s \, j_s^\vee \mod
\, (q\!-\!1) \, \cO_q(G)^\vee  \quad }  $$
where  $ \, \gamma_2^\vee := {(q-1)}^{-2} \gamma_2 \in {(q-1)}^{-2}
{\big( J_G^\vee \big)}^2 \subseteq \cO_q(G)^\vee \, $;  thus
$ \mathfrak{t} $  is a Lie subalgebra of  $ \cO_1(G)^\vee \, $.
But then we have  $ \, \cO_1(G)^\vee \cong U(\mathfrak{t}) \, $
as Hopf algebras, by the above description of  $ \cO_1(G)^\vee $
and PBW theorem.
 \vskip5pt
   Next, the specialization map  $ \; p^\vee \colon \, \cO_q(G)^\vee
\relbar\joinrel\relbar\joinrel\twoheadrightarrow \cO_1(G)^\vee =
U(\mathfrak{t}) \; $  restricts to  $ \; \eta \, \colon \, J_G^\vee
\relbar\joinrel\relbar\joinrel\relbar\joinrel\twoheadrightarrow
\mathfrak{t} \, := J_G^\vee \! \mod (q\!-\!1) \, \cO_q(G)^\vee
\, $.  Moreover, multiplication by  $ {(q\!-\!1)}^{-1} $
yields a  $ \bk'_q $--module  isomorphism  $ \, \mu
\, \colon \, J_G \, {\buildrel \cong \over
{\lhook\joinrel\relbar\joinrel\relbar\joinrel\twoheadrightarrow}}
\, J_G^\vee \; $.  Consider the natural projection map  $ \; \rho
\, \colon \, J_1 \relbar\joinrel\twoheadrightarrow J_1 \big/
J_1^{\,2} = \fg^* \, $,  and let  $ \; \nu \, \colon \, \fg^*
\lhook\joinrel\longrightarrow J_1 \, $  be a section of  $ \rho
\, $.  The specialization map  $ \; p \, \colon \, \cO_q(G)
\relbar\joinrel\twoheadrightarrow \cO(G) \; $  restricts to
$ \; p' \colon \, J_G \relbar\joinrel\twoheadrightarrow J_1
\; $,  and we fix a section  $ \, \gamma \, \colon \, J_1 \!
\lhook\joinrel\relbar\joinrel\rightarrow J_G \, $  of  $ p' \, $.
Then the composition map  $ \, \sigma := \eta \circ \mu \circ
\gamma \circ \nu \, \colon \, \fg^* \longrightarrow \mathfrak{t}
\, $  is a well-defined Lie bialgebra isomorphism, independent
of the choice of  $ \nu $  and  $ \gamma \, $.
 \vskip5pt
   So far we did not exploit our special choice of the subset  $ \,
\{ j_1, \dots, j_n \} \, $:  \, we do it now to prove that  $ \,
\mathfrak{t} = \mathfrak{h}^\perp \; $.  In fact, the analysis above
to prove that  $ \, \sigma \, \colon \, \fg^* \cong \mathfrak{t} \, $
shows also that the unital subalgebra
  $$  \cO_1\big(G\big/H\big)^{\vee} \; := \,\; \Oqgh^{\vee} \! \mod
(q\!-\!1) \, \cO_q(G)^{\vee}  $$
of  $ U\big(\fg^*\big) $  is generated by  $ \eta \big( \big\{
j_1^{\,\vee}, \dots, j_{n-h}^{\,\vee} \big\} \big) $,  and
  $$  \displaylines{
   \eta \big( \big\{ j_1^{\,\vee}, \dots, j_{n-h}^{\,\vee} \big\} \big)
\, = \, (\eta \circ \mu) \big( \{ j_1 , \dots, j_{n-h} \} \big) \, = \,
(\eta \circ \mu \circ \gamma) \big( \{ y_1 , \dots, y_{n-h} \} \big) \,
=   \hfill  \cr
   \hfill   = \, (\eta \circ \mu \circ \gamma \circ \nu) \big( \big\{
\overline{y}_1 , \dots, \overline{y}_{n-h} \} \big) \, = \, \sigma
\big( \big\{ \overline{y}_1 , \dots, \overline{y}_{n-h} \} \big) }  $$
 \eject
\noindent
 where  $ \, \overline{y}_s := y_s \! \mod J_1^{\,2} \, $  ($ \, s = 1,
\dots, n\!-\!h \, $).  Therefore  $ \cO_1 \big( G\big/H \big)^{\vee} $
is the subalgebra of  $ U\big(\fg^*\big) $  generated by the
$ \bk $--span  of  $ \big\{ \overline{y}_1 , \dots,
\overline{y}_{n-h} \} \, $.
                                          \par
   Finally, the  $ \bk $--span  of  $ \big\{ \overline{y}_1 ,
\dots, \overline{y}_{n-h} \} $  coincides with the subspace
$ \mathfrak{h}^\perp $  of  $ \fg^* $.  Indeed, as  $ \Ogh $
is the algebra of semi-invariant functions on  $ G \, $,  every
$ y_s $  is a  ($ H $--)semi-invariant  function on  $ G \, $:  but
it also vanishes at  $ e \in H $,  hence by  $ H $--semi-invariance
it vanishes on all of  $ H \, $.  When mapping  $ y_s $  to
$ \, \overline{y}_s \in J_1 \, $,  then, it is mapped into
$ \mathfrak{h}^\perp $.  Thus the whole  $ \bk $--span  of
$ \big\{ \overline{y}_1 , \dots, \overline{y}_{n-h} \} $  is
contained in  $ \mathfrak{h}^\perp $,  hence coincides with
it by dimension equality.
                                          \par
   The outcome is that  $ \cO_1\big(G\big/H\big)^{\vee} $
is the subalgebra of  $ U\big(\fg^*\big) $  generated by
$ \mathfrak{h}^\perp $,  which is a Lie subalgebra of
$ \fg^* $,  so  $ \; \cO_1\big(G\big/H\big)^{\vee}
= \, U\big(\mathfrak{h}^\perp\big) \; $.
\end{proof}

\medskip

   We now wish to explore the nature of  $ \widehat{\cO_q}
\big(G/H\big)^\vee $  as a ``quantum homogeneous space''.  We
start with an important observation on the extensions of the
comultiplication  $ \Delta $  in  $ \cO_q(G) $  to the new
algebras we have defined.

\medskip

\begin{remark}
  Let $ \; \Delta : \cO_q(G) \lra \cO_q(G) \otimes \cO_q(G) \; $  be
the comultiplication in  $ \cO_q(G) \, $.  Then  $ \Delta $  can be
uniquely (and canonically) extended to a coassociative morphism of
topological algebras
  $$  \widetilde{\Delta} : \, \widetilde{\cO_q}(G) \relbar\joinrel\lra
\widetilde{\cO_q}(G) \,\widetilde{\otimes}\, \widetilde{\cO_q}(G)  $$
where again  $ \widetilde{\cO_q}(G) $  is the  $ I_G $--adic
completion of  $ \cO_q(G) $,  and  $ \, \widetilde{\otimes}\, $  is
the  $ {I_G}_\otimes $--adic  completion of  $ \, \cO_q(G) \otimes
\cO_q(G) \, $,  with
  $ \; {I_G}_\otimes \, := \, I_G \otimes \cO_q(G) \, + \,
\cO_q(G) \otimes I_G \; $.
%
%
Even more, such a  $ \widetilde{\Delta} $  actually restricts to a
coassociative algebra morphism (we use the same symbol to denote it):
  $$  \widetilde{\Delta} : \, \cO_q(G)\big[d^{-1}\big] \relbar\joinrel\lra
\cO_q(G)\big[d^{-1}\big] \,\widetilde{\otimes}\; \cO_q(G)\big[d^{-1}\big]
\quad .  $$
In fact, as  $ d $  is a quantum section we have (see  Definition
\ref{q-sect})
  $$  \Delta(d)  \; = \;  d \otimes d \, + \, {\textstyle \sum_i}
\, h_i \otimes k_i \;\; ,  \quad  \text{for some \ \ } h_i \in
\cO_q(G) \, , \; k_i \in I_q(H) \; .   \eqno (4.2)  $$
Since  $ d $  is Ore, we can re-write  $ \, \tDelta(d) = \Delta(d)
\big( d^{-1} \otimes d^{-1} \big) (d \otimes d) \, $,  \, which reads
  $$  \tDelta(d) \; = \; \big( 1 \otimes 1 \, + \, {\textstyle \sum_i}
\, h_i \, d^{-1} \otimes k_i \, d^{-1} \big) (d \otimes d) \;\; .  $$
This in turn implies
  $$  \hskip-5pt \begin{array}{cl}
   \tDelta\big(d^{-1}\big)  \; = \;  {\tDelta(d\,)}^{-1}  \, =
\;  {(d \otimes d\,)}^{-1}  &  \hskip-4pt \big( 1 \otimes 1 +
{\textstyle \sum_i} \, h_i \, d^{-1} \otimes k_i \, d^{-1} \big)^{-1}
\; =   \hfill  \\
   \hfill   = \;  \big( d^{-1} \otimes d^{-1} \big)  &
\hskip-4pt {\textstyle \sum\limits_{n=0}^{+\infty}} \,
(-1)^n \Big( {\textstyle \sum_i} \, h_i \, d^{-1}
\otimes k_i \, d^{-1} \Big)^n
\end{array}   \hskip13pt \hfill (4.3)  $$
where the bottom term does belong to  $ \cO_q(G)\big[d^{-1}\big]
\,\widetilde{\otimes}\, \cO_q(G)\big[d^{-1}\big] $,  as expected,
because  $ \, k_i \in I_q(H) \subset I_G \, $  (for every  $ i \, $),
hence the last formal series above is convergent in the
$ {I_G}_\otimes $--adic  topology.
\end{remark}

\medskip

   Let us now turn our attention to the algebra  $ {\cO_q}{(G)}^\vee $
and its  $ (q-1) $--adic  completion  $ \widehat{\cO_q}{(G)}^\vee \, $.
   By construction (cf.~\cite{ga2}),  the coproduct of  $ {\cO_q(G)}^\vee $,
hence of  $ \widehat{\cO_q}{(G)}^\vee $  too, is induced by  the coproduct
$ \Delta $  of  $ \cO_q(G) \, $.  Note that the coproduct  $ \widehat{\Delta} $
of  $ \widehat{\cO_q}{(G)}^\vee $  takes values in the  {\sl topological\/}
tensor product  $ \, \widehat{\cO_q}{(G)}^\vee \widehat{\otimes}\,
\widehat{\cO_q}{(G)}^\vee \, $,  which by definition is the  $ (q \!
- \! 1) $--adic  completion of the  {\sl algebraic\/}  tensor product
$ \, \widehat{\cO_q}{(G)}^\vee \! \otimes \widehat{\cO_q}{(G)}^\vee \, $
--- and coincides, moreover, with the  $ (q \! - \! 1) $--adic  completion
of  $ \, \cO_q{(G)}^\vee \! \otimes \cO_q{(G)}^\vee \, $.

\vskip11pt

   We are ready to move another key step.

\medskip

\begin{proposition} \label{quant2}
 $ \, \widehat{\cO_q}\big(G/H\big)^\vee \, $
is a left coideal of  $ \; \widehat{\cO_q}{(G)}^\vee \, $.
\end{proposition}

{\it Proof.}
 We want to show that the coproduct  $ \widehat{\Delta} $  maps
$ \widehat{\cO_q}\big(G/H\big)^\vee $  into the topological tensor
product  $ \; \widehat{\cO_q}{(G)}^\vee \, \widehat{\otimes} \;
\widehat{\cO_q}\big(G/H\big)^\vee \, $.

\medskip

   We first observe that
  $$  \tDelta\big(d^{-1}\big) \, \in \,
\cO_q(G)\big[d^{-1}\big] \,\widetilde{\otimes}\,
\Oqgh\big[d^{-1}\big] \;\; .  $$
This is because  $ \Oqgh $  is a left coideal of  $ \cO_q(G) \, $
--- cf.~Proposition  \ref{Oqgh-coid},  Theorem \ref{Oqgh-graded}
---   hence we have that the elements  $ k_i $  occurring in formula
(4.2) can be taken to belong to  $ \cO_q\big(G/H\big) \, $.
                                              \par
   Even more precisely, as the  $ \cO_q(G) $--coaction  on  $ \Oqgh $
via  $ \Delta $  is graded (by  Theorem \ref{Oqgh-graded}{\it (c)\/}),
all the  $ k_i $'s  have degree 1, like  $ d $  itself.  Thus, the
series occurring in (4.3) in fact belongs to  $ \cO_q(G)\big[d^{-1}\big]
\,\widetilde{\otimes}\, \cO_q^{loc}\big(G/H\big) \, $.  To sum up,
  $$  \widetilde{\Delta}\big(d^{-1}\big)  =
\big( d^{-1} \otimes d^{-1} \big) \cdot \delta_\otimes \; ,  \quad
\text{with \ }  \delta_\otimes \in \cO_q(G)\big[d^{-1}\big]
\,\widetilde{\otimes}\, \cO_q^{loc}\big(G/H\big) \; .   \eqno
(4.4)
$$
   \indent   Since the coaction  $ \; \Delta \, \colon \, \Oqgh
\longrightarrow \cO_q(G) \otimes \Oqgh \; $  is grading-preserving
and product-preserving, the definitions of  $ \cO_q^{loc}\big(G/H\big) $  and  $ \cO_q^{loc}(G) $  and (4.4) together yield
  $$  \widetilde{\Delta}\Big(\cO_q^{loc}\big(G/H\big)\!\Big)
\; \subseteq \;  \cO_q(G)\big[d^{-1}\big] \,\widetilde{\otimes}\,
\cO_q^{loc}\big(G/H\big) \;\; .   \eqno (4.5)  $$

\smallskip

   Now, we described above the completions of the algebras  $ A $
and  $ A^\vee $   --- for  $ \, A \in \big\{ \cO_q\big(G/H\big), \cO_q(G)
\big\} \, $  ---   w.r.t.~the  $ I_G $--adic  or the  $ (q-1) $--adic
topology.  Using that, or an entirely similar analysis, we see
also that
  $$  \cO_q(G)\big[d^{-1}\big] \,\widetilde{\otimes}\;
\cO_q^{\text{\it loc}}\big(G/H\big)  \,\; \subseteq \;\,
\widehat{\cO_q}{(G)}^\vee \,\widehat{\otimes}\;
\widehat{\cO_q}\big(G/H\big)^\vee \;\; .  $$
In short, this is because  $ \, I_G \subseteq (q\!-\!1) \, \cO_q{(G)}^\vee
\, $  and  $ \, I_{G/H} \subseteq (q\!-\!1) \, \cO_q\big(G/H\big)^\vee
\, $.  Also, it is easily seen that
  $$  \widetilde{\Delta} \big( I_{G/H}^{\text{\it loc}} \big)
\,\; \subseteq \;\,  \cO_q(G)\big[d^{-1}\big]
\,\widetilde{\otimes}\, I_{G/H}^{\text{\it loc}} \, + \, I_G
\;\widetilde{\otimes}\, \cO_q^{\text{\it loc}}\big(G/H\big)  \quad
;
$$
this along with (4.5) immediately implies
  $$  \widehat{\Delta} \big( (q\!-\!1)^{-1} I_{G/H}^{\text{\it loc}} \big)
\,\; \subseteq \;\,  \cO_q(G)\big[d^{-1}\big]
\,\widehat{\otimes}\, (q\!-\!1)^{-1} I_{G/H}^{\text{\it loc}} \, +
\, (q\!-\!1)^{-1} I_G \;\widehat{\otimes}\, \cO_q^{\text{\it
loc}}\big(G/H\big)  $$
which in turn yields, by the very definition of  $ \cO_q{(G)}^\vee $
and  $ \Oqgh^\vee $,
  $$  \widehat{\Delta} \Big( \Oqgh^\vee \Big)  \,\; \subseteq \;\,
\cO_q{(G)}^\vee \,\widehat{\otimes}\; \Oqgh^\vee \;\; .  $$
Finally, taking  $ (q-1) $--adic  completions on both sides,
and also noting that  $ \, \cO_q{(G)}^\vee \,\widehat{\otimes}\;
\Oqgh^\vee = \, \widehat{\cO_q}{(G)}^\vee \,\widehat{\otimes}\;
\widehat{\cO_q}\big(G/H\big)^\vee \, $,  \, we get
  $$  \widehat{\Delta} \Big( \widehat{\cO_q}\big(G/H\big)^\vee \Big)
\,\; \subseteq \;\,  \cO_q{(G)}^\vee \,\widehat{\otimes}\;
\widehat{\cO_q}\big(G/H\big)^\vee  \quad .   \eqno \qed  $$

\medskip

   In the end, we get the main result of this section:

\medskip

\begin{theorem}  \label{maintheorem}
  $ \, \widehat{\cO_q}\big(G\big/H\big)^\vee \, $  is a quantization
of  $ \, U\big(\fh^\perp\big) $  as a subalgebra and left coideal
of  $ \, U(\fg^*) \, $.  In other words,  $ \widehat{\cO_q}
\big(G\big/H\big)^\vee $  is an  {\sl infinitesimal}  quantization
of the coisotropic subgroup  $ H^\perp $  of  $ \, G^* \, $.
\end{theorem}

\begin{proof}
 Just collect the previous results.  First we have
  $$  \widehat{\cO_q}\big(G\big/H\big)^\vee \,\;{\textstyle
\bigcap}\;\, (q-1) \, \widehat{\cO_q}{(G)}^\vee  \,\; = \;\,
(q-1) \, \widehat{\cO_q}\big(G\big/H\big)^\vee  $$
as an easy consequence of  Lemma \ref{coqgv-strict}.  Then, by
Theorem \ref{quant1}  and by the fact that  $ \; \widehat{\cO_q}
\big( G\big/H \big)^\vee{\Big|}_{q=1} = \, \Oqgh^\vee{\Big|}_{q=1}
\; $,  \, we have that the specialization of  $ \widehat{\cO_q}
\big(G\big/H\big)^\vee $  is  $ U\big(\fh^\perp\big) \, $.
Moreover,  Proposition  \ref{quant2}  proves that the subalgebra
$ \widehat{\cO_q}\big(G\big/H\big)^\vee $  of  $ \, \widehat{\cO_q}
{(G)}^\vee $  is also a left coideal.  Therefore,  $ \widehat{\cO_q}
\big(G\big/H\big)^\vee $  is an infinitesimal quantization of
$ H^\perp $,  in the standard sense.
\end{proof}

\medskip

\begin{remark}  \label{functoriality}
 The construction of  $ \, \widehat{\cO_q}\big(G\big/H\big)^\vee \, $
is functorial, in the following sense.  For a fixed  $ \cO_q(G) \, $,
every $ \Oqgh $  is uniquely characterized by the pair
$ \big( \pi_H , d_H \big) $  given by the projection  $ \; \pi_H \colon \cO_q(G) \relbar\joinrel\twoheadrightarrow \cO_q(H) \; $  and the quantum section  $ \, d_H \in \cO_q(G) \, $.  The natural notion of morphism among
such pairs, say  $ \, \big( \pi_{H'} , d_{H'} \big) \longrightarrow \big( \pi_{H''} , d_{H''} \big) \, $, can be cast into the form a Hopf algebra endomorphism  $ \phi $  of  $ \cO_q(G) $  such that  $ \, \phi \big( \text{\it Ker}\,(\pi_{H'}) \big) \subseteq \text{\it Ker}\,(\pi_{H''}) \, $   --- or  $ \, \phi \big( I_q(H') \big) \subseteq I_q(H'') \, $  ---   and
$ \, \phi(d_{H'}) = d_{H''} \, $.  Then, one defines  $ \phantom{\Big|} \! {(\ )}^\vee $  on morphisms by scalar extension followed by restriction; proving the functoriality is a matter of bookkeeping.  More in general,
one might decide not to fix  $ \cO_q(G) \, $,  nor even  $ G \, $.  Then morphisms  $ \, \phi \colon \cO_q(G') \longrightarrow  \cO_q(G'') \, $
take the place of the endomorphisms of (the single)  $ \cO_q(G) $  in
the recipe above, yet  $ {(\ )}^\vee $  is defined again on morphisms
via scalar extension and restriction   --- and one has to exploit the functoriality of  $ {(\ )}^\vee $  over quantum groups  $ \cO_q(G) \, $.
\end{remark}

\bigskip

\section{Examples: Quantum Grassmannians and quantum flag varieties}  \label{Q-Grass}

\medskip

In this section we want to examine in detail some examples of quantum
homogeneous spaces and apply the quantum duality principle recipe to
them. We start with the quantum Grassmannian.

\subsection{The quantum Grassmannian as quantum projective homogenous space}
 \label{class-case}

Let us recall the classical setting.

\medskip

Let  $ \, G = GL_n(\bk) \, $  and let
$ \, H = P \, $  a (standard) maximal parabolic subgroup, say
  $$  P \, = \, \left\{\,
   \begin{pmatrix}
      A  &  B  \\
      0  &  D
   \end{pmatrix}
\,\bigg|\; A \in GL_r(\Bbbk) \, , \, B \in \text{\it Mat}_{r,n-r}(\Bbbk) \, ,
\, D \in GL_{n-r}(\Bbbk)  \,\right\}  \quad .  $$

\medskip

   We want first to give a very ample line bundle on the homogeneous
space  $ G\big/ P $   --- the Grassmann variety ---   that realizes the
classical Pl{\"u}cker embedding into a projective space, following
the recipe in  \S \ref{proj-case}.
%
 \vskip7pt
   Let  $ \, I = (i_1, \dots, i_r) \, $,  $ \, 1 \leq i_1 < \dots < i_r
\leq n \, $.  Define
  $$  d^I \, :  \; g = \big( x_{ij} \big)  \; \mapsto
\;  d^I(g)  \, := \,  {\textstyle \sum_{\sigma \in
\,\mathcal{S}_r}} (-1)^{\ell(\sigma)} \, x_{i_1,\,\sigma(1)}
\cdots x_{i_r,\,\sigma(r)}   \eqno (5.1)  $$
as the function corresponding to the determinant of the minor of a matrix
$ \, g = \big( x_{ij} \big) \in GL_n(\bk) \, $  obtained by taking rows
$ i_1, \dots, i_r \, $  and columns  $ \, 1, \dots, r \, $.  Then  $ \,
d^I \in \cO\big(GL_n(\bk)\big) \, $  for all  $ I \, $,  i.e.~these
are regular functions on  $ GL_n(\bk) \, $.

\smallskip

   If  $ \, I_0 := (1,\dots,r) \, $,  then  $ d^{I_0} $  restricts to a map (with same name)
  $$  d^{I_0} \, : \; P \longrightarrow \bk^\times \; ,  \qquad
   M := \begin{pmatrix}
      A  &  B  \\
      0  &  D
        \end{pmatrix}  \, \mapsto \,
d^{I_0}\big(M\big) = \text{\it det}\big(A\big)  $$
which is a character of  $ P \, $.  One checks that the line bundle
$ \cL $  associated to  $ d^{I_0} $  is very ample, and it provides an
embedding of  $ G\big/P $  into a projective space, following the recipe
in  \S \ref{proj-case}.  Algebraically, this means that the graded algebra
$ \Ogp $  is realized as embedded into  $ \cO(G) $  as
  $$  \Ogp  \, = \,
{\textstyle \bigoplus_{n \geq 0}} \, \Ogp_n
\, = \, {\textstyle \bigoplus_{n \geq 0}} \,
H^0\big(G\big/P, \cL^{\otimes n}\big)  \quad .  $$
In particular, one can easily verify, for any set  $ I $  of  $ r $
rows, that
  $$  \quad  d^I(gp) \, = \, d^{I_0}(p) \, d^I(g)  \hskip43pt
\forall \;\; g \in GL_n(\bk) \, , \; p \in P  $$
i.e.~$ d^I $  is  $ d^{I_0} $--semi-invariant of degree 1.  In addition,
one proves that the  $ d^I $'s  form a  $ \bk $--basis  of the space
$ \Ogp_1 $ of semi-invariants of degree 1  (cf.~\cite{ls}).
                                                       \par
   On the other hand, the spaces  $ \, \Ogp_n = H^0 \big(G\big/P,
\cL^{\otimes n}\big) \, $  are in one-to-one correspondence   --- up to twisting by any integral power of  $ \text{\it det} \, $  (i.e., by any character of  $ GL_n(\Bbbk) \, $)  ---  with the irreducible representations of  $ GL_n(\bk) \, $.

\medskip

   We will now see that this picture extends to the quantum setup.

\medskip

   Let  $ \, \cO_q(M_n) \, $  be the unital
associative algebra over  $ \, \bk_q = \bk\big[q,q^{-1}\big] \, $
with generators  $ x_{ij} \, $  (for  $ \, 1 \leq i, j \leq n \, $)
and relations
  $$  \displaylines{
   x_{ij} \, x_{ik} \; = \; q \, x_{ik} \, x_{ij} \; ,  \qquad
 x_{ji} \, x_{ki} \; = \; q \, x_{ki} \, x_{ji}  \; \qquad
\forall \;\; j < k \; , \quad \forall \;\, i  \cr
   x_{ij} \, x_{kl} \; = \; x_{kl} \, x_{ij}  \; \qquad \qquad  \forall
\;\; i < k \, , \, j > l  \hbox{\ \ or \ }  i > k \, , \, j < l  \cr
   x_{ij} \, x_{kl} \, - \, x_{kl} \, x_{ij} \; = \;
\big( q - q^{-1} \big) \, x_{kj} \, x_{il}  \; \qquad
\forall \;\; i < k \, , \, j < l \; . }  $$
(the so-called ``Manin relations'').  This algebra bears also a structure
of  $ \bk_q $--bialgebra,  whose coproduct and counit are given by
  $$  \Delta (x_{ij})  \; = \;  {\textstyle \sum_{k=1}^n} \,
x_{ik} \otimes x_{kj} \;\; ,  \qquad  \epsilon(x_{ij})
\; = \;  \delta_{ij}  \qquad  \forall \;\;\; i \, , j \; .  $$
\indent
   Define the ``quantum determinant'' (of order  $ n \, $)  $ det_q $  as
  $$  det_q  \; := \;  {\textstyle \sum_{\sigma \in \,\mathcal{S}_n}}
(-q)^{\ell(\sigma)} \, x_{1,\,\sigma(1)} \cdots x_{n,\,\sigma(n)} \;
\in \; \cO_q(M_n)  \quad .  $$
One proves that  $ det_q $  belongs to the center of  $ \cO(M_n) \, $,
and it is group-like, i.e.~$ \; \Delta(det_q) = det_q \otimes det_q \; $
and  $ \; \epsilon(det_q) = 1 \; $.
                                                    \par
   More in general, for any  $ \, 1 \leq r \leq n \, $  and for any
choice of  $ r $--tuples  of  {\sl increasing\/}  indices  $ \, I =
(i_1,\dots,i_r) \, $  and  $ \, J = (j_1,\dots,j_r) \, $,  we define
the ``quantum determinant of the minor  $ (I,J) $'',  i.e.~of the
minor (of the matrix with entries the  $ x_{ij} $'s) whose sets of
rows and columns are  $ I $  and  $ J $,  namely
  $$  D^I_J  \; := \;
{\textstyle \sum_{\sigma \in \,\mathcal{S}_r}}
(-q)^{\ell(\sigma)} \, x_{i_1,\,j_{\sigma(1)}} \cdots
x_{i_r,\,i_{\sigma(r)}}  \quad .   \eqno (5.2)  $$
   \indent   These satisfy (cf.~\cite{ksc}, \S 9.2.2) the following
quantum analogue of well-known classical identities (e.g., the first
one is analogous to Binet theorem):
  $$  \Delta\big(D^I_J\big)  \, = \,
{\textstyle \sum_K} D^I_K \otimes D^K_J \quad ,  \qquad
\epsilon\big(D^I_J\big)  \, = \,  \delta_{I,\,J} \quad .
\eqno (5.3)  $$

\medskip

   Since  $ \text{\it det}_q $  is central in  $ \cO_q(M_n) \, $,  it is
a Ore element as well, and we can consider the enlarged algebra  $ \, \cO_q(GL_n) := \cO_q(M_n)\big[{\text{\it det}_q}^{\!-1}\big] \, $  obtained from  $ \cO_q(M_n) $  by formally inverting  $ \text{\it det}_q \, $.
Then   --- see  \cite{ksc}  again ---   the bialgebra structure of
$ \cO_q(M_n) $  uniquely extends to  $ \cO_q(GL_n) \, $;  even more,
the latter then is a  {\sl Hopf algebra\/}  indeed.  In fact, by construction
$ \cO_q(GL_n) $  {\sl is a quantum group},  namely a quantization of
$ GL_n(\bk) \, $,  in the sense of  Definition \ref{quant-G}.
                                                    \par
   We shall again denote by  $ x_{ij} $  the images in  $ \cO_q(GL_n) $
of the generators  $ x_{ij} $  of  $ \cO_q(M_n) \, $.  Similarly, we
shall again denote by  $ D^I_J $  the images in  $ \cO_q(GL_n) $  of
the ``quantum minors'' of  $ \cO_q(M_n) \, $:  then they again enjoy
(5.2) and (5.3).  Letting  $ \, J_0 := (1,\dots,r) =: I_0 \, $,
hereafter we shall set  $ \, D^I := D^I_{J_0} \, $.
                                                    \par
   The specialization (at  $ \, q = 1 $)  of any quantum minor
$ D^I_J $  is the corresponding classical minor  $ d^I_J $  (on
the same sets of rows and columns); in particular, every  $ D^I $
specializes to  $ d^I $   --- see (5.1) ---   and, among them, $
D^{I_0} $  to  $ d^{I_0} \, $.

\medskip

   We define the  \textit{quantum parabolic subgroup}  $ \cO_q(P) $  as
the quotient algebra
  $$  \cO_q(P)  \; := \;  \cO_q(GL_n) \bigg/ \Big( \big\{\, x_{ij}
\;\big|\, r\!+\!1 \leq i \leq n \, ; \, 1 \leq j \leq r \,\big\} \Big)
\quad .  $$
One can easily check that this  $ \cO_q(P) $  is in fact a  {\sl Hopf
algebra\/}  quotient.  Thus the natural projection map  $ \; \pi : \,
\cO_q(G) \relbar\joinrel\relbar\joinrel\twoheadrightarrow \cO_q(P) \; $
is a Hopf algebra epimorphism.  Therefore,  $ \cO_q(P) $  {\sl is a
quantum Poisson subgroup of}  $ \, \cO_q(G) = \cO_q\big(GL_n(\bk)\big)
\, $,  in the sense of  Definition \ref{qcoisosg},  whose specialization
is  $ \cO(P) \, $.

\medskip

   We are now in a position to appreciate the first important fact
--- in the present setting ---   about quantum minors:

\medskip

\begin{lemma}  \label{D_q-sect}
  The quantum minor  $ \, D^{I_0} $  is a quantum section of the line
bundle on  $ G\big/P $  given by  $ \, d^{I_0} $,  in the sense
of  Definition \ref{q-sect}.
\end{lemma}

\begin{proof}  Using the first identity in (5.3) one gets
  $$  \Delta_\pi\big(D^{I_0}\big)  \, = \,
\big((\text{\it id} \otimes \pi) \circ \Delta \big)\big(D^{I_0}\big)
\, = \,  (\text{\it id} \otimes \pi) \Big( {\textstyle \sum_K}
D^{I_0}_K \otimes D^K_{J_0} \Big)  \, = \,  {\textstyle \sum_K}
D^{I_0}_K \otimes \overline{D^K_{J_0}}  $$
and then from this
  $$  \Delta_\pi\big(D^{I_0}\big)  \; = \;
D^{I_0}_{I_0} \otimes \overline{D^{I_0}_{J_0}}
\; = \;  D^{I_0} \otimes \overline{D^{I_0}}  $$
because  $ \; \overline{D^K_{J_0}} := \pi\big(D^K_{J_0}\big) =
\delta_{K,\,I_0} \, \overline{D^{I_0}_{J_0}} \; $,  \, by definition of
$ \pi \, $,  and  $ \, D^{I_0}_{I_0} = D^{I_0}_{J_0} = D^{I_0} \; $.
Therefore  (Definition \ref{pre-q})  $ D^{I_0} $  is a pre-quantum
section; but  $ \cO_q(P) $  is a quantum subgroup, so  (Proposition
\ref{crit_q-sect})  $ D^{I_0} $  is a quantum section.
\end{proof}

\medskip

   Using  $ D^{I_0} \, $,  we can perform the construction of the
algebra  $ \Oqgp $  of  $ D^{I_0} $--semi-invariants  (or simply
{\sl semi-invariants\/}),  as in \S 3.  First we have

\medskip

\begin{lemma}   \label{D_semi-inv}
  The quantum minors  $ \, D^I \! $  are all semi-invariants
of degree 1, that is to say  $ \; D^I \in \Oqgp_1 \; $  for
every set of rows  $ \; I = (i_1,\dots,i_r) \; $.
\end{lemma}

\noindent
 {\it Proof.}  Arguing as in the proof of  Lemma \ref{D_q-sect},
we prove the claim by
  $$  \displaylines{
   \quad  \Delta_\pi\big(D^I\,\big)  \, = \,
\big((\text{\it id} \otimes \pi) \circ \Delta \big)
\big(D^I\,\big)  \, = \,
(\text{\it id} \otimes \pi) \Big( {\textstyle \sum_K}
D^I_K \otimes D^K_{J_0} \Big)  \, =   \hfill  \cr
   \hfill   = \,  {\textstyle \sum_K} D^I_K \otimes \overline{D^K_{J_0}}
\, = \, D^I_{I_0} \otimes \overline{D^{I_0}_{J_0}}  \, = \,
D^I \otimes \overline{D^{I_0}} \quad .  \qquad \square }  $$

\medskip

   Roughly speaking, the outcome of this last result is that the line
bundle on  $ G\big/P $  given by  $ d^{I_0} $  {\sl has enough ``quantum
sections'' to provide a ``quantum projective embedding''}.  To be precise,
the following holds:

\medskip

\begin{corollary}  The space  $ \, \Oqgp $  of all
$ \, D^{I_0} $--semi-invariants  is a quantization of
$ \, \Ogp $,  in the sense of  Definition \ref{q-prhspace}.
\end{corollary}

\begin{proof}   By construction, every quantum minor  $ D^I $  specializes
to the corresponding classical minor  $ d^I $.  By  \S \ref{class-case},
the latter form a basis of  $ \Ogp_1 \, $.
 \eject
\noindent
 This along with  Lemma
\ref{D_semi-inv}  proves that the natural embedding
  $$  \Oqgp_1 \Big/ (q-1) \, \Oqgp_1 \;
\lhook\joinrel\relbar\joinrel\relbar\joinrel\relbar\joinrel\longrightarrow
\; \Ogp_1  $$
is onto.  But then, as noticed in  Remark \ref{rem-qprhsp},  this is
enough to conclude.
\end{proof}

\medskip

   Actually, we can prove the following, much more precise result:

\medskip

\begin{proposition}  \label{semi-Grass}
  The algebra  $ \, \Oqgp $  is generated by the  $ D^I $'s.
\end{proposition}

\begin{proof}  By  Lemma \ref{D_semi-inv},  the  $ D^I $'s  belong to
$ \Oqgp \, $.  Therefore, we are only left to prove that, conversely,
every semi-invariant is contained in the  $ \bk $--subalgebra  of
$ \cO_q(G) $  generated by the  $ D^I $'s.
                                                  \par
   To this end, Theorems 1.2 and 1.3 in \cite{gl} give us immediately
the result if we take  $ \bk(q) $  as ground ring instead of  $ \, \bk_q
:= \bk\big[q,q^{-1}\big]\, $.  Then Lemmas 3.9, 3.10, 3.11 in  \cite{fi4}
give us our result.  We now see that in detail.

\smallskip

   We start by rewriting the Proposition 1.1 in  \cite{gl}
in our notation:

\smallskip

\textit{Let  $ \, A \stackrel{\phi} \lra B \stackrel{\psi} \lra C \, $
be a complex of  $ \bk_q $--modules,  such that  $ C $  is torsion free.
Suppose there are  $ \bk_q$--module  decompositions  $ \, A = \oplus_j
A_j \, $,  $ \, B = \oplus_j B_j \, $,  $ \, C = \oplus_j C_j \, $  such
that  $ B_j $  is finitely generated, and the maps  $ \phi $  and
$ \psi $  respect the decomposition, that is  $ \, \phi(A_j)
\subseteq B_j \, $  and  $ \, \psi(B_j) \subseteq C_j \, $.
Then if the sequence  $ \;  \overline{A} \lra \overline{B} \lra \overline{C} \; $  obtained by reduction modulo  $ (q-1) $  is exact, then so is also}
  $$  \bk(q) \otimes_{\bk_q} A \stackrel{\phi} \lra
\bk(q) \otimes_{\bk_q} B \stackrel{\psi} \lra
\bk(q) \otimes_{\bk_q} C  $$

\smallskip

   Let's apply this result to our special situation.
                                             \par
   The subalgebra  $ \, A := \bk_q\big[D^I\big] \, $  generated in
$ \cO_q({SL}_n) $ by quantum determinants is a  $ \bk_q $--graded
module (by the degree).  This fact is non trivial and depends on the
explicit form of this algebra in terms of generators and relations,
see  \cite{fi1}.  We have that an element  $ \, a \in \cO_q({SL}_n)
\, $  is in  $ \, A_j \, $  iff
  $$  {\textstyle \sum_{(a)}} \, a_{(1)} \otimes \overline{a}_{(2)}
\, = \,  a \otimes \overline{d}^{\,j} \; ,  \qquad  \text{where \ }
\Delta(a)  \, = \,  {\textstyle \sum_{(a)}} \,
a_{(1)} \otimes a_{(2)}  $$
where  $ \overline{x} $  denotes reduction of  $ x $  modulo
$ I_q(P) $  (see notation in section \ref{q-lin-bun}).
                                             \par
   So we can set up maps
  $$  A \stackrel{\phi} \lra \cO_q({SL}_n) \stackrel{\psi}
\lra \cO_q({SL}_n) \otimes \cO_q(P)  $$
where  $ \phi $  is the inclusion and  $ \psi $  is defined by
  $$  \psi(a)  \, = \,  {\textstyle \sum_{(a)}} \,
a_{(1)} \otimes \overline{a}_{(2)} - a \otimes \bd^{\,j}
\qquad  \forall \;\; a \in A_j \;\; .  $$
   \indent   One can check that all the hypothesis of the previous
result, for  $ \, B:= \cO_q({SL}_n) \, $  and  $ \, C := \cO_q({SL}_n)
\otimes \cO_q(P) \, $,  are satisfied, hence we obtain that  $ \; \bk(q) \otimes_{\bk} A  \, \cong \,  \text{\it Ker}\, \big( id \otimes \psi \big) \; $.
In other words, the semi-invariants coincide with the subalgebra
generated by the quantum determinants over the ring  $ \bk(q) \, $.

\medskip

   We now obtain the result over  $ \bk_q $  by Lemma 3.11
in  \cite{fi2},  namely

\medskip

   \hfill   \textit{If  $ \, wX \in \bk_q\big[D^I\big] \, $,  $ \, w \in \bk_q \, $,  $ \, w \neq 0 \, $,  then  $ \, X \in \bk_q\big[D^I\big] \; $.}   \hfill
\end{proof}

\smallskip

\begin{remark}
   Thus, using our own recipe, we have constructed the quantum homogeneous
space  $ \Oqgp \, $.  It is immediate to see that this is the same as the deformation of the algebra of the classical Grassmannian, along with its
classical Pl{\"u}cker embedding, as it is described in  \cite{fi1}  or in
\cite{tt}.
\end{remark}

\smallskip

   Finally, for the  $ \cO_q(G) $--comodule  structure of the space of
semi-invariants of degree 1, we have also the following analogue of a
classical result:

\medskip

\begin{proposition}
 $ \, \Oqgp_1 \cong \bigwedge_q \! \big(\bk_q^n\big) \, $  as left
$ \, \cO_q(G) $--comodules.
\end{proposition}

\begin{proof} This is a direct calculation. Let's sketch it.  By all
the previous analysis, we already know that  $ \Oqgp_1 $  has basis
the set of all the  $ D^I $'s,  and the left  $ \cO_q(G) $--coaction
on  $ \Oqgp_1 $  is given by
  $$  D^I  \, \mapsto \,
{\textstyle \sum_K} \, D^I_K \otimes D^K  \quad .  $$
   \indent   Now consider the coaction of  $ \cO_q(G) $  on
$ \bigwedge_q \! \big( \bk_q^n \big) \, $,  \, given by
  $$  \displaylines{
   \xi_{i_1} \cdots \xi_{i_r} \, \mapsto \,
{\textstyle \sum_{k_1, \dots, \, k_r}} \, g_{i_1\,k_1} \cdots g_{i_r\,k_r}
\otimes \xi_{k_1} \cdots \xi_{k_r}  \; =   \hfill  \cr
   \hfill   = \;  {\textstyle \sum_\sigma} (-q)^{\ell(\sigma)} g_{i_1\,k_1}
\cdots g_{i_r\,k_r} \otimes \xi_{k_1^o} \cdots \xi_{k_r^o}  \; = \;
{\textstyle \sum_K} D^I_K \otimes \xi_{k_1^o} \cdots \xi_{k_r^o} }  $$
where $\sigma$ is the permutation reordering  $ k_1, \dots, k_r $
and  $ \, K = \big( k_1^o, \dots, k_r^o \big) $  are the same indices,
but reordered. Now the result follows.
\end{proof}

\smallskip

\begin{remark}  \label{q-flag-var}
  Similar arguments can be used to prove that any quantum flag variety
is a ``quantum projective homogeneous space'' in the sense of  Defini\-tion
\ref{q-prhspace}  (for details about quantum flag varieties, we refer to
\cite{fi3}).
                                                        \par
   For the flag of type  $ (m_1, \dots, m_s) $,
the quantum section  $ d $  to start with is
 \vskip-5pt
  $$  d \, := \, D^{(m_1)} \cdots D^{(m_n)}  $$
 \vskip-1pt
\noindent
 where the  $ D^{(m_j)} $'s  are the principal quantum minors of size
$ m_j \, $.
                                                        \par
   The proofs of all results go over exactly as in the Grassmannian
case.
\end{remark}

   We now turn to the construction of the quantum big cell ring, that
will be crucial for the explicit construction of the QDP functor.

\medskip

\begin{definition}
  Let  $ \, I_0 = (1 \dots r) \, $,  $ \, D_0 := D^{I_0} \, $.  Define
  $$  \cO_q(G)\big[D_0^{-1}\big] \; := \;
\cO_q(G)[T] \Big/ \big( T \, D_0 - 1 \, , D_0 \, T - 1 \big)  $$
Moreover, we define the  {\sl big cell ring\/}  $ \, \cO_q^{\,\text{\it
loc}}\big( G\big/P \big) \, $  to be the  $ \Bbbk_q $--subalgebra  of
$ \, \cO_q(G)\big[D_0^{-1}\big] \, $  generated by the elements
  $$  {\ } \qquad  t_{ij} \, := \, {(-q)}^{r-j} \, \Dij \, D_0^{-1}
\eqno \forall \;\;\, i \, , j \; : \; 1 \leq j \leq r < i \leq n  $$
where  $ \, \Delta_{i{}j} := D^{1 \cdots \widehat{j} \cdots r \, i} \, $,  \, for all  $ i $,  $ j $  as above  (see \cite{fi2} for more details).
\end{definition}

\medskip

   As in the commutative setting, we have the following result:

\medskip

\begin{proposition} \label{classical}
  $ \; \displaystyle{ \cO_q^{\,\text{\it loc}}\big( G\big/P \big)
\, \cong \, \cO_q\big( G\big/P \big)\big[D_0^{-1}\big]_{proj} } \; $,
\; where the right-hand side is the degree-zero component of
$ \; \cO_q\big( G\big/P \big)[T] \Big/ \big(T D_0 - 1 \, ,
D_0 \, T - 1 \big) \, $.
\end{proposition}

\begin{proof}  In the classical setting, the analogous result is proved by
this argument: one uses the so-called ``straightening relations'' to get rid of the extra minors (see, for example,  \cite{dep},  \S 2).  Here the argument works essentially the same, using the  {\sl quantum straightening\/
{\rm (or} Pl{\"u}cker\/{\rm )}  relations\/}  (see  \cite{fi1},  \S 4,
\cite{tt}, formula (3.2)(c) and Note I, Note II).
\end{proof}

\medskip

\begin{remark}  As before, we have that
  $$  \cO_q^{\,\text{\it loc}}\big( G \big/ P \big)
\,\;{\textstyle \bigcap}\;\, (q-1) \, \cO_q^{\,\text{\it loc}}(G) \,\;
= \;\, (q-1) \, \cO_q^{\,\text{\it loc}}\big( G \big/ P \big)  $$
This can be easily deduced from Remark  \ref{spec_semi-inv-n},
taking into account Proposition  \ref{classical}.  As a consequence,
the map
  $$  \cO_q^{\,\text{\it loc}}\big( G \big/ P \big) \Big/ (q-1)
\, \cO_q^{\,\text{\it loc}}\big( G \big/ P \big) \;
\relbar\joinrel\relbar\joinrel\relbar\joinrel\longrightarrow
\; \cO_q^{\,\text{\it loc}}(G) \Big/ (q-1) \,
\cO_q^{\,\text{\it loc}}(G)  $$
is  {\sl injective},  so that the specialization map
  $$  \pi_{G/P}^{\,\text{\it loc}} \colon \,
\cO_q^{\,\text{\it loc}}\big( G \big/ P \big) \;
\relbar\joinrel\relbar\joinrel\relbar\joinrel\twoheadrightarrow
\; \cO_q^{\,\text{\it loc}}\big( G \big/ P \big) \Big/
(q-1) \, \cO_q^{\,\text{\it loc}}\big( G \big/ P \big)  $$
coincides with the restriction of the specialization map
  $$  \pi_G^{\,\text{\it loc}} \; \colon \; \cO_q^{\,\text{\it loc}}(G)
\; \relbar\joinrel\relbar\joinrel\relbar\joinrel\twoheadrightarrow
\; \cO_q^{\,\text{\it loc}}(G) \Big/ (q-1) \,
\cO_q^{\,\text{\it loc}}(G)  \quad .  $$
\end{remark}

\medskip

   The following proposition gives a description of the algebra
$ \cO_q^{\,\text{\it loc}}\big( G \big/ P \big) \, $:

\medskip

\begin{proposition} \label{bigcell}
  The big cell ring is isomorphic to a matrix algebra
  $$  \begin{array}{cccl}
   \qquad  \cO_q^{\,\text{\it loc}}\big( G\big/P \big)  &  \lra
&  \cO_q\big(M_{(n-r) \times r}\big)  &  \\
   \qquad  t_{ij}  &  \mapsto  &  x_{ij}  &  \qquad  \forall \;\;
1 \leq j \leq r < i \leq n
      \end{array}  $$
i.e.~the generators  $ t_{ij} $'s  satisfy the Manin relations.
\end{proposition}

\begin{proof}  See \cite{fi2},  Proposition 1.9.
\end{proof}

\medskip

\begin{remark}
   The Grassmannian  $ \, GL_n \big/ P \, $  can also be realized as a similar quotient of  $ SL_n $  by a suitable parabolic  $ P' $  (corresponding to  $ P $,  say).  Then one can also perform all related quantum constructions   --- the previous and the later ones ---   using
$ SL_n $  instead of  $ GL_n \, $,  and modifying each step as needed.  To begin with, one considers
  $$  \cO_q(SL_n)  \, := \,  \cO_q(GL_n) \Big/ \big( \text{\it det}_q - 1 \big)  \, \cong \, \cO_q(M_n) \Big/ \big( \text{\it det}_q - 1 \big)  $$
--- where  $ \, \big( \text{\it det}_q - 1 \big) \, $  is the (two-sided) ideal generated by  $ \, \text{\it det}_q - 1 \, $  ---   which is again a Hopf algebra, for the quotient structure from either  $ \cO_q(GL_n) $  or
$ \cO_q(M_n) \, $.  This is a quantization of  $ SL_n(\bk) \, $,  in the sense of  Definition \ref{quant-G},  for which we can consider again quantum minors and a corresponding  $ \cO_q(P) $  as before.  Then all this can be used to give an alternative definition of  $ \, \cO_q \big( G\big/P \big) =
\cO_q \big( SL_n \big/ P' \big) \, $  and of all was considered above.  Similarly, all constructions and results of section  \ref{qdp-grass}  hereafter can be carried on using  $ \cO_q(SL_n) $   --- and its related gadgets ---   instead of  $ \cO_q(GL_n) \, $.
                                                \par
   Finally, similar considerations hold as well for the quantum flag varieties mentioned in  Remark \ref{q-flag-var}.
\end{remark}

\medskip

\subsection{QDP for quantum Grassmannians}
\label{qdp-grass}

   We now turn to the quantum duality principle applied explicitly to
the quantum homogeneous spaces constructed above. We start with Grassmannians.
                                                          \par
   Let us first explicitly describe the Poisson structure of the algebraic group  $ GL_n \, $.  Starting from  $ \cO_q({GL}_n) \, $,  as usual the classical algebra
$ \cO({GL}_n) $  inherits from the former a Poisson bracket, which makes it into a Poisson Hopf algebra, so that  $ {GL}_n $  becomes a Poisson
group (see  Remark \ref{pois-bra}{\it (2)\/}).  We want to describe now this Poisson bracket.  Recall that
  $$  \cO({GL}_n)  \; = \;  \Bbbk \big[ {\{\, \bar{x}_{ij} \,\}}_{i, j
= 1, \dots, n} \,\big]\big[\text{\it det}^{-1}\big]  \; = \;
\Bbbk \big[ {\{\, \bar{x}_{ij} \,\}}_{i, j = 1, \dots, n} \,\big][t]
\Big/ \big( t \, \text{\it det} - 1 \big)  $$
where  $ \; \text{\it det} := \text{\sl det}\, \big( \bar{x}_{i,j} \big)_{i, j = 1, \dots, n} \; $  is the usual determinant.  Setting  $ \, \bar{x}
= p(x) \, $  for  $ \; p : \cO_q(GL_n) \lra \cO(GL_n) \; $,  \,
the Poisson structure is given (as usual) by
  $$  \big\{ \bar{a} \, , \bar{b} \,\big\} \; := \; {(q-1)}^{-1}
\, ( a \, b - b \, a )\Big|_{q=1}   \eqno \forall \;\; \bar{a} \, ,
\bar{b} \in \cO({GL}_n) \;\; .  $$
In terms of generators, we have
  $$  \displaylines{
   \big\{ \bar{x}_{ij} \, , \bar{x}_{ik} \big\} \, = \, \bar{x}_{ij}
\, \bar{x}_{ik}  \quad  \forall \;\; j < k \; ,   \phantom{\Big|}
 \quad  \big\{ \bar{x}_{ij} \, , \bar{x}_{\ell k} \big\} \, =
\, 0  \qquad  \forall \;\; i < \ell \, , k < j  \cr
   \big\{ \bar{x}_{ij} \, , \bar{x}_{\ell j} \big\} \, = \,
\bar{x}_{ij} \, \bar{x}_{\ell j}  \quad  \forall \;\; i < \ell \, ,
\phantom{\Big|}  \qquad  \big\{ \bar{x}_{ij} \, , \bar{x}_{\ell k} \big\}
\, = \, 2 \, \bar{x}_{ij} \, \bar{x}_{\ell k}  \quad  \forall \;\;
i < \ell \, , j < k  \cr
   \big\{ \text{\it det}^{-1} , \bar{x}_{ij} \big\} \, = \,0 \; ,
\phantom{\Big|}  \quad  \big\{ \text{\it det} \, , \bar{x}_{ij} \big\}
\,= \, 0  \qquad \forall \;\; i, j = 1, \dots, n \, .  \cr }  $$
   \indent   As  $ {GL}_n $  is a Poisson Lie group, its Lie algebra
$ \fgl_n $  has a Lie bialgebra structure (see  \cite{cp},  {pg.~\!}24).
To describe it, let us denote with  $ \m_{ij} $  the elementary matrices, which form a basis of  $ \fgl_n \, $.  Define  ($ \, \forall \; i = 1,
\dots, n-1 \, $,  $ \, j = 1, \dots, n \, $)
  $$  e_i := \m_{i,i+1} \; ,  \quad  g_j := \m_{j,j} \; ,  \quad
f_i := \m_{i+1,i} \; ,  \quad  h_i := g_{i} - g_{i+1}  $$
Then  $ \; \big\{\, e_i \, , \, f_i \, , \, g_j \;\big|\; i = 1, \dots, n-1,
\, j = 1, \dots, n \,\big\} \; $  is a set of Lie algebra generators of
$ \fgl_n \, $,  and a Lie cobracket is defined on  $ \fgl_n $  by
  $$  \delta(e_i) \, = \, h_i \otimes e_i - e_i \otimes h_i \, ,  \quad
\delta(g_j) \, = \, 0 \, ,  \quad  \delta(f_i) \, = \, h_i \otimes f_i
- f_i \otimes h_i \;   \eqno \forall \;\, i , j .  $$
This cobracket makes  $ \fgl_n $  itself into a  {\sl Lie bialgebra\/}:
this is the so-called  {\sl standard\/}  Lie bialgebra structure on
$ \fgl_n \, $.  It follows immediately that  $ \, U(\fgl_n) \, $  is
a co-Poisson Hopf algebra, whose co-Poisson bracket is the (unique)
extension of the Lie cobracket of  $ \fgl_n $  while the Hopf
structure is the standard one.

\medskip

   Similar constructions hold for the group  $ {SL}_n \, $.  One simply
drops the generator  $ d^{-1} \, $,  and imposes the relation  $ \, d \!
= \! 1 \, $,  in the description of  $ \cO(SL_n) \, $,  and replaces the
$ g_s $'s  with the  $ h_i $'s  ($ \, i = 1, \dots, n \, $)  when
describing  $ \, \fsl_n \, $.

\medskip

   Since  $ \fgl_n $ is a Lie bialgebra, its dual space  $ \fgl_n^{\,*} $
admits a Lie bialgebra structure, dual to the one of  $ \fgl_n \, $.  Let
$ \, \big\{\, \e_{ij} := \m_{ij}^{\,*} \;\big|\; i, j = 1, \dots, n
\,\big\} \, $  be the basis of  $ \, \fgl_n^{\,*} \, $  dual to the
basis of elementary matrices for  $ \fgl_n \, $.  As a Lie algebra,
$ \fgl_n^{\,*} $  can be realized as the subset of  $ \, \fgl_n
\oplus \fgl_n \, $  of all pairs
  $$  \left( \! \begin{pmatrix}
   \! -m_{11} \!  &  \! 0 \!  &  \! \cdots \!  & \!   0 \!   \\
   \! m_{21} \!  &  \! -m_{22} \!  &  \! \cdots & \!  0 \!   \\
   \! \vdots \!  &  \! \vdots \!  &  \! \vdots \!  &  \! \vdots \!   \\
   \! m_{n-1,1} \!  &  \! m_{n-1,2} \!  &  \! \cdots \!  &  \! 0 \!   \\
   \! m_{n,1} \!  &  \! m_{n,2} \!  &  \! \cdots \!  &  \! -m_{n,n}
      \end{pmatrix} , \,
      \begin{pmatrix}
   m_{11} \!  &  \! m_{12} \!  &  \! \cdots \!  &  \! m_{1,n-1} \!
&  \! m_{1,n} \!   \\
   \! 0 \!  &  \! m_{22} \!  &  \! \cdots \!  &  \! m_{2,n-1} \!  &
\! m_{2,n} \!   \\
   \! \vdots \!  &  \! \vdots \!  &  \! \vdots \!  &  \! \vdots \!  &
\! \vdots \!   \\
   \! 0 \!  &  \! 0 \!  &  \! \cdots \!  &  \! m_{n-1,n-1} \!  &
\! m_{n-1,n} \!   \\
   \! 0 \!  &  \! 0 \!  &  \! \cdots \!  &  \! 0 \!  &  \! m_{n,n} \!
      \end{pmatrix} \! \right)  $$
with its natural structure of Lie subalgebra of  $ \, \fgl_n \oplus
\fgl_n \, $.  In fact, the elements  $ \, \e_{ij} \, $  correspond to
elements in  $ \, \fgl_n \oplus \fgl_n \, $  in the following way:
  $$  \e_{ij} \cong \big( \m_{ij} \, , 0 \big)  \hskip5pt  \forall \;
i \! > \! j \, ,  \hskip7pt  \e_{ij} \cong \big(\! -\m_{ij} \, , +\m_{ij}
\big)  \hskip5pt  \forall \; i \! = \! j \, ,  \hskip7pt  \e_{ij} \cong
\big( 0 \, , \m_{ij} \big)  \hskip5pt  \forall \; i < j \, .  $$
Then the Lie bracket of  $ \, \fgl_n^{\,*} \, $  is given by
  $$  \begin{array}{clc}
   \big[ \e_{i,j} \, , \, \e_{h,k} \big]  \hskip-5pt  &  = \, \delta_{j,h}
\, \e_{i,k} - \delta_{k,i} \, \e_{h,j} \;\; ,   &   \forall  \;\; i \!
\leq \! j \, , \, h \! \leq \! k  \;\;\,  \text{and}  \;\; \forall \;\;
i \! > \! j \, , \, h \! > \! k  \\
                                 \\
   \big[ \e_{i,j} \, , \, \e_{h,k} \big]  \hskip-5pt  &  = \, \delta_{k,i}
\, \e_{h,j} - \delta_{j,h} \, \e_{i,k} \;\; ,  &   \forall \;\; i \! = \!
j \, , \, h \! > \! k  \;\;\,  \text{and}  \;\; \forall \;\; i \! > \! j
\, , \, h \! = \! k  \\
                                 \\
   \big[ \e_{i,j} \, , \, \e_{h,k} \big]  \hskip-5pt  &  = \, 0 \; ,  &
\forall\;\; i \! < \! j \, , \, h \! > \! k  \;\;\,  \text{and}  \;\;
\forall \;\; i \! > \! j \, , \, h \! < \! k
      \end{array}  $$
   \indent   Note that the elements  ($ \, 1 \leq i \leq n \! - \! 1 \, $,
$ \, 1 \leq j \leq n \, $)
  $$  \e_i \, = \, e_i^{\,*} \, = \, \e_{i,i+1} \;\; ,  \qquad
\f_i \, = \, f_i^{\,*} \, = \, \e_{i+1,i} \;\; ,  \qquad
\g_j \, = \, g_j^{\,*} \, = \, \e_{jj}  $$
are Lie algebra generators of  $ \fgl_n^{\,*} \, $.
In terms of them, the Lie bracket reads
  $$  \big[ \e_i \, , \f_j \big] \, = \, 0 \; ,  \qquad  \big[ \g_i \, ,
\e_j \big] \, = \, \delta_{ij} \, \E_i \; ,  \qquad  \big[ \g_i \, , \f_j
\big] \, = \, \delta_{ij} \, \f_j   \eqno \forall \;\; i, j \; .  $$

\medskip

   On the other hand, the Lie cobracket structure of  $ \fgl_n^{\,*} $
is given by
  $$  \delta\big(\e_{i,j}\big) \, = \,
{\textstyle \sum\limits_{k=1}^n} \, \e_{i,k} \wedge \e_{k,j}
\eqno \forall \;\; i, j = 1, \dots, n  $$
where  $ \, x \wedge y := x \otimes y - y \otimes x \; $.

\medskip

   Finally, all these formul{\ae}  also provide a presentation of
$ U\big(\fgl_n^{\,*}\big) $  as a co-Poisson Hopf algebra.

\medskip

   A similar description holds for  $ \, \fsl_n^{\,*} \! = \fgl_n^{\,*}
\Big/ Z\big(\fgl_n^{\,*}\big) \, $,  where  $ \, Z \big( \fgl_n^{\,*}
\big) \, $  is the centre of  $ \fgl_n^{\,*} \, $,  generated by  $ \,
\mathfrak{l}_n := \g_1 \! + \cdots + \g_n \, $.  The construction is
immediate by looking at  the embedding  $ \, \fsl_n \hookrightarrow
\fgl_n \, $.

\bigskip

We now turn to the construction of the QDP functor.

\medskip

\begin{observation} \label{qdp}
$ \cO_q(G)^\vee $  (see Definition \ref{OqGvee}  and \S  \ref{class-case})  is generated, as a unital subalgebra of  $ \, \cO_q(G) \otimes_{\Bbbk_q} \Bbbk(q) \, $,  by the elements
  $$  \D_- \, := \, {(q-1)}^{-1} \, \big( D_q^{-1} - 1 \big) \; ,
\quad \hskip5pt  \chi_{ij} \, := \, {(q-1)}^{-1} \, \big( x_{ij}
- \delta_{ij} \big)   \eqno \forall \; i, j = 1, \dots, n  $$
where the  $ x_{ij} $'s  are the generators of  $ \cO_q(G) \, $.  As
$ \, x_{ij} = \delta_{ij} + (q-1) \, \chi_{ij} \in \cO_q(G)^\vee \, $,
we have an obvious embedding of  $ \cO_q(G) $  into  $ \cO_q(G)^\vee
\, $.
\end{observation}

   Following  Definition \ref{cOqghvee},  we define
  $$  \cO_q\big( G \big/ P \big)^\vee  \, := \,  \big\langle \,
{(q-1)}^{-1} \, J_{G/P}^{\,\text{\it loc}} \, \big\rangle \, = \,
{\textstyle \sum\limits_{n=0}^\infty} \, {(q-1)}^{-n} \, {\big(
J_{G/P}^{\,\text{\it loc}} \big)}^n  \;\;\; .  $$

   We can provide a concrete description of
$ {\cO_q\big( G\big/P \big)}^\vee \, $:

\medskip

\begin{proposition}  We have
 \vskip-5pt
  $$  {\cO_q\big( G\big/P \big)}^\vee  \; = \;\,
\Bbbk_q \Big\langle {\{\, \mu_{i{}j} \,\}}_{i = r+1,
\dots, n}^{j = 1, \dots, r} \Big\rangle \bigg/ I_M  $$
 \vskip-1pt
\noindent
 where  $ \; \mu_{i{}j} := {(q-1)}^{-1} \, t_{i{}j} \; $  (for all
$ i $  and  $ j \, $),  $ \, I_M \, $  is the ideal of the Manin
relations among the  $ \mu_{i{}j} $'s, and  $ \; t_{ij} :=
{(-q)}^{r-j} \, \D_{ij} \, D_0^{-1} \, $  (for all  $ i $
and  $ j $).
\end{proposition}

\begin{proof}  Trivial from definitions and Proposition \ref{bigcell}.
\end{proof}

\medskip

   We want to see explicitly what is  $ \, {\cO_q\big( G\big/P \big)}^\vee
{\Big|}_{q=1} $  inside  $ \, U\big({\fgl_n}^{\!*}\big) \, $.
In other words, we want to understand what is the space that
$ {\cO_q\big( G\big/P \big)}^\vee $  is quantizing.  We check now
by direct inspection that this is  $ \, U\big(\fp^\perp\big) \, $,  as
already prescribed by  Theorem \ref{quant1}.

\medskip

\begin{proposition}  \label{quantumspec}
  $$  {\cO_q\big( G \big/ P \big)}^\vee{\Big|}_{q=1}
\, = \; U\big(\fp^\perp\big)  $$
as a subalgebra of  $ \; {\cO_q(G)}^\vee{\Big|}_{q=1} \! =
U\big({\fgl_n}^{\!*}\big) \, $,  where  $ \, \fp^\perp \, $
is the orthogonal subspace to  $ \, \fp := \text{\it Lie}\,(P)
\, $  inside  $ \, {\fgl_n}^{\!*} \, $.
\end{proposition}

\begin{proof}  Thanks to the previous discussion, it is enough to
show that
  $$  \pi_G^\vee\Big({\cO_q\big( G \big/ P \big)}^\vee\Big) \, =
\, U\big(\fp^\perp\big) \, \subseteq \, U \big( {\fgl_n}^{\!*}
\big) \, = \, {\cO_q(G)}^\vee{\Big|}_{q=1}  \quad .  $$
To do this, we describe the isomorphism  $ \, {\cO_q(G)}^\vee
{\Big|}_{q=1} \! \cong U\big({\fgl_n}^{\!*}\big) \, $  (cf.~\cite{ga2}).
According to  Remark \ref{qdp},  the algebra  $ {\cO_q(G)}^\vee $  is generated by the elements
  $$  \D_- \, := \, {(q-1)}^{-1} \, \big( D_q^{-1} - 1 \big) \; ,
\quad \hskip5pt  \chi_{ij} \, := \, {(q-1)}^{-1} \, \big( x_{ij}
- \delta_{ij} \big)   \eqno \forall \; i, j = 1, \dots, n  $$
inside  $ \, \cO_q(G) \otimes_{\Bbbk_q} \Bbbk(q) \, $.
In terms of these  generators, the isomorphism reads
  $$  \displaylines{
   {\cO_q(G)}^\vee{\Big|}_{q=1} \,
\relbar\joinrel\relbar\joinrel\relbar\joinrel\longrightarrow
\; U\big({\fgl_n}^{\!*}\big)  \cr
   \overline{{\D}_-} \mapsto -(\e_{1,1} + \cdots + \e_{n,n})
\; ,  \qquad  \overline{\chi_{i,j}} \mapsto \e_{i,j} \qquad
\forall \;\; i \, , j \; .  \cr }  $$
where we used notation  $ \; \overline{X} := X \mod (q-1) \,
{\cO_q(G)}^\vee \; $.  Indeed, from  $ \, \overline{\chi_{i,j}}
\mapsto \e_{i,j} \, $  and  $ \, {(q-1)}^{-1} \, \big( D_q - 1
\big) \in {\cO_q(G)}^\vee \, $,  one gets  $ \, \overline{D_q}
\mapsto 1 \, $  and  $ \, \overline{{(q-1)}^{-1} \, \big( D_q
- 1 \big)} \mapsto \e_{1,1} + \cdots + \e_{n,n} \, $.  Moreover,
the relation  $ \, D_q \, D_q^{-1} = 1 \, $  in  $ \cO_q(G) $
implies  $ \; D_q \, \D_- = - {(q-1)}^{-1} \, \big( D_q - 1 \big)
\; $  in  $ {\cO_q(G)}^\vee \, $,  whence clearly  $ \; \overline{\D_-}
\mapsto -(\e_{1,1} + \cdots + \e_{n,n}) \; $  as claimed.

\medskip

   In other words, the specialization  $ \, p_G^\vee \, \colon
{\cO_q(G)}^\vee \relbar\joinrel\relbar\joinrel\twoheadrightarrow
U\big({\fgl_n}^{\!*}\big) \, $  is given by
  $$  p_G^\vee \big( \D_- \big) \, = \, -(\e_{1,1} + \cdots + \e_{n,n})
\; ,  \qquad  p_G^\vee \big( \chi_{i,j} \big) \, = \, \e_{i,j}  \qquad
\forall \;\; i \, , j \; .  $$
   \indent   If we look at  $ \widehat{{\cO_q(G)}^\vee} $,  things are
even simpler.  Since
  $$  D_q \, \in \, \Big( 1 + (q-1) \, {\cO_q(G)}^\vee \Big) \subset
\Big( 1 + (q-1) \, \widehat{{\cO_q(G)}^\vee} \Big) \; ,  $$
then  $ \, D_q^{-1} \in \Big( 1 + (q-1) \, \widehat{{\cO_q(G)}^\vee}
\Big) \, $, and the generator  $ \, \D_- \, $  can be dropped.
The specialization map  $ \widehat{p_{G/P}^\vee} $  of course is
still described by  formul\ae{}  as above.

\medskip

   Now let us compute  $ \, p_{G/P}^\vee \Big( {\cO_q\big( G\big/P
\big)}^\vee \Big) = \widehat{p_G^\vee} \Big( {\cO_q\big( G\big/P
\big)}^\vee \Big) \, $.  Recall that  $ {\cO_q\big( G\big/P
\big)}^\vee $  is generated by the  $ \mu_{i{}j} $'s,  with
  $$  \mu_{i{}j} \; := \; {(q-1)}^{-1} \, t_{i{}j} \; = \;
{(q-1)}^{-1} \, {(-q)}^{r-j} \, \Dij \, D_0^{-1}  $$
for  $ \, i = r+1, \dots, n \, $,  and  $ \, j = 1, \dots, r \, $;
thus we must compute  $ \, \widehat{p_G^\vee} \big( \mu_{i{}j}
\big) \, $.
                                           \par
   By definition, for every  $ \, i \not= j \, $  the element
$ \, x_{i{}j} = (q-1) \, \chi_{i{}j} \, $  is mapped to 0 by
$ \widehat{p_G^\vee} \, $.  Instead, for each  $ \ell $  the
element  $ \, x_{\ell\,\ell} = 1 + (q-1) \, \chi_{\ell\,\ell} \, $
is mapped to 1 (by  $ \widehat{p_G^\vee} $  again).  But then,
expanding the  $ q $--determinants  one easily finds that
  $$  \displaylines{
   \widehat{p_G^\vee} \, \Big( {(q-1)}^{-1} \, \Dij \Big) \; = \;
\Big( {(q-1)}^{-1} \, {\textstyle \sum_{\sigma \in \, \mathcal{S}_r}}
\, {(-q)}^{\ell(\sigma)} \, x_{1 \, \sigma(1)} \cdots x_{r \, \sigma(r)}
\Big)  \, =   \hfill  \cr
   \hfill   = \; \widehat{p_G^\vee} \, \Big( {(q-1)}^{-1} \, {\textstyle
\sum\limits_{\sigma \in \, \mathcal{S}_r}} {(-q)}^{\ell(\sigma)} \,
\big( \delta_{1 \, \sigma(1)} + (q-1) \chi_{1 \, \sigma(1)}) \cdots
\big( \delta_{1 \, \sigma(r)} + (q-1) \, \chi_{1 \, \sigma(r)})
\Big)  \cr }  $$
The only term in  $ (q-1) $  in the expansion of  $ \Dij $  comes
from the product
  $$  \big( 1 + (q \! - \! 1) \chi_{1 \, 1}) \cdots \big( 1 +
(q \! - \! 1) \chi_{r \,r} \big) \, (q \! - \! 1) \chi_{i \, j}
\equiv  (q \! - \! 1) \chi_{i \, j}  \! \mod {(q \! - \! 1)}^2
\cO\big(G \! \big/ \! P\big)  $$
Therefore, from the previous analysis we get
  $$  \displaylines{
   \widehat{p_G^\vee} \, \Big( {(q-1)}^{-1} \, \Dij \Big)
\; = \; \widehat{p_G^\vee} \, \big( \chi_{i,j} \big) \;
= \; \e_{i,j}  \cr
   \widehat{p_G^\vee} \, \big( D_0 \big) \; = \;
\widehat{p_G^\vee} \big( 1 \big) \; = \; 1 \;\; ,  \qquad
\widehat{p_G^\vee} \, \big( D_0^{-1} \big) \; = \;
\widehat{p_G^\vee} \big(1 \big) \; = \; 1  \cr }  $$
hence we conclude that  $ \; \widehat{p_G^\vee} \big( \mu_{i{}j} \big)
= {(-1)}^{r-j} \, \e_{i,j} \; $,  \, for all  $ \, 1 \leq j \leq r < i
\leq n \, $.

\medskip

   The outcome is that  $ \; p_{G/P}^\vee \Big( {\cO_q\big( G\big/P
\big)}^\vee \Big) \, = \, U(\mathfrak{h}) \; $,  \; where
  $$  \mathfrak{h} \, := \, \text{\it Span}\,\big( \big\{\, \e_{i,j}
\,\big|\, r+1 \leq i \leq n \, , \; 1 \leq j \leq r \,\big\} \big)
\;\; .  $$
On the other hand, from the very definitions and our description
of  $ \, {\fgl_n}^{\!*} \, $  one easily finds that
$ \; \mathfrak{h} = \fp^\perp \; $,  \; for  $ \; \fp :=
\text{\it Lie}\,(P) \; $.  The claim follows.
\end{proof}

\begin{proposition}  $ \, \widehat{{\cO_q\big( G \big/ P \big)}^\vee} \, $
is a left coideal of  $ \; \widehat{{\cO_q(G)}^\vee} \, $.
\end{proposition}

\begin{proof} This is the same as Proposition \ref{quant2}.
\end{proof}

Hence for the quantum Grassmannian we have proved directly the
following result:

\begin{theorem}
  $ \, \widehat{{\cO_q\big( G \big/ P \big)}^\vee} \, $
is a quantum homogeneous  $ \, G^* $--space,  which is
an infinitesimal quantization of the homogeneous  $ \,
G^* $--space  $ \fp^\perp \, $.
\end{theorem}

\subsection{Quantum generalized flag varieties for
simple gro\-ups as quantum projective homogeneous spaces}

  We now turn to a more general example of quantum projective space:
namely the quantization of a generalized flag variety for any simple
group, following  \cite{ko}.  As before, we begin with a
brief description of the classical setting.

\medskip

   Let  $ G $  be a connected, simply connected, complex simple Lie group,  and
   let  $ \fg $  its Lie algebra.  Let  $ S $  be a subset of simple roots of  $ \fg $  and
   let  $ \, \la = \sum_{\al_i \not\in S} \omega_i \, $  be a weight of  $ \fg \, $,  where the
   $ \omega_i $'s  are fundamental weights.

\medskip

   Let  $ V(\la) $  be the highest weight representation of  $ \fg $  (and of  $ G \, $)  associated with the weight  $ \la \, $,  and let  $ v_\la $  be a non-zero highest weight vector of  $ V(\la) \, $.  We have the following morphism of algebraic varieties:
  $$  \begin{array}{ccc}
   G  &  \lra  &  \bP(V(\la))  \\
   g  & \mapsto  &  g \cdot v_\la
      \end{array}  $$
This induces a projective embedding of the flag variety  $ \, G\big/P_S \, $
into the projective space  $ \bP\big(V(\la)\big) \, $,  where  $ \, P_S :=
\text{\it Stab}_G(v_\la) \, $  is the parabolic subgroup associated to the
set  $ S \, $.  The graded algebra of regular functions on  $ G\big/P_S $
relative to this embedding is given by
  $$  \cO\big(G\big/P_S\big)  \; = \;
{\textstyle \bigoplus}_{n \in \N} \, {V(n\la)}^*   \eqno (5.4)  $$
where the grading is given by  $ \, {\cO\big(G\big/P_S\big)}_n := {V(n\la)}^*
\, $  and the multiplication is via the Cartan multiplication (see  \cite{fi4} for more details).

\medskip

   We are now going to identify  $ \cO\big(G\big/P_S\big) $  with a graded
subalgebra of  $ \cO(G) \, $.  Indeed, the algebra  $ \cO(G) $  is in (Hopf)
duality with  $ U(\fg) $  and it can be thought of as the linear span inside
$ {U(\fg)}^* $  of the functionals  $ \; c^\la_{f,v} : U(\fg) \lra \C \, $
(the so-called ``matrix coefficients'') given by
  $$  \qquad  c^\la_{f,v}(u) := f(u.v)  \;\qquad \forall \;\; u \in U(\fg)
\, ,  \quad f \in V(\la)^* \, ,  \quad v \in V(\la) \; .  $$
Then  $ \cO\big(G\big/P_S\big) $  can be realized, equivalently, as the subalgebra
of  $ \cO(G) $  generated by the  $ c^\la_{f,v_\la} $'s,  for all  $ \, f \in
{V(\la)}^* \, $.

\medskip

   This point of view carries over without changes to the quantum setting:
it is considered in  \cite{lr},  \cite{so}  and by several others; a quick
review of this construction can be found in \cite{ko}.  The key point is that
every finite dimensional representation  $ V(\mu) $  of  $ U(\fg) $  of highest
weight  $ \mu $  has a (standard) quantization, which is a representation for
$ U_q(\fg) \, $,  call it  $ V_q(\mu) \, $.  In particular, one can again define
matrix coefficient functionals  $ \, c^\mu_{f,v} \, $   --- for every  $ \, f \in
{V_q(\mu)}^* \, $  and  $ \, v \in V_q(\mu) \, $  ---   for all dominant weights
$ \mu $  of  $ \fg \, $.  Their  $ \bk_q $--span  inside  $ {U_q(\fg)}^* $  is, by definition, the algebra  $ \cO_q(G) \, $,  which is a quantization of  $ \cO(G)
\, $.  In fact, the latter follows because  $ V_q(\mu) \, $,  respectively
$ {V_q(\mu)}^* \, $,  specializes to  $ V(\mu) \, $,  respectively to
$ {V(\mu)}^* \, $,  at  $ \, q = 1 \, $,  hence ``quantum'' and ``classical'' matrix coefficients bijectively correspond to each other   --- via  $ \, c^\mu_{f,v} \mapsto c^\mu_{f,v} \, $,  say ---   under specialization at
$ \, q = 1 \, $.

\smallskip

   For later use, we point out how matrix coefficient behave under the coproduct.  For any dominant weight  $ \mu $  of  $ \fg \, $,  let  $ \, \{v_1, \dots, v_r\} \, $  be a  $ \bk_q $--basis  of  $ V_q(\mu) $   --- with  $ \, r = \text{\it dim}\,\big(V_q(\la)\big) \, $  ---   and let  $ \, \{f_1, \dots, f_r\} \, $  be the dual basis of  $ {V_q(\mu)}^* \, $.  Then every matrix coefficient
$ c^{\,\mu}_{f,v} $  has coproduct given by
  $$  \Delta \big( c^{\,\mu}_{f,v} \big)  \; = \;
{\textstyle \sum}_{i=1}^r \, c^{\,\mu}_{f,v_i} \otimes c^{\,\mu}_{f_i,v}
\eqno (5.5)  $$
(just because  $ U_q(\fg) $  is acting on  $ V_q(\mu) $  via matrices!).

\medskip

   From the quantization  $ \cO_q(G) $  the group  $ G $  inherits a Poisson group structure   --- cf.~Remark \ref{pois-bra}{\it (2)\/} ---   for which  $ P_S $  {\sl is a Poisson subgroup}.  On the other hand, let  $ \, I_q(P_S) \, $  be the
 two-sided ideal of  $ \cO_q(G) $  generated by the set of matrix coefficients
$ \, \Big\{\, c^{\,\mu}_{f,v} \;\Big|\;\, \forall \; n \in \N \, , \, \mu \not= n \la \, \vee \, v \not\in \big( \bk_q \! \setminus \! \{0\} \big) \! \cdot \! v_{n\la} \,\Big\} \, $.  Then, using (5.5), one easily shows that  $ I_q(P_S) $  {\sl is a Hopf ideal of}  $ \, \cO_q(G) \, $;  therefore the quotient  $ \bk_q $--module  and the canonical map
  $$  \cO_q(P_S)  \; := \;  \cO_q(G) \Big/ I_q(P_S)  \quad ,  \qquad
\pi \colon \cO_q(G) \relbar\joinrel\twoheadrightarrow \cO_q(P_S)
\eqno (5.6)  $$
are respectively a quotient Hopf algebra and a Hopf algebra epimorphism.  Indeed, this provides  {\sl a quantization of}  $ \, P_S \, $,  as a Poisson subgroup of
$ G \, $,  namely  $ \cO_q(P_S) $  {\sl is a quantum subgroup of}  $ \, G $
in the sense of  Definition \ref{qcoisosg}.

\bigskip

   In sight of the above construction of  $ \cO_q(G) $  and of the classical description of  $ \cO\big(G\big/P_S\big) $  in (5.4), we define
  $$  {\cO_q\big(G\big/P_S\big)}'  \, := \,
{\textstyle \bigoplus}_{n \in \N} \, {V_q(n\la)}^*  \, = \;
\bk_q\text{\sl --span  of\ }  \; {\big\{ c^{\,n\la}_{f,v_{n\la}}
\big\}}_{f \in {V_q(n\la)}^*, \, n \in \N}   \eqno (5.7)  $$
%
%
where  $ v_\la $  is a non-zero element of weight  $ \la $  in  $ V_q(\la) \, $.
Then, as the quantum matrix coefficients specialize to the classical ones, comparing (5.7) with (5.4) we see at once that
  $$  {\cO_q\big(G\big/P_S\big)}' \Big/ (q-1) \, {\cO_q\big(G\big/P_S\big)}'
\,\; \cong \;\,  \cO\big(G\big/P_S\big)   \eqno (5.8)  $$
so that  $ {\cO_q\big(G\big/P_S\big)}' $  is a quantization, as a  $ \bk_q $--module,  of  $ \cO\big(G\big/P_S\big) \, $.

\medskip

   We are now going to show that this  $ {\cO_q\big(G\big/P_S\big)}' $  is in fact a quantum homogeneous space in the sense of  Definition \ref{q-prhspace},  in particular it can be realized as the space of semi-invariants inside  $ \cO_q(G) $  with respect to a suitable quantum section.  Indeed, we shall find  $ \, {\cO_q\big(G\big/P_S\big)}' = \cO_q\big(G\big/P_S\big) \, $  where the latter
is the space of all semi-invariants (for a suitable quantum section) as in  Definition \ref{def_semi-inv}  and  Definition \ref{q-prhspace}.

\smallskip

   First of all, let  $ \lambda $  be the dominant weight fixed above, and let
$ v_\la $  be the (up to a scalar factor) uniquely determined non-zero element of weight  $ \la $  in  $ V_q(\la) \, $.  Fix a  $ \bk_q $--basis  $ \, \{v_1, \dots, v_r\} \, $  of  $ V_q(\la) $  with  $ \, v_1 = v_\la \, $,  and let  $ \, \{f_1, \dots, f_r\} \, $  be the dual basis of  $ {V_q(\mu)}^* \, $,  for which we set
$ \, f_\la := f_1 \, $.

\medskip

\begin{proposition}  \label{c_q-sect}
  The element  $ \, c^{\,\la}_{f_\la,v_\la} $  is a quantum section in
$ \, \cO_q(G) \, $.
\end{proposition}

\begin{proof}  Thanks to  Proposition \ref{crit_q-sect}{\it (b)\/}  we only need to show that  $ c^{\,\la}_{f_\la,v_\la} $  is a pre-quantum section, with respect to the setup of (5.6), i.e.
  $$  \Delta_\pi \big( c^{\,\la}_{f_\la,v_\la} \big)  \, = \,
c^{\,\la}_{f_\la,v_\la} \otimes \pi \big( c^{\,\la}_{f_\la,v_\la} \big)
\quad .  $$
But this follows at once from the identity (5.5), applied to the bases chosen above, once we notice in addition that  $ \, c^{\,\la}_{f_i,v_\la} \in I_q(P_S)
\, $  for all  $ \, i \not= 1 \, $.
\end{proof}

\smallskip

\begin{proposition}  \label{1_semi-inv}
  The space of  $ c^{\,\la}_{f_\la,v_\la} $--semi-invariants  of
degree 1 inside  $ \cO_q(G) \, $,  that is  $ \; \cO_q\big(G\big/P_S\big)_1 \,
:= \, \Big\{\, f \in \cO_q(G) \;\Big|\; \Delta_\pi(f) = f \otimes
\pi\big(c^{\,\la}_{f_\la,v_\la}\big) \Big\} \, $,  \; is just the  $ \,
\bk_q $--span  of  $ \, \big\{\, c^{\,\la}_{f,v_\la} \;\big|\;
f \in {V_q(\la)}^* \,\big\} \, $.  In other words,
  $$  \cO_q\big(G\big/P_S\big)_1  \,\; = \;\,  \bk_q\text{\sl --span  of \ }
\big\{\, c^{\,\la}_{f,v_\la} \;\big|\; f \in {V_q(\la)}^* \,\big\}  \quad .  $$
\end{proposition}

\begin{proof}  This again is immediate as before.  Consider any  $ \bk_q $--linear  combination of several  $ c^{\,\mu}_{f,v} $'s  which is semi-invariant of degree 1, with respect to the quantum section  $ \, d := c^{\,\la}_{f_\la,v_\la} \, $.  We can assume these  $ c^{\,\mu}_{f,v} $'s  to be linearly independent over
$ \bk_q \, $,  and so the semi-invariance of their linear combination as a whole also implies the semi-invariance of each of the  $ c^{\,\mu}_{f,v} $'s  on its own.
                                               \par
   Now, assume that a single matrix coefficient  $ c^{\,\mu}_{f,v} $  is
semi-invariant of degree 1 (with respect to  $ c^{\,\la}_{f_\la,v_\la} \, $).  Then (5.5) implies at once that  $ \, \mu = \la \, $.  Moreover, choosing bases
$ \, \big\{\, v_1 = v_\la \, , \, v_2 \, , \, \dots \, , \, v_r \,\big\} \, $  and  $ \, \big\{\, f_1 = f_\la \, , \, f_2 \, , \, \dots \, , \, f_r \,\big\} \, $  as before  Pro\-position \ref{c_q-sect},  the identity (5.5) also gives
  $$  c^{\,\la}_{f,v} \otimes \pi\big(c^{\,\la}_{f_\la,v_\la}\big)
\; = \;  \Delta_\pi \big( c^{\,\la}_{f,v} \big)  \; = \;
c^{\,\la}_{f,v_\la} \otimes \pi\big(c^{\,\la}_{f_\la,v}\big)
\, + \,  {\textstyle \sum}_{i=2}^r \, c^{\,\la}_{f,v_i} \otimes \pi\big(c^{\,\la}_{f_i,v}\big)  $$
This forces  $ \, c^{\,\la}_{f_i,v} \in \text{\it Ker}\,(\pi) = I_q(P_S) \, $  for all  $ \, i > 1 \, $,  so that  $ \, v \in \bk_q \! \cdot \! v_\la \, $,  say
$ \, v = \kappa \, v_\la \, $  for some  $ \, \kappa \in \bk \, $,  whence eventually
  $$  c^{\,\la}_{f,v}  \, = \,  c^{\,\la}_{f,\,\kappa\,v}  \, = \,
\kappa \, c^{\,\la}_{f,v}  \, \in \,  \bk_q\text{\sl --span  of \ }
\big\{\, c^{\,\la}_{f,v_\la} \;\big|\; f \in {V_q(\la)}^* \,\big\}  \quad .  $$
   \indent   This proves that  $ \cO_q\big(G\big/P_S\big)_1 $  is indeed contained in the prescribed  $ \bk_q $--span.  The converse is clear, just reversing the previous argument.
\end{proof}

\medskip

\begin{proposition}  \label{n_semi-inv}
    The space of  $ c^\la_{f_\la,v_\la} $--semi-invariants
of degree n inside  $ \cO_q(G) \, $,  that is  $ \; \cO_q\big(G\big/P_S\big)_n
\, := \, \Big\{\, f \in \cO_q(G) \;\Big|\; \Delta_\pi(f) = f \otimes
\pi\Big(\!\big(c^\la_{f_\la,v_\la}\big)^n\Big) \Big\} \, $,  is just
the  $ \, \bk_q $--span  of  $ \, \big\{\, c^{\,n\la}_{f,v_{n\la}} \;\big|\;
f \in {V_q(n\la)}^* \,\big\} \, $.  In other words,
  $$  \cO_q\big(G\big/P_S\big)_n  \,\; = \;\,  \bk_q\text{\sl --span  of \ }
\, \big\{\, c^{\,n\la}_{f,v_{n\la}} \;\big|\; f \in {V_q(n\la)}^* \,\big\}
\;\; .  $$
\end{proposition}

\begin{proof}  This follows from an argument which closely mimics the one
used in the proof of  Proposition \ref{1_semi-inv}.  One takes into account,
in addition, the following two remarks:
                                                \par
   {\it (a)} \,  the vector  $ v_\la^{\,\otimes n} $  has weight  $ \, n \la
\, $  inside  $ {V_q(\la)}^{\otimes n} \, $;  thus it can be canonically identified with a (non-zero) highest weight vector, say  $ \, v_{n\la} \, $,
in  $ V_q(n\la) \, $,  hence it can be chosen as  $ \, v_1 := v_{n\la} \, $,
the first element of a suitable  $ \bk_q $--basis  of  $ V_q(n\la) $  to be used in that argument;
                                                \par
   {\it (b)} \,  with notation as above, the  $ n $--th  power function  $ {\big( c^{\,\la}_{f_\la,v_\la} \big)}^n $  inside  $ \cO_q(G) $  is nothing but a matrix coefficient again, namely  $ \, {\big( c^{\,\la}_{f_\la,v_\la} \big)}^n = c^{\,n\la}_{f_{n\la},v_{n\la}} \, $.

\smallskip

   These two remarks, drafted into an argument totally similar to the one
used for  Proposition \ref{1_semi-inv},  eventually yield the claim.
\end{proof}

\smallskip

   We are now ready for the main result of this subsection:

\smallskip

\begin{theorem}
  Let  $ \cO_q\big(G\big/P_S\big) $  be defined as in  Definition
\ref{def_semi-inv},  with respect to the quantum section  $ \, d :=
c^\la_{f_e,v_\la} \in \cO_q(G) \, $.  Then  $ \cO_q\big(G\big/P_S\big) $
is a quantum projective homogeneous space, namely, it is a quantization of
$ \, \cO\big(G\big/P_S\big) \, $,  in the sense of  Definition \ref{q-prhspace}.  \end{theorem}

\begin{proof}  This follows at once by putting together the previous results, i.e.~Pro\-positions \ref{c_q-sect},  \ref{1_semi-inv}  and  \ref{n_semi-inv},  and the specialization formula (5.8).
\end{proof}

\medskip

\begin{remark}  {\ }

\smallskip

   {\it (1)} \,  With some extra work, one can also show that  $ \cO_q \big( G\big/P_S \big) $  is generated   --- as a graded algebra ---   in degree 1, i.e.~by
$ \, {\cO_q\big(G\big/P_S\big)}_1 \; $.

\smallskip

   {\it (2)} \,  Our last remark is that in this setup of quantum generalized flag varieties one can
   also apply the QDP, following the general recipe of \S 4.
                                              \par
   Indeed, in \cite{ko}, \S 3.4, it is noticed that the quantum section  $ \,
d := c^\la_{f_e,v_\la} \, $  {\sl is a Ore element in}  $ \cO_q(G) \, $.
Therefore, as pointed out in  \S \ref{qdp-philosophy},  we can define the
localizations
  $$  \cO_q\big(G\big/P_s\big)\Big[{\big(c^\la_{f_e,v_\la}\big)}^{-1}\Big]
\, \subseteq \,  \cO_q(G)\Big[{\big(c^\la_{f_e,v_\la}\big)}^{-1}\Big]
\quad .  $$
and we can then apply the QDP   --- according to \S 4 ---   to this setting.
\end{remark}

\medskip
\subsection{The coisotropic case}
 One could argue whether the generality we are working with is necessary or not.
 In this section we will describe how to find families of coisotropic parabolic subgroups
 inside semisimple algebraic groups with the standard
 multiplicative Poisson structure.

Such families give rise to smooth Poisson bivectors on the
projective quotients which cannot be obtained as quotient by
Poisson parabolic subgroups. It would be therefore interesting to
investigate whether in such examples it is possible to find and
quantize a graded Poisson algebra associated to them as described
in section 2.

 We shall start with a low-dimensional example and describe in a very explicit manner the objects we are interested in and then we will generalize such results.
                                       \par
  Let us consider the group  $ {SL}_2(\mathbb C) $  with the standard
Poisson algebraic structure, normalized as follows: if  $ \, a, b,
c, d \, $ are matrix coefficients in positions  $ \,
\begin{pmatrix} a & b
\\  c & d  \end{pmatrix} \, $,  we let
%
%
 $ \; \{ a \, , b \} \, = \, a \, b \; $,  $ \; \{ a \, , c \} \, = \, a \, c \; $,  $ \; \{ b \, , d \} \, = \, b \, d \; $,  $ \; \{ c \, , d \} \, = \, c \, d \; $,  $ \; \{ b \, , c \big\} \, = \, 0 \; $,  $ \; \{ a \, , d \big\} \, = \, 2 \, a \, d \; $  (this is the opposite normalization to that in  \cite{kos}).  We take the standard parabolic subgroup of upper triangular matrices
  $$  P \; = \; \left\{ \left(\begin{array}{cc}a&b\\ 0&d\end{array}\right) \,\bigg|\; a, b, d \in \C \,\right\} $$
This is a Poisson subgroup in  $ {SL}_2(\mathbb C) \, $;  thus,
the quotient  $ \; \bP^1 \C \simeq {SL}_2(\C) \big/ P \; $  is
endowed with the (homogenous) quotient smooth Poisson bivector  $
\pi_0 \, $.
                                                      \par
   Let us now consider the following element
  $$  g_{\varepsilon} \; := \, \left(\!\! \begin{array}{cc} \sqrt\varepsilon & \sqrt{1-\varepsilon}  \\
-\sqrt{1-\varepsilon}  &  \sqrt{\varepsilon}  \end{array} \right)
\;\; ,  \qquad  \varepsilon\in [0,1]  $$
and let  $ \; P_\varepsilon  \, := \, g_\varepsilon \, P \,
g_\varepsilon^{-1} \; $.  Then $ P_\varepsilon $  is defined
inside the group  $ {SL}_2(\C) $ by the equation:
  $$  \sqrt{\varepsilon(1-\varepsilon)} \, (a-d\,)  \; = \;  (\varepsilon-1) \, b \, + \, \varepsilon \, c  \quad .  $$
The infinitesimal generators of its Lie algebra are
  $$  \displaylines{
   H_\varepsilon \, := \, g_\varepsilon \, H \, g_\varepsilon^{-1} = \, (2 \, \varepsilon - 1) \, H \, - \, 2 \, \sqrt{\varepsilon (1 - \varepsilon)} \, \big( X^+ + X^- \big) \;\; ,  \cr
   X_\varepsilon \, := \, g_\varepsilon \, X^+ \, g_\varepsilon^{-1} = \, \sqrt{\varepsilon( 1 - \varepsilon)} \, H \, + \, \varepsilon \, X^+ - (1 - \varepsilon) \, X^- \;\; .  }  $$
It is then easily verified, through the infinitesimal criterion of
Proposition \ref{coiso-inf},  that  $ P_\varepsilon $  is
coisotropic, because
  $$  \delta(H_\varepsilon) \; = \; H_\varepsilon \wedge H  \quad ,  \qquad
\delta(X_\varepsilon) \; = \; X_\varepsilon \wedge H \quad .  $$
This means that on  $ \bP^1\C $  there is an induced Poisson
bivector $ \pi_\varepsilon $  as quotient  $ \, {SL}_2(\C) \big/
P_\varepsilon \, $.  That this Poisson bivector is different form
$ \pi_0 $  follows considering the image of the diagonal subgroup
of  $ {SL}_2(\C) \, $,  which induces a single 0-dimensional
Poisson leaf with respect to  $ \pi_0 $  and an  $ {\mathbb S}^1
$--family  of 0-dimensional leaves with respect to  $
\pi_\varepsilon \, $.
                                                        \par
   This phenomenon, as said above, is not specific only of  $ \bP^1 \C $  but can be, for example, shown to take place for all complex Grassmannians.
                                                        \par
   Let  $ G $  be a complex semisimple algebraic group and let  $ K $  be its
real compact form.  Up to a factor  $ \imath \, $,  which is
inessential in what follows, the standard Poisson structure on  $ G $  can be identified with the one which is automatically defined on it when it is identified with the Drinfeld's double of  $ K \, $.  Let now  $ H $  be a coisotropic subgroup of  $ K $  and let
us consider the subgroup
 $ H K^* $  of  $ G $  (here  $ \, K^* = A N \, $  is the Manin dual of  $ K $  inside  $ G \, $).  Then  $ \, P = H K^* \, $  is parabolic in  $ G \, $,  $ \, H = P \cap K \, $  and  $ \, K \big/ H \simeq G \big/ P \, $  as smooth manifolds.  It can be shown quite easily that the coisotropy of  $ K $  implies the coisotropy of  $ P \, $,  and furthermore, via Theorem 4.1 in  \cite{cish},  that  $ K \big/ H $  and $ G \big/ P $  are also Poisson diffeomorphic.  Thus in order to check whether  $ P $  is coisotropic it is enough to check whether  $ P \cap K $  is coisotropic w.r.~to the standard Poisson structure on the compact real group  $ K \, $.  There we can rely on results in  \cite{cish},  where a 1-parameter family of coisotropic subgroups $ \, H_\varepsilon \subseteq SU(n) \, $  was given.  Such subgroups induce a
$ 1 $--parameter  family of homogeneous Poisson quotients on
complex Grassmannians.


\bigskip


\begin{thebibliography}{99}

\bibitem{bcgst} F.~Bonechi, N.~Ciccoli, R.~Giachetti, E.~Sorace and
M.~Tarlini, {\it The coisotropic subgroup structure of quantum
$ SL(2,{\mathbb R}) $},  Journ.~Geom.~Phys. {\bf 37}, 190--200 (2001).

\bibitem{bo} A.~Borel, {\it Linear algebraic groups},  Graduate Texts
in Mathematics  {\bf 126},  Springer-Verlag, New York-Heidelberg, 1991.

\bibitem{cp} V. ~Chari, A. ~Pressley, {\it A guide to quantum groups},
Cambridge Press, 1994.

\bibitem{cic} N.~Ciccoli, \textit{Quantization of Co-isotropic
Subgroups},  Lett.~Math.~Phys.~{\bf 42}  (1997), 123--138.

\bibitem{cg1} N.~Ciccoli, F.~Gavarini,  {\it A quantum duality principle
for coisotropic subgroups and Poisson quotients},  Adv.~Math.~{\bf 199}
(2006), 104--135.

\bibitem{cg2} N.~Ciccoli, F.~Gavarini,  {\it A global quantum duality
principle for subgroups and homogeneous spaces},  work in progress.

\bibitem{cish} N.~Ciccoli, A. J.--L.~Sheu, {\it Covariant Poisson structures on complex Grassmannians}, Comm.~Anal.~Geom., {\bf 14}  (2006), 443--474.

\bibitem{dep}  C. De Concini, D. Eisenbud, C. Procesi,  {\it Young
Diagrams and Determinantal Varieties},  Invent.~Mathematic{\ae}
{\bf 56}  (1980), 129--165.

\bibitem{fi1} R.~Fioresi,  {\it Quantization of the Grassmannian manifold},
J.~Algebra  {\bf 214}  (1999), 418--447.

\bibitem{fi2} R.~Fioresi,  {\it A deformation of the
big cell inside the Grassmannian manifold $G(r,n)$},
Rev.~Math.~Phy.~{\bf 11}  (1999), 25--40.

\bibitem{fi3} R.~Fioresi {\it Quantum deformation of the flag variety}
Comm.~Algebra, Vol. 27, n. 11 (1999).

%
%

\bibitem{fi4} R.~Fioresi,  {\it Quantum coinvariant theory for the
quantum special linear group and Quantum Schubert varieties},
J.~Algebra  {\bf 242}  (2001), 433--446.

\bibitem{ga1} F.~Gavarini,  {\it The quantum duality principle},
Ann.~Inst.~Fourier  {\bf 52}  (2002), 809--834.

\bibitem{ga2} F.~Gavarini, {\it The global quantum duality principle},
Journal f\"ur die reine und angewandte Mathematik  {\bf 612}  (2007),
17--33   ---   see also the expanded version  {\tt http://arxiv.org/abs/math/0303019}  (2003).

\bibitem{ga3} F.~Gavarini, {\it The global quantum duality
principle: theory, examples, and applications},
{\tt http://arxiv.org/abs/math/0303019}  (2003).

\bibitem{gl} K.~Goodearl, T.~Lenagan, {\it Quantized coinvariants at
trascendental $q$}, in ``Hopf algebras in non commutative
geometry and physics'', 155--165,
Lecture Notes in Pure and Appl.~Math.
{\bf 239}, Dekker, 2005.

\bibitem{gh} P.~Griffiths, J.~Harris, {\it Principles of
algebraic geometry}, Wiley Interscience, 1994.

\bibitem{ha} R.~Hartshorne, {\it Algebraic Geometry},
Graduate Texts in Mathematics  {\bf 52},  Springer-Verlag,
New York-Heidelberg, 1977.

\bibitem{ko} S. Kolb, {\it The AS-Cohen-Macaulay property for
quantum flag manifolds of minuscule weight},  J.~Algebra (2008),
to appear   --- see also  {\tt http://arxiv.org/abs/0707.1389}  (2007).

\bibitem{kos} L.~Korogodsky, Y.~Soibelman,  {\it Algebras of
functions on quantum groups, Part I},  Math.~Surv.~and Monographs
{\bf 56},  A.M.S., Providence (RI), 1998.

\bibitem{ksc} A.~Klimyk, K.~Schm{\"u}dgen,  {\it Quantum groups
and their representations},  Texts and Monographs in Physics,
Springer-Verlag, Berlin, 1997.

\bibitem{lr} V. ~Lakshmibai, N. ~Reshetekhin
{\it Quantum flag and Schubert schemes},
(Amherst, MA, 1990), 145--181, Contemp. Math., 134, Amer. Math. Soc.,
Providence, RI, (1992).

\bibitem{ls} V.~Lakshmibai, C.~S.~Seshadri, C.~Musili,
{\it Geometry of G/P. IV. Standard monomial theory for
classical types},  Proc.~Indian Acad.~Sci.~Sect.~A
Math.~Sci.~{\bf 88}  (1979), 279--362.

\bibitem{lu} J.~H.~Lu,  {\it Multiplicative and affine Poisson
structures on Lie groups},  Ph.D.~thesis University of California,
Berkeley, 1990   ---   see also  {\tt http://hkumath.hku.hk/\~{}jhlu/thesis.tex}.

\bibitem{ma1} Y.~Manin, {\it Topics in non commutative geometry},
M.~B.~Porter Lectures, Princeton University Press, Princeton, 1991.

\bibitem{mo}  S.~Montgomery,  {\it Hopf Algebras and Their Actions
on Rings},  CBMS Regional Conference Series in Mathematics  {\bf 82},
AMS, Providence, 1993.

\bibitem{so} Y. S. ~Soibelman {\it On the quantum flag manifold},
Func.~Ana.~Appl.~{\bf 25}, 225-227, (1992).

\bibitem{tt} E.~Taft, J.~Towber, {\it Quantum deformation
of flag schemes and Grassmann schemes. I. A $q$-deformation
of the shape-algebra for ${\rm GL}(n)$},  J.~Algebra  {\bf 142}
(1991), 1--36.

\bibitem{za} S.~Zakrzewski,  {\it Poisson homogeneous spaces},
in: J.~Lukierski, Z.~Popowicz, J.~Sobczyk (eds.),  {\sl Quantum
groups (Karpacz, 1994)},  629--639, PWN, Warsaw, 1995.

\end{thebibliography}
\end{document}